\documentclass[11pt]{article}
\UseRawInputEncoding
\usepackage{currfile}
\usepackage{color}
\usepackage{mathrsfs}
\usepackage{caption}
\usepackage{subfigure}
\usepackage{setspace}
\usepackage{epstopdf}

\usepackage{amscd}
\usepackage{amssymb,amsfonts,amsmath,psfrag,amscd,cite}
\usepackage{amsmath}
  \usepackage{paralist}
  \usepackage{graphics}
  \usepackage{epsfig}
  \usepackage[colorlinks=true]{hyperref}
  \usepackage[all]{xy}
\newtheorem{theo}{Theorem}

\newtheorem{defi}[theo]{Definition}
\newtheorem{lemm}[theo]{Lemma}
\newtheorem{rema}[theo]{Remark}

\makeatletter \@addtoreset{equation}{section}
\@addtoreset{theo}{section} \makeatother

\begin{document}


\date{}
\title{Geometric Reductions, Dynamics and Controls \\ for Hamiltonian System with Symmetry}
\author{Hong Wang \\
School of Mathematical Sciences and LPMC,\\
Nankai University, Tianjin 300071, P.R.China\\
E-mail: hongwang@nankai.edu.cn}
\date{\emph{ Dedicated to 80th anniversary of the birth of Professor Jerrold E. Marsden}\\
August 26, 2022 } \maketitle

{\bf Abstract:} This is a survey article,
from the viewpoint of the completeness of
the Marsden-Weinstein reduction,
to introduce briefly some recent developments
of the symmetric reductions and Hamilton-Jacobi theory of
the regular controlled Hamiltonian systems,
the nonholonomic controlled Hamiltonian systems
and the controlled magnetic Hamiltonian systems.
These research reveal the deeply internal
relationships of the geometrical structures of phase spaces,
the nonholonomic constraint, the dynamical
vector fields and the controls of these systems.\\




\tableofcontents

\section{Marsden-Weinstein Reduction on a Cotangent Bundle}

It is well-known that the reduction theory for the mechanical system with symmetry
is an important subject and that it is widely studied in the theory of
mathematics and mechanics, as well as applications; see
Abraham and Marsden \cite{abma78}, Arnold \cite{ar89},
Libermann and Marle \cite{lima87}, Marsden \cite{ma92},
Marsden \textit{et al.} \cite{mamiorpera07, mamora90},
Marsden and Perlmutter \cite{mape00},
Marsden and Ratiu \cite{mara99}, Marsden and
Weinstein \cite{mawe74}, Meyer \cite{me73},
Nijmeijer and Van der Schaft \cite {nivds90}
and Ortega and Ratiu \cite{orra04}
and so on., for more details and development.\\

In mechanics, the phase space of a
Hamiltonian system is often the cotangent bundle $T^*Q$ of a
configuration manifold $Q$, and the reduction theory on the cotangent
bundle of a configuration manifold is a very important special case
of general symplectic reduction theory. In the following we first
give the Marsden-Weinstein reduction for a Hamiltonian system
with symmetry on the cotangent bundle of a smooth
configuration manifold with a canonical symplectic form;
see Abraham and Marsden \cite{abma78} and
Marsden and Weinstein \cite{mawe74}.\\

Let $Q$ be a smooth manifold and $TQ$ the tangent bundle, $T^* Q$
is the cotangent bundle with a canonical symplectic form $\omega$.
Assume that $\Phi: G\times Q \rightarrow Q$ is a left smooth action
of a Lie group $G$ on the manifold $Q$. The cotangent lift is the
action of $G$ on $T^\ast Q$, $\Phi^{T^\ast}:G\times T^\ast
Q\rightarrow T^\ast Q$ given by $g\cdot
\alpha_q=(T\Phi_{g^{-1}})^\ast\cdot \alpha_q,\;\forall\;\alpha_q\in
T^\ast_qQ,\; q\in Q$. The cotangent lift of any proper (resp. free)
$G$-action is proper (resp. free). Moreover, assume that the cotangent lift action is
symplectic with respect to the canonical symplectic form $\omega$,
and that it has an $\operatorname{Ad}^\ast$-equivariant momentum map
$\mathbf{J}:T^\ast Q\to \mathfrak{g}^\ast$ given by
$<\mathbf{J}(\alpha_q),\xi>=\alpha_q(\xi_Q(q)), $ where $\xi\in
\mathfrak{g}$, $\xi_Q(q)$ is the value of the infinitesimal
generator $\xi_Q$ of the $G$-action at $q\in Q$, $<,>:
\mathfrak{g}^\ast \times \mathfrak{g}\rightarrow \mathbb{R}$ is the
duality pairing between the dual $\mathfrak{g}^\ast $ and
$\mathfrak{g}$. Assume that $\mu\in \mathfrak{g}^\ast $ is a regular
value of the momentum map $\mathbf{J}$, and $G_\mu=\{g\in
G|\operatorname{Ad}_g^\ast \mu=\mu \}$ is the isotropy subgroup of
the coadjoint $G$-action at the point $\mu$. From the
Marsden-Weinstein reduction, we know that the reduced space
$((T^*Q)_\mu, \omega_\mu)$ is a symplectic manifold.\\

In the following we give further a precise analysis of the
geometrical structure of the symplectic reduced space of $T^* Q$.
From Marsden and Perlmutter \cite{mape00} and Marsden \textit{et al.}
\cite{mamiorpera07}, we know that the classification of symplectic
reduced space of the cotangent bundle $T^* Q$ as follows: (1) If
$\mu=0$, the symplectic reduced space of cotangent bundle $T^\ast Q$
at $\mu=0$ is given by $((T^\ast Q)_\mu, \omega_\mu)= (T^\ast(Q/G),
\omega_0)$, where $\omega_0$ is the canonical symplectic form of the
reduced cotangent bundle $T^\ast (Q/G)$. Thus, the symplectic
reduced space $((T^\ast Q)_\mu, \omega_\mu)$ at $\mu=0$ is a
symplectic vector bundle; (2) If $\mu\neq0$, and $G$ is Abelian,
then $G_\mu=G$, in this case the Marsden-Weinstein symplectic reduced
space $((T^*Q)_\mu, \omega_\mu)$ is symplectically diffeomorphic to
symplectic vector bundle $(T^\ast (Q/G), \omega_0-B_\mu)$, where
$B_\mu$ is a magnetic term; (3) If $\mu\neq0$, and $G$ is not
Abelian and $G_\mu\neq G$, in this case theMarsden-Weinstein symplectic
reduced space $((T^*Q)_\mu, \omega_\mu)$ is symplectically
diffeomorphic to a symplectic fiber bundle over $T^\ast (Q/G_\mu)$
with fiber to be the coadjoint orbit $\mathcal{O}_\mu$, see the
cotangent bundle reduction theorem---bundle version, also see
Marsden and Perlmutter \cite{mape00} and Marsden \textit{et al.}
\cite{mamiorpera07}.\\

Thus, from the above discussion, we know that
the Marsden-Weinstein symplectic reduced space of a Hamiltonian system
defined on the cotangent bundle of a configuration manifold may not be a
cotangent bundle. Therefore, the symplectic reduced system of a
Hamiltonian system with symmetry defined on the cotangent bundle
$T^*Q$ may not be a Hamiltonian system on a cotangent bundle; that
is, the set of Hamiltonian systems with symmetries on the cotangent
bundle is not complete under the Marsden-Weinstein reduction.

\section{Hamilton-Jacobi Equations}

The Hamilton-Jacobi theory is an important research subject
in mathematics and analytical mechanics.
On the one hand, it provides a characterization
of the generating functions of certain time-dependent canonical
transformations, such that for a given Hamiltonian system in such a form,
its solutions are extremely easy to find by reduction to the
equilibrium. On the other hand,
it is possible in many cases that the Hamilton-Jacobi equation provides an
immediate way to integrate the equation of the motion of a system, even
when the problem of the Hamiltonian system itself has not been or cannot
be solved completely; see Abraham \textit{et al.} \cite{abma78, abmara88},
Arnold \cite{ar89} and Marsden and Ratiu \cite{mara99}.
In addition, the Hamilton-Jacobi equation is
also fundamental in the study of the quantum-classical relationship
in quantization, and it plays an important role
in the study of stochastic dynamical systems; see
Woodhouse \cite{wo92}, Ge and Marsden \cite{gema88},
and L\'{a}zaro-Cam\'{i} and Ortega \cite{laor09}.
For these reasons, the equation is described as a useful tool in the study of
Hamiltonian system theory, which has been extensively developed in
recent years and become one of the most active subjects in the
study of modern applied mathematics and analytical mechanics.\\

The Hamilton-Jacobi theory, from the
variational point of view, was originally developed by Jacobi in 1866,
and it states that the integral of the Lagrangian of a mechanical system along the
solution of its Euler-Lagrange equation satisfies the
Hamilton-Jacobi equation. The classical description of this problem
from the generating function and the geometrical point of view was
given by Abraham and Marsden in \cite{abma78} as follows:
letting $Q$ be a smooth manifold and $TQ$
the tangent bundle, $T^* Q$ is the cotangent bundle with a canonical
symplectic form $\omega$, and the projection $\pi_Q: T^* Q
\rightarrow Q $ induces the map $
T\pi_{Q}: TT^* Q \rightarrow TQ. $

\begin{theo}
Assume that the triple $(T^*Q,\omega,H)$ is a Hamiltonian system
with the Hamiltonian vector field $X_H$, and $W: Q\rightarrow
\mathbb{R}$ is a given generating function. Then the following two assertions
are equivalent:\\
\noindent $(\mathrm{i})$ For every curve $\sigma: \mathbb{R}
\rightarrow Q $ satisfying that $\dot{\sigma}(t)= T\pi_Q
(X_H(\mathbf{d}W(\sigma(t))))$, $\forall t\in \mathbb{R}$,
$\mathbf{d}W \cdot \sigma $ is an integral curve of the Hamiltonian
vector field $X_H$.\\
\noindent $(\mathrm{ii})$ $W$ satisfies the Hamilton-Jacobi equation
$H(q^i,\frac{\partial W}{\partial q^i})=E, $ where $E$ is a
constant.
\end{theo}

From the proof of the Theorem2.1 given in
Abraham and Marsden \cite{abma78}, we know that
the assertion $(\mathrm{i})$, with equivalent to
Hamilton-Jacobi equation $(\mathrm{ii})$ by the generating function,
gives a geometric constraint condition of the canonical symplectic form
on the cotangent bundle $T^*Q$
for the Hamiltonian vector field of the system.
Thus, the Hamilton-Jacobi equation reveals the deeply internal relationships of
the generating function, the canonical symplectic form,
and the dynamical vector field of a Hamiltonian system.\\

However, from the discussion in \S 1, we know that
the set of Hamiltonian systems with symmetries on a cotangent
bundle is not complete under the Marsden-Weinstein reduction,
and the symplectic reduced system of a
Hamiltonian system with symmetry defined on the cotangent bundle
$T^*Q$ may not be a Hamiltonian system on a cotangent bundle,
so we cannot give the Hamilton-Jacobi equation for the Marsden-Weinstein
reduced Hamiltonian system as in Theorem 2.1;
we have to look for a new way.\\

It is worthy of noting that if we take that $\gamma=\mathbf{d}W$ in
Theorem 2.1, then $\gamma$ is a closed one-form on $Q$, and
the equation $\mathbf{d}(H \cdot \mathbf{d}W)=0$ is equivalent to
the Hamilton-Jacobi equation $H(q^i,\frac{\partial W}{\partial
q^i})=E$, where $E$ is a constant, which was called the classical
Hamilton-Jacobi equation. This result was used the
formulation of a geometric version of the Hamilton-Jacobi theorem for
Hamiltonian systems, see Cari\~{n}ena \textit{et al.} \cite{cagrmamamuro06,
cagrmamamuro10}. Moreover, noting that Theorem 2.1 was also
generalized in the context of a time-dependent Hamiltonian system by
Marsden and Ratiu in \cite{mara99}, and the Hamilton-Jacobi equation may
be regarded as a nonlinear partial differential equation for some
generating function $S$. Thus, the problem becomes how to choose a
time-dependent canonical transformation $\Psi: T^*Q\times \mathbb{R}
\rightarrow T^*Q\times \mathbb{R} $ which transforms the dynamical
vector field of a time-dependent Hamiltonian system to equilibrium
such that the generating function $S$ of $\Psi$ satisfies the
time-dependent Hamilton-Jacobi equation. In particular, for the
time-independent Hamiltonian system, one may look for a symplectic
map as the canonical transformation. This work offers an important
idea: that one can use the dynamical vector field of a Hamiltonian
system to describe the Hamilton-Jacobi equation. As a consequence, if we assume that
$\gamma: Q \rightarrow T^*Q$ is a closed one-form on $Q$, and define
that $X_H^\gamma = T\pi_{Q}\cdot X_H \cdot \gamma$, where $X_{H}$ is
the dynamical vector field of the Hamiltonian system $(T^*Q,\omega,H)$,
then we have that $X_H^\gamma$ and $X_H$ are $\gamma$-related; that
is, $T\gamma\cdot X_H^\gamma= X_H\cdot \gamma ,$ is equivalent to that
$\mathbf{d}(H \cdot \gamma)=0, $ which was given in Cari\~{n}ena
\textit{et al.} \cite{cagrmamamuro06, cagrmamamuro10}.
Motivated by the above research,
Wang in \cite{wa17} gave
an important notion and proved a key lemma, which is
a modification of the corresponding result of Abraham and Marsden
in \cite{abma78}.

\begin{defi}
The one-form $\gamma$ is said to be closed with respect to $T\pi_{Q}:
TT^* Q \rightarrow TQ $ if, for any $v, w \in TT^* Q, $ we have that
$\mathbf{d}\gamma(T\pi_{Q}(v),T\pi_{Q}(w))=0. $
\end{defi}

From the above definition we know that if $\gamma$ is a closed one-form,
then it must be closed with respect to $T\pi_{Q}: TT^* Q \rightarrow
TQ. $ Conversely, if $\gamma$ is closed with respect to
$T\pi_{Q}: TT^* Q \rightarrow TQ, $ then it may not be closed.
The following lemma is very important for our research, its proof
was given in Wang \cite{wa17}

\begin{lemm}
Assume that $\gamma: Q \rightarrow T^*Q$ is a one-form on $Q$, and
that $\lambda=\gamma \cdot \pi_{Q}: T^* Q \rightarrow T^* Q .$ Then
we have that the following two assertions hold:\\
\noindent $(\mathrm{i})$ For any $x, y \in TQ, \;
\gamma^*\omega(x,y)= -\mathbf{d}\gamma (x,y),$ and for any $v, w \in
TT^* Q, \; \lambda^*\omega(v,w)=\\ -\mathbf{d}\gamma(T\pi_{Q}(v), \;
T\pi_{Q}(w)),$
since $\omega$ is the canonical symplectic form on $T^*Q$. \\
\noindent $(\mathrm{ii})$ For any $v, w \in TT^* Q, \;
\omega(T\lambda \cdot v,w)= \omega(v, w-T\lambda \cdot
w)-\mathbf{d}\gamma(T\pi_{Q}(v), \; T\pi_{Q}(w)). $
\end{lemm}

By using the above Lemma 2.3, we can derive precisely
the geometric constraint conditions of
the Marsden-Weinstein reduced symplectic form for the
dynamical vector field of the Marsden-Weinstein reducible Hamiltonian
system; these are called the Type I and Type II Hamilton-Jacobi equations.
For convenience, the maps involved in
the following theorem and its proof are shown in Diagram-1.

\begin{center}
\hskip 0cm \xymatrix{ \mathbf{J}^{-1}(\mu) \ar[r]^{i_\mu}
& T^* Q \ar[d]_{X_{H\cdot \varepsilon}} \ar[dr]^{X_H^\varepsilon} \ar[r]^{\pi_Q}
& Q \ar[d]^{X_H^\gamma} \ar[r]^{\gamma}
& T^*Q \ar[d]_{X_H} \ar[dr]_{X_{h_\mu \cdot\bar{\varepsilon}}}\ar[r]^{\pi_\mu}
& (T^* Q)_\mu \ar[d]^{X_{h_\mu}} \\
& T(T^*Q)  & TQ \ar[l]^{T\gamma}
& T(T^*Q) \ar[l]^{T\pi_Q} \ar[r]_{T\pi_\mu} & T(T^* Q)_\mu }
\end{center}
$$\mbox{Diagram-1}$$

\begin{theo} ( Hamilton-Jacobi Theorem for the Marsden-Weinstein Reduced Hamiltonian System)
For the Marsden-Weinstein reducible Hamiltonian system
$(T^*Q,G,\omega, H)$ with the canonical symplectic form $\omega$ on $T^*Q$ and with the
reduced Hamiltonian system $((T^\ast Q)_\mu, \omega_\mu,h_\mu)$,
assume that $\gamma: Q \rightarrow T^*Q$ is a one-form
on $Q$, and that $\lambda=\gamma \cdot \pi_{Q}: T^* Q \rightarrow T^* Q $, and that
$\varepsilon: T^* Q \rightarrow T^* Q $ is a $G_\mu$-invariant
symplectic map, where $G_\mu$ is
the isotropy subgroup of the co-adjoint $G$-action at the point $\mu$.
Denote that $X_H^\gamma = T\pi_{Q}\cdot X_H \cdot \gamma$,
and that $X_H^\varepsilon = T\pi_{Q}\cdot X_H \cdot \varepsilon$,
where $X_{H}$ is the dynamical vector field of
the Marsden-Weinstein reducible Hamiltonian system $(T^*Q,G,\omega,H)$. Moreover,
assume that $\mu \in \mathfrak{g}^\ast $ is a regular value of the momentum
map $\mathbf{J}$, and that $\textmd{Im}(\gamma)\subset
\mathbf{J}^{-1}(\mu), $ and that $\gamma$ is $G_\mu$-invariant, and
that $\varepsilon(\mathbf{J}^{-1}(\mu)) \subset \mathbf{J}^{-1}(\mu). $
Denote that
$\bar{\gamma}=\pi_\mu(\gamma): Q \rightarrow (T^* Q)_\mu, $ and that
$\bar{\lambda}=\pi_\mu(\lambda): \mathbf{J}^{-1}(\mu) (\subset T^*Q) \rightarrow (T^* Q)_\mu $,
and that $\bar{\varepsilon}=
\pi_\mu(\varepsilon): \mathbf{J}^{-1}(\mu) (\subset T^*Q) \rightarrow (T^* Q)_\mu $.
Then the following two assertions hold:\\
\noindent $(\mathbf{i})$
If the one-form $\gamma: Q \rightarrow T^*Q $ is closed with respect to
$T\pi_Q: TT^* Q \rightarrow TQ, $
then $\bar{\gamma}$ is a solution of the Type I Hamilton-Jacobi equation
$T\bar{\gamma}\cdot X_H^\gamma= X_{h_\mu}\cdot \bar{\gamma} $
for the Marsden-Weinstein
reduced Hamiltonian system $((T^\ast Q)_\mu, \omega_\mu,h_\mu)$.\\
\noindent $(\mathbf{ii})$
The $\varepsilon$ and $\bar{\varepsilon} $ satisfy the Type II Hamilton-Jacobi equation
$T\bar{\gamma}\cdot X_H^\varepsilon= X_{h_\mu}\cdot \bar{\varepsilon}, $
if and only if they satisfy the equation
$T\bar{\varepsilon}\cdot(X_{h_\mu \cdot \bar{\varepsilon}})= T\bar{\lambda}\cdot X_H \cdot \varepsilon, $
where $X_{h_\mu}$ and$X_{h_\mu \cdot \bar{\varepsilon}} \in TT^* Q$
are the Hamiltonian vector fields of the Marsden-Weinstein reduced Hamiltonian
functions $h_{(\mu,a)}$ and $h_\mu \cdot \bar{\varepsilon}: T^* Q
\rightarrow \mathbb{R}, $ respectively.
\end{theo}
See the proof and the more details in Wang \cite{wa17}.
It is worth noting that
the Type I Hamilton-Jacobi equation
$T\bar{\gamma}\cdot X_H^\gamma= X_{h_\mu}\cdot
\bar{\gamma} $ is the equation of the
Marsden-Weinstein reduced differential one-form $\bar{\gamma}$, and that
the Type II Hamilton-Jacobi equation $T\bar{\gamma}\cdot X_H^\varepsilon
= X_{h_\mu}\cdot \bar{\varepsilon} $ is the equation of the symplectic
diffeomorphism map $\varepsilon$ and the Marsden-Weinstein reduced
symplectic diffeomorphism map $\bar{\varepsilon}. $
The reason that these are called the Type I and Type II Hamilton-Jacobi equations,
is that they are the development of the classical Hamilton-Jacobi equation
given by Theorem 2.1; (see Abraham and Marsden \cite{abma78}
and Wang \cite{wa17}). \\

We know that the orbit reduction for a Hamiltonian system with symmetry
is an alternative approach to symplectic reduction
given by Marle \cite{ma76},
and Kazhdan, Kostant and Sternberg \cite{kakost78},
and that it is different from the Marsden-Weinstein reduction.
We also derive precisely
the geometric constraint conditions of
the regular orbit reduced symplectic form for the
dynamical vector field of the regular orbit reducible Hamiltonian
system; that is, the Type I and Type II Hamilton-Jacobi equations
for the regular orbit reduced Hamiltonian system.
See the more details in Wang \cite{wa17}.

\section{Regular Controlled Hamiltonian System and Its Dynamics}

Since the symplectic reduced space of a Hamiltonian system
defined on the cotangent bundle of a configuration manifold may not be a
cotangent bundle, hence, the set of Hamiltonian systems with
symmetries on the cotangent bundle is not complete under the
Marsden-Weinstein reduction. This is a serious problem. If we define
directly a controlled Hamiltonian system with symmetry on a
cotangent bundle, then it is possible that the Marsden-Weinstein
reduced controlled Hamiltonian system may not have a definition.
However, from the classification of symplectic
reduced space of the cotangent bundle $T^* Q$, we know that
the Marsden-Weinstein symplectic
reduced space $((T^*Q)_\mu, \omega_\mu)$ is symplectically
diffeomorphic to a symplectic fiber bundle over $T^\ast (Q/G_\mu)$
with fiber to be the coadjoint orbit $\mathcal{O}_\mu$.
Thus, if we may define an RCH system on a symplectic fiber bundle,
then it is possible to describe uniformly the RCH system on $T^*Q$ and its
regular reduced RCH systems on the associated reduced spaces, and
we can study regular
reduction theory of the RCH systems with symplectic structures and
symmetries, as an extension of the Marsden-Weinstein reduction theory of
Hamiltonian systems under regular controlled Hamiltonian equivalence
conditions. This is why the authors in Marsden \textit{et al.} \cite{mawazh10}
set up the regular reduction theory of the RCH system on a symplectic fiber
bundle, by using momentum map and the associated reduced symplectic
form and from the viewpoint of completeness of regular symplectic
reduction.

\subsection{Regular Controlled Hamiltonian System}

First, our idea in Marsden \textit{et al.} \cite{mawazh10},
was that we first define a controlled Hamiltonian
system on $T^*Q$ by using the symplectic form, and the such system is
called a regular controlled Hamiltonian (RCH) system.
Next, we regard a Hamiltonian system on $T^*Q$
as a spacial case of an RCH system without external force and
control, and hence, the set of Hamiltonian systems on $T^*Q$ is a subset
of the set of RCH systems on $T^*Q$. On the other hand,
note that the symplectic reduced space on a cotangent bundle
is not complete under the Marsden-Weinstein reduction,
and hence, the symplectic reduced system of a Hamiltonian system with symmetry
defined on the cotangent bundle $T^*Q$ may not be a Hamiltonian system
on a cotangent bundle.
In order to describe, uniformly, the RCH systems defined on a cotangent
bundle and on the regular reduced spaces, we
first define an RCH system on a symplectic fiber bundle,
then we can obtain the RCH system on the cotangent bundle of a configuration
manifold as a special case.\\

Let $(\mathbb{E},M,\pi)$ be a fiber bundle, and for each point $x \in M$,
assume that the fiber $\mathbb{E}_x=\pi^{-1}(x)$ is a smooth submanifold of $\mathbb{E}$
with a symplectic form $\omega_{\mathbb{E}}(x)$; that is,
$(\mathbb{E}, \omega_{\mathbb{E}})$ is a
symplectic fiber bundle. If, for any Hamiltonian function $H: \mathbb{E} \rightarrow
\mathbb{R}$, we have a Hamiltonian vector field $X_H$
which satisfies the Hamilton's equation; that is,
$\mathbf{i}_{X_H}\omega_{\mathbb{E}}=\mathbf{d}H$,
then $(\mathbb{E}, \omega_{\mathbb{E}}, H )$ is a
Hamiltonian system. Moreover, considering the external force and
the control, we can define a kind of regular controlled Hamiltonian
(RCH) system on the symplectic fiber bundle $\mathbb{E}$ as follows:

\begin{defi} (RCH system) An RCH system on $\mathbb{E}$ is a 5-tuple
$(\mathbb{E}, \omega_{E}, H, F, W)$,
where $(\mathbb{E}, \omega_{\mathbb{E}}, H )$ is a
Hamiltonian system. and the function $H: \mathbb{E} \rightarrow \mathbb{R}$
is called the Hamiltonian, a fiber-preserving map $F: \mathbb{E}\rightarrow
\mathbb{E}$ is called the (external) force map, and a fiber sub-manifold
$W$ of $\mathbb{E}$ is called the control subset.
\end{defi}
Sometimes, $W$ is also denoted as the set of fiber-preserving
maps from $\mathbb{E}$ to $W$. When a feedback control law $u:
\mathbb{E}\rightarrow W$ is chosen,
the 5-tuple $(\mathbb{E}, \omega_{\mathbb{E}}, H,
F, u)$ is a closed-loop dynamical system. In particular, when $Q$
is a smooth manifold, and $T^\ast Q$ is its cotangent bundle with a
symplectic form $\omega$ (not necessarily canonical symplectic
form), then $(T^\ast Q, \omega )$ is a symplectic vector bundle. If
we take that $\mathbb{E}= T^* Q$, from the above definition we can obtain an RCH
system on the cotangent bundle $T^\ast Q$; that is, 5-tuple $(T^\ast
Q, \omega, H, F, W)$. Because
the fiber-preserving map $F: T^*Q\rightarrow T^*Q$ is the (external)
force map, which is the reason that the fiber-preserving map $F:
\mathbb{E}\rightarrow \mathbb{E}$ is called an (external) force map in above definition.
Here for convenience, we assume that all
controls appearing in this paper are the admissible controls.

\subsection{The Dynamics of an RCH System}

In order to describe the dynamics of an RCH system, we have to give a
good expression of the dynamical vector field of the RCH system, by using
the notation of the vertical lifted map of a vector along a fiber; see
Marsden \textit{et al.} \cite{mawazh10}.\\

First, for the notations of the vertical lifts along fiber,
we need to consider three case:
(1) $\pi: \mathbb{E} \rightarrow M$ is a fiber bundle;
(2) $\pi: \mathbb{E} \rightarrow M$ is a vector bundle;
(3) $\pi: \mathbb{E} \rightarrow M, \; \mathbb{E}= T^*Q, \; M=Q, $ is
a cotangent bundle, which is a special vector bundle.
For the cases (2) and (3), we can use the standard definition
of the vertical lift operator given in Marsden and Ratiu \cite{mara99}.
But for the case (1), the above operator cannot be used.
This question was found by one of referees who give his opinion
in a review report of our manuscript.
In order to deal with uniformly
the three cases, we have to give a new definition of the vertical
lifted maps of a vector along a fiber, and make it to be
not conflict with that given in Marsden and Ratiu \cite{mara99},
and it is not and cannot be an
extension of the definition of Marsden and Ratiu.\\

It is worthy of noting that there are two aspects in our new definition.
First, for two different points, $a_x,\; b_x $ in the fiber $\mathbb{E}_x, \; x\in M$,
how define the moving vertical part of a vector in one point $b_x$ to
another point $a_x$; Second, for a fiber-preserving map
$F: \mathbb{E} \rightarrow \mathbb{E}, $ we know that $a_x$ and $F_x(a_x)$
are the two points in $\mathbb{E}_x$,
how define the moving vertical part of a tangent vector in image
point $F_x(a_x)$ to $a_x$. The eventual goal is to give a good
expression of the dynamical vector field of RCH system by using the
notation of the vertical lift map of a vector along a fiber. Our
definitions are reasonable and clear, and should be stated
explicitly.\\

We first consider that $\mathbb{E}$ and $M$ are smooth manifolds, their tangent bundles
$T\mathbb{E}$ and $TM$ are vector bundles,
and for the fiber bundle $\pi: \mathbb{E} \rightarrow M$, we consider
the tangent mapping $T\pi: T\mathbb{E} \rightarrow TM$ and its kernel $ker
(T\pi)=\{\rho \in T\mathbb{E} \; | T\pi(\rho)=0\}$, which is a vector subbundle of
$T\mathbb{E}$. Denote that $V\mathbb{E}:= ker(T\pi)$, which is called a vertical bundle
of $\mathbb{E}$. Assume that there is a metric on $\mathbb{E}$, and that we take a
Levi-Civita connection $\mathcal{A}$ on $T\mathbb{E}$, and denote that $H\mathbb{E}:=
ker(\mathcal{A})$, which is called a horizontal bundle of $\mathbb{E}$, such
that $T\mathbb{E}= H\mathbb{E} \oplus V\mathbb{E}. $
Hence, for any $x\in M, \; a_x, b_x \in \mathbb{E}_x, $
any tangent vector $\rho(b_x)\in T_{b_x}\mathbb{E}$ can be split into
horizontal and vertical parts; that is, $\rho(b_x)=
\rho^h(b_x)\oplus \rho^v(b_x)$, where $\rho^h(b_x)\in H_{b_x}\mathbb{E}$ and
$\rho^v(b_x)\in V_{b_x}\mathbb{E}$. Let $\gamma$ be a geodesic in $\mathbb{E}_x$
connecting $a_x$ and $b_x$, and denote by $\rho^v_\gamma(a_x)$ a
tangent vector at $a_x$, which is a parallel displacement of the
vertical vector $\rho^v(b_x)$ along the geodesic $\gamma$ from $b_x$
to $a_x$. Since the angle between two vectors is invariant under a
parallel displacement along a geodesic, then
$T\pi(\rho^v_\gamma(a_x))=0, $ and hence, $\rho^v_\gamma(a_x) \in
V_{a_x}\mathbb{E}. $ Now, for $a_x, b_x \in \mathbb{E}_x $ and tangent vector
$\rho(b_x)\in T_{b_x}\mathbb{E}$, we can define the vertical lift map of a
vector along a fiber given by
$$\textnormal{vlift}: T\mathbb{E}_x \times \mathbb{E}_x \rightarrow T\mathbb{E}_x; \;\;
\textnormal{vlift}(\rho(b_x),a_x) = \rho^v_\gamma(a_x). $$
It is easy to check, from the basic fact in differential geometry,
that this map does not depend on the choice of the geodesic $\gamma$.\\

If $F: \mathbb{E} \rightarrow \mathbb{E}$ is a fiber-preserving map, for any $x\in M$, we have
that $F_x: \mathbb{E}_x \rightarrow \mathbb{E}_x$ and $TF_x: T\mathbb{E}_x \rightarrow T\mathbb{E}_x$,
then for any $a_x \in \mathbb{E}_x$ and $\rho\in T\mathbb{E}_x$, the vertical lift of
$\rho$ under the action of $F$ along a fiber is defined by
$$(\textnormal{vlift}(F_x)\rho)(a_x)
=\textnormal{vlift}((TF_x\rho)(F_x(a_x)), a_x)= (TF_x\rho)^v_\gamma(a_x), $$
where $\gamma$ is a geodesic in $\mathbb{E}_x$ connecting $F_x(a_x)$ and
$a_x$.\\

In particular, when $\pi: \mathbb{E} \rightarrow M$ is a vector bundle, for
any $x\in M$, the fiber $\mathbb{E}_x$ is a vector space. In this case, we
can choose the geodesic $\gamma$ to be a straight line, and the
vertical vector is invariant under a parallel displacement along a
straight line; that is, $\rho^v_\gamma(a_x)= \rho^v(b_x).$ Moreover,
when $\mathbb{E}= T^*Q$, by using the local trivialization of $TT^*Q$, we
have that $TT^*Q\cong TQ \times T^*Q$, locally. Note that $\pi: T^*Q
\rightarrow Q$, and $T\pi: TT^*Q \rightarrow TQ$, then in this case,
for any $\alpha_x, \; \beta_x \in T^*_x Q, \; x\in Q, $ we know that
$(0, \beta_x) \in V_{\beta_x}T^*_x Q, $ and hence, we can get that
$$ \textnormal{vlift}((0, \beta_x)(\beta_x), \alpha_x) = (0, \beta_x)(\alpha_x)
=\left.\frac{\mathrm{d}}{\mathrm{d}s}\right|_{s=0}(\alpha_x+s\beta_x) , $$
which is consistent with the definition of vertical lift operator along
a fiber given in Marsden and Ratiu \cite{mara99}.\\

For a given RCH System $(T^\ast Q, \omega, H, F, W)$, the dynamical
vector field of the associated Hamiltonian system $(T^\ast Q,
\omega, H) $ is $X_H$, which satisfies the equation
$\mathbf{i}_{X_H}\omega=\mathbf{d}H $.
Considering the external force $F: T^*Q \rightarrow T^*Q, $
which is a fiber-preserving map, by using the above
notation of the vertical lift map of a vector along a fiber, the
change of $X_H$ under the action of $F$ is such that
$$\textnormal{vlift}(F)X_H(\alpha_x)
= \textnormal{vlift}((TFX_H)(F(\alpha_x)), \alpha_x)
= (TFX_H)^v_\gamma(\alpha_x),$$
where $\alpha_x \in T^*_x Q, \; x\in Q $ and $\gamma$ is a straight
line in $T^*_x Q$ connecting $F_x(\alpha_x)$ and $\alpha_x$. In the
same way, when a feedback control law $u: T^\ast Q \rightarrow W,$
which also is a fiber-preserving map, is chosen,
the change of $X_H$ under the action of $u$ is such that
$$\textnormal{vlift}(u)X_H(\alpha_x)
= \textnormal{vlift}((TuX_H)(F(\alpha_x)), \alpha_x)
= (TuX_H)^v_\gamma(\alpha_x).$$
As a consequence, we can give an expression of the dynamical vector
field of the RCH system as follows:
\begin{theo}
The dynamical vector field of an RCH system $(T^\ast
Q,\omega,H,F,W)$ with a control law $u$ is the synthesis
of Hamiltonian vector field $X_H$ and its changes under the actions
of the external force $F$ and the control $u$; that is,
$$X_{(T^\ast Q,\omega,H,F,u)}(\alpha_x)
= X_H(\alpha_x)+ \textnormal{vlift}(F)X_H(\alpha_x)
+ \textnormal{vlift}(u)X_H(\alpha_x),$$ for any $\alpha_x \in T^*_x
Q, \; x\in Q $. For convenience, that is simply written as
\begin{equation}X_{(T^\ast Q,\omega,H,F,u)}
=X_H +\textnormal{vlift}(F) +\textnormal{vlift}(u), \label{3.1}
\end{equation}
\end{theo}
where $\textnormal{vlift}(F)=\textnormal{vlift}(F)X_H$,
and $\textnormal{vlift}(u)=\textnormal{vlift}(u)X_H$ are the
changes of $X_H$ under the actions of $F$ and $u$.
We also denote that $\textnormal{vlift}(W)=
\bigcup\{\textnormal{vlift}(u)X_H | \; u\in W\}$. It is
worthy of noting that in order to deduce and calculate easily, we
always use the simple expressions of the dynamical vector
field $X_{(T^\ast Q,\omega,H,F,u)}$.\\

From the expression (3.1) of the dynamical vector
field of an RCH system, we know that under the actions of the external force $F$
and the control $u$, in general, the dynamical vector
field is not Hamiltonian, and hence the RCH system is not
yet a Hamiltonian system. However,
it is a dynamical system closed relative to a
Hamiltonian system, and it can be explored and studied by the extending
methods for the external force and the control
in the study of Hamiltonian system.

\subsection{Controlled Hamiltonian Equivalence}

For two given Hamiltonian systems $(T^\ast
Q_i,\omega_i,H_i),$ $ i= 1,2,$ we say them to be
equivalent, if there exists a
diffeomorphism $\varphi: Q_1\rightarrow Q_2$ such that
their Hamiltonian vector fields $X_{H_i}, \; i=1,2 $ satisfy
the condition $X_{H_1}\cdot \varphi^\ast
=T(\varphi^\ast) X_{H_2}$, where the map
$\varphi^\ast= T^\ast \varphi:T^\ast Q_2\rightarrow T^\ast Q_1$
is the cotangent lifted map of $\varphi$, and the
map $T(\varphi^\ast):TT^\ast Q_2\rightarrow TT^\ast Q_1$ is the
tangent map of $\varphi^\ast$. From Marsden and Ratiu \cite{mara99},
we know that the condition $X_{H_1}\cdot \varphi^\ast
=T(\varphi^\ast) X_{H_2}$ is equivalent to the fact that
the map $\varphi^\ast:T^\ast Q_2\rightarrow T^\ast Q_1$
is symplectic with respect to the canonical symplectic forms
on $T^*Q_i, \; i=1,2.$ \\

For two given RCH systems $(T^\ast
Q_i,\omega_i,H_i,F_i,W_i),$ $ i= 1,2,$ we also want to define
their equivalence; that is, to look for a diffeomorphism
$\varphi: Q_1\rightarrow Q_2$ such that the condition
$X_{(T^\ast Q_1,\omega_1,H_1,F_1,W_1)}\cdot \varphi^\ast
=T(\varphi^\ast) X_{(T^\ast Q_2,\omega_2,H_2,F_2,W_2)}$ holds. However,
note that when an RCH system is given, the force map $F$ is
determined, but the feedback control law $u: T^\ast Q\rightarrow
W $ could be chosen. In order to emphasize explicitly the impact of the external force
and the control in the study of the RCH systems, by using
the above expression (3.1) of the dynamical vector field of the RCH system,
we can describe the feedback
control law to modify the structure of the RCH system, and the controlled
Hamiltonian matching conditions and RCH-equivalence are induced as
follows:
\begin{defi}
(RCH-equivalence) Suppose that we have two RCH systems $(T^\ast
Q_i,\omega_i,H_i,F_i,W_i),$ $ i= 1,2,$ we say that they are
RCH-equivalent, or simply, that $(T^\ast
Q_1,\omega_1,H_1,F_1,W_1)\stackrel{RCH}{\sim}\\ (T^\ast
Q_2,\omega_2,H_2,F_2,W_2)$, if there exists a
diffeomorphism $\varphi: Q_1\rightarrow Q_2$ such that the
following controlled Hamiltonian matching conditions hold:

\noindent {\bf RCH-1:} The control subsets $W_i, \; i=1,2$ satisfy
the condition $W_1=\varphi^\ast (W_2),$ where the map
$\varphi^\ast= T^\ast \varphi:T^\ast Q_2\rightarrow T^\ast Q_1$
is cotangent lifted map of $\varphi$.

\noindent {\bf RCH-2:} For each control law $u_1:
T^\ast Q_1 \rightarrow W_1, $ there exists the control law $u_2:
T^\ast Q_2 \rightarrow W_2 $  such that the two
closed-loop dynamical systems produce the same dynamical vector fields; that is,
$X_{(T^\ast Q_1,\omega_1,H_1,F_1,u_1)}\cdot \varphi^\ast
=T(\varphi^\ast) X_{(T^\ast Q_2,\omega_2,H_2,F_2,u_2)}$,
where the map $T(\varphi^\ast):TT^\ast Q_2\rightarrow TT^\ast Q_1$
is the tangent map of $\varphi^\ast$.
\end{defi}

From the expression (3.1) of the dynamical vector field of the RCH system
and the condition $X_{(T^\ast Q_1,\omega_1,H_1,F_1,u_1)}\cdot \varphi^\ast
=T(\varphi^\ast) X_{(T^\ast Q_2,\omega_2,H_2,F_2,u_2)}$, we have that
$$
(X_{H_1}+\textnormal{vlift}(F_1)X_{H_1}+\textnormal{vlift}(u_1)X_{H_1})\cdot \varphi^\ast
= T(\varphi^\ast)[X_{H_2}+\textnormal{vlift}(F_2)X_{H_2}+\textnormal{vlift}(u_2)X_{H_2}].
$$
By using the notation of vertical lift map of a vector along a fiber,
for $\alpha_x \in T_x^\ast Q_2, \; x \in Q_2$, we have that
\begin{align*}
T(\varphi^\ast)\textnormal{vlift}(F_2)X_{H_2}(\alpha_x)
&=T(\varphi^\ast)\textnormal{vlift}((TF_2X_{H_2})(F_2(\alpha_x)), \alpha_x)\\
&=\textnormal{vlift}(T(\varphi^\ast)\cdot TF_2\cdot T(\varphi_\ast)X_{H_2}
(\varphi^\ast F_2 \varphi_\ast(\varphi^\ast \alpha_x)),\varphi^\ast \alpha)\\
&=\textnormal{vlift}(T(\varphi^\ast F_2 \varphi_\ast)X_{H_2}
(\varphi^\ast F_2 \varphi_\ast(\varphi^\ast \alpha_x)),\varphi^\ast \alpha)\\
&=\textnormal{vlift}(\varphi^\ast
F_2\varphi_\ast)X_{H_2}(\varphi^\ast \alpha_x),
\end{align*}
where the map $\varphi_\ast=(\varphi^{-1})^\ast: T^\ast Q_1\rightarrow T^\ast Q_2$.
In the same way, we have that
$T(\varphi^\ast)\textnormal{vlift}(u_2)X_{H_2}=\textnormal{vlift}(\varphi^\ast
u_2\varphi_\ast)X_{H_2}\cdot \varphi^\ast$.
Thus, the explicit relation
between the two control laws $u_i \in W_i, \; i=1,2$ in {\bf RCH-2} is given by
\begin{align}
& (\textnormal{vlift}(u_1) -\textnormal{vlift}(\varphi^\ast u_2\varphi_\ast))\cdot \varphi^\ast \nonumber \\
& = -X_{H_1}\cdot \varphi^\ast +T\varphi^\ast (X_{H_2})+
(-\textnormal{vlift}(F_1)+\textnormal{vlift}(\varphi^\ast F_2
\varphi_\ast))\cdot \varphi^\ast.
\label{3.2}\end{align}

From the above relation (3.2) we know that, when two RCH systems $(T^\ast
Q_i,\omega_i,H_i,F_i,W_i),$ $ i= 1,2,$ are RCH-equivalent with respect to $\varphi^*$,
the corresponding Hamiltonian systems $(T^\ast Q_i,\omega_i,H_i),$
$ i= 1,2,$ may not be equivalent with respect to $\varphi^*$. If
two corresponding Hamiltonian systems are also equivalent with respect to $\varphi^*$,
then the control laws $u_i : T^*Q \rightarrow W_i, \; i=1,2$ and the
external forces $F_i: T^*Q_i \rightarrow T^*Q_i, \; i=1,2$ in {\bf RCH-2}
must satisfy the following condition
\begin{equation}
\textnormal{vlift}(u_1) -\textnormal{vlift}(\varphi^\ast u_2\varphi_\ast)
= -\textnormal{vlift}(F_1)+\textnormal{vlift}(\varphi^\ast F_2 \varphi_\ast).
\label{3.3}\end{equation}

\section{Regular Reductions for the RCH System}

We know that when the external force and the control of an RCH
system $(T^*Q,\omega,H,F,W)$ are both zeros, in this case the RCH system
is just a Hamiltonian system $(T^*Q,\omega,H)$.
Thus, we can regard a Hamiltonian system defined on $T^*Q$
as a spacial case of an RCH system without external force and
control. As a consequence, the set of Hamiltonian systems with symmetries
on $T^*Q$ is a subset of the set of RCH systems with symmetries on $T^*Q$.
If we admit the Marsden-Weinstein symplectic reduction of a Hamiltonian system
with symmetry and a momentum map, then we may study the regular
point reduction of an RCH system with symmetry and a momentum map,
as an extension of the Marsden-Weinstein symplectic reduction of
a Hamiltonian system under regular controlled Hamiltonian equivalence
conditions. In order to do these, in this section
we shall give the regular point reduction
and the regular orbit reduction of an RCH system,
by using the Marsden-Weinstein reduction and
regular orbit reduction for a Hamiltonian system, respectively.

\subsection{Regular Point Reduction of an RCH System}

In the subsection, we consider the RCH system with symmetry and a momentum map,
and give the regular point reduced RCH system and
the RpCH-equivalence for the regular point reducible RCH system,
and prove the regular point reduction theorem.\\

At first, we consider the regular point reducible RCH system.
Let $Q$ be a smooth manifold and $T^\ast Q$ is its cotangent bundle
with the symplectic form $\omega$. Let $\Phi: G\times Q\rightarrow
Q$ be a smooth left action of a Lie group $G$ on $Q$, which is free
and proper, so the cotangent lifted left action $\Phi^{T^\ast}:
G\times T^\ast Q\rightarrow T^\ast Q$ is also free and
proper. Assume that the cotangent lifted action is symplectic
with respect to the symplectic form $\omega$,
and that it admits an
$\operatorname{Ad}^\ast$-equivariant momentum map $\mathbf{J}:T^\ast
Q\rightarrow \mathfrak{g}^\ast$, where $\mathfrak{g}$ is the Lie
algebra of $G$ and $\mathfrak{g}^\ast$ is the dual of
$\mathfrak{g}$. Let $\mu\in\mathfrak{g}^\ast$ be a regular value of
$\mathbf{J}$ and denote that $G_\mu$ is the isotropy subgroup of the
coadjoint $G$-action at the point $\mu\in\mathfrak{g}^\ast$, and that it
is defined by $G_\mu=\{g\in G|\operatorname{Ad}_g^\ast \mu=\mu \}$.
Since $G_\mu (\subset G)$ acts freely and properly on $Q$ and on
$T^\ast Q$, then $Q_\mu=Q/G_\mu$ is a smooth manifold and that the
canonical projection $\rho_\mu:Q\rightarrow Q_\mu$ is a surjective
submersion. It follows that $G_\mu$ acts also freely and properly on
$\mathbf{J}^{-1}(\mu)$, so that the space $(T^\ast
Q)_\mu=\mathbf{J}^{-1}(\mu)/G_\mu$ is a symplectic manifold with the
symplectic form $\omega_\mu$ uniquely characterized by the relation that
\begin{equation}\pi_\mu^\ast \omega_\mu=i_\mu^\ast
\omega. \label{4.1}\end{equation} The map
$i_\mu:\mathbf{J}^{-1}(\mu)\rightarrow T^\ast Q$ is the inclusion
and the map $\pi_\mu:\mathbf{J}^{-1}(\mu)\rightarrow (T^\ast Q)_\mu$ is the
projection. The pair $((T^\ast Q)_\mu,\omega_\mu)$ is called the
Marsden-Weinstein reduced space of $(T^\ast Q,\omega)$ at $\mu$.\\

Assume that $H: T^\ast Q\rightarrow \mathbb{R}$ is a $G$-invariant
Hamiltonian, the flow $F_t$ of the Hamiltonian vector field $X_H$
leaves the connected components of $\mathbf{J}^{-1}(\mu)$ invariant
and commutes with the $G$-action, so it induces a flow $f_t^\mu$ on
$(T^\ast Q)_\mu$ defined by $f_t^\mu\cdot \pi_\mu=\pi_\mu \cdot
F_t\cdot i_\mu$, and the vector field $X_{h_\mu}$ generated by the
flow $f_t^\mu$ on $((T^\ast Q)_\mu,\omega_\mu)$ is Hamiltonian with
the associated regular point reduced Hamiltonian function
$h_\mu:(T^\ast Q)_\mu\rightarrow \mathbb{R}$ defined by
$h_\mu\cdot\pi_\mu=H\cdot i_\mu$, and the Hamiltonian vector fields
$X_H$ and $X_{h_\mu}$ are $\pi_\mu$-related.
Moreover, assume that the fiber-preserving map $F:T^\ast Q\rightarrow T^\ast
Q$ and the control subset $W$ of\; $T^\ast Q$ are both $G$-invariant.
In order to get the $R_p$-reduced RCH system, we also assume that
$F(\mathbf{J}^{-1}(\mu))\subset \mathbf{J}^{-1}(\mu)£¬$ and that $W \cap
\mathbf{J}^{-1}(\mu)\neq \emptyset $.
Thus, we can introduce a regular point
reducible RCH system as follows (see Marsden \textit{et al.} \cite{mawazh10},
Wang \cite{wa18} and Wang \cite{wa13}):
\begin{defi}
(Regular Point Reducible RCH System) A 6-tuple $(T^\ast Q, G,
\omega, H, F, W)$ with
the canonical symplectic form $\omega$ on $T^*Q$,
where the Hamiltonian $H:T^\ast Q\rightarrow
\mathbb{R}$, the fiber-preserving map $F:T^\ast Q\rightarrow T^\ast
Q$ and the fiber submanifold $W$ of\; $T^\ast Q$ are all
$G$-invariant, is called a regular point reducible RCH system if
there exists a point $\mu\in\mathfrak{g}^\ast$ which is a regular
value of the momentum map $\mathbf{J}$ such that the regular point
reduced system; that is, the 5-tuple $((T^\ast Q)_\mu,
\omega_\mu,h_\mu,f_\mu,W_\mu)$, where $(T^\ast
Q)_\mu=\mathbf{J}^{-1}(\mu)/G_\mu$, $\pi_\mu^\ast
\omega_\mu=i_\mu^\ast\omega$, $h_\mu\cdot \pi_\mu=H\cdot i_\mu$,
$F(\mathbf{J}^{-1}(\mu))\subset \mathbf{J}^{-1}(\mu) $, $f_\mu\cdot
\pi_\mu=\pi_\mu \cdot F\cdot i_\mu$, $W \cap
\mathbf{J}^{-1}(\mu)\neq \emptyset $ and $W_\mu=\pi_\mu(W\cap
\mathbf{J}^{-1}(\mu))$, is an RCH system, which is simply written as
the $R_p$-reduced RCH system. Where $((T^\ast Q)_\mu,\omega_\mu)$ is the
$R_p$-reduced space, the function $h_\mu:(T^\ast Q)_\mu\rightarrow
\mathbb{R}$ is called the $R_p$-reduced Hamiltonian, the fiber-preserving
map $f_\mu:(T^\ast Q)_\mu\rightarrow (T^\ast Q)_\mu$ is called the
$R_p$-reduced (external) force map, and $W_\mu$ is a fiber submanifold of
\;$(T^\ast Q)_\mu$ and is called the $R_p$-reduced control subset.
\end{defi}

It is worthy of noting that for the regular point reducible RCH system
$(T^\ast Q,G,\omega,H,F,W)$, the $G$-invariant external force map
$F: T^*Q \rightarrow T^*Q $ has to satisfy the conditions
$F(\mathbf{J}^{-1}(\mu))\subset \mathbf{J}^{-1}(\mu), $ and
$f_\mu\cdot \pi_\mu=\pi_\mu \cdot F\cdot i_\mu $ such that we can
define the $R_p$-reduced external force map $f_\mu:(T^\ast
Q)_\mu\rightarrow (T^\ast Q)_\mu. $ The condition $W \cap
\mathbf{J}^{-1}(\mu)\neq \emptyset $ in above definition makes that
the $G$-invariant control subset $W\cap \mathbf{J}^{-1}(\mu)$ can be
reduced and the $R_p$-reduced control subset is $W_\mu= \pi_\mu(W\cap
\mathbf{J}^{-1}(\mu))$. \\

In the following we consider the RCH system with symmetry and a momentum map,
and give the RpCH-equivalence for the regular point reducible RCH systems,
and prove the regular point reduction theorem.
Denote by $X_{(T^\ast Q,G,\omega,H,F,u)}$ the dynamical vector field of the
regular point reducible RCH system $(T^\ast Q,G,\omega,
H,F,W)$ with a control law $u$. Assume that it can be expressed by
\begin{equation}X_{(T^\ast Q,G,\omega,H,F,u)}
=X_H+\textnormal{vlift}(F)+\textnormal{vlift}(u).
\label{4.2}\end{equation}

By using the above expression (4.2),
for the regular point reducible RCH system we can
introduce the regular point reducible controlled Hamiltonian
equivalence (RpCH-equivalence) as follows:

\begin{defi}(RpCH-equivalence)
Suppose that we have two regular point reducible RCH systems
$(T^\ast Q_i, G_i,\omega_i,H_i, F_i, W_i),\; i=1,2$, we
say that they are RpCH-equivalent, or simply, that $(T^\ast Q_1,
G_1,\omega_1,H_1,F_1,W_1)\stackrel{RpCH}{\sim}(T^\ast
Q_2,G_2,\omega_2,H_2,F_2,W_2)$, if there exists a
diffeomorphism $\varphi:Q_1\rightarrow Q_2$ such that the following
regular point reducible controlled Hamiltonian matching conditions hold:

\noindent {\bf RpCH-1:} For $\mu_i\in \mathfrak{g}^\ast_i $, the
regular reducible points of the RCH systems $(T^\ast Q_i, G_i,\omega_i,
H_i, F_i, W_i),\\ i=1,2$, the map
$\varphi_\mu^\ast=i_{\mu_1}^{-1}\cdot\varphi^\ast\cdot i_{\mu_2}:
\mathbf{J}_2^{-1}(\mu_2)\rightarrow \mathbf{J}_1^{-1}(\mu_1)$ is
$(G_{2\mu_2},G_{1\mu_1})$-equivariant,
and $W_1\cap \mathbf{J}_1^{-1}(\mu_1) = \varphi_\mu^\ast (W_2\cap
\mathbf{J}_2^{-1}(\mu_2))$, where $\mu=(\mu_1,
\mu_2)$, and denote by $i_{\mu_1}^{-1}(S)$ the pre-image of a subset
$S\subset T^\ast Q_1$ for the map
$i_{\mu_1}:\mathbf{J}_1^{-1}(\mu_1)\rightarrow T^\ast Q_1$.

\noindent {\bf RpCH-2:} For each control law $u_1:
T^\ast Q_1 \rightarrow W_1, $ there exists the control law $u_2:
T^\ast Q_2 \rightarrow W_2 $  such that the two
closed-loop dynamical systems produce the same dynamical vector fields; that is,
$X_{(T^\ast Q_1,G_1,\omega_1,H_1,F_1,u_1)}\cdot \varphi^\ast
=T(\varphi^\ast) X_{(T^\ast Q_2,G_2,\omega_2,H_2,F_2,u_2)}$.
\end{defi}

It is worthy of noting that for the regular point reducible RCH
system, the induced equivalent map $\varphi^*$
keeps the equivariance of the $G$-action
at the regular reducible point.
If an $R_p$-reduced feedback control law $u_\mu:(T^\ast
Q)_\mu\rightarrow W_\mu$ is chosen, the $R_p$-reduced RCH
system $((T^\ast Q)_\mu, \omega_\mu, h_\mu, f_\mu, u_\mu)$ is a
closed-loop regular dynamic system with a control law $u_\mu$.
Assume that its vector field $X_{((T^\ast Q)_\mu, \omega_\mu, h_\mu,
f_\mu, u_\mu)}$ can be expressed by
\begin{equation}X_{((T^\ast Q)_\mu, \omega_\mu, h_\mu, f_\mu, u_\mu)}
=X_{h_\mu}+\textnormal{vlift}(f_\mu)+\textnormal{vlift}(u_\mu),
\label{4.3}
\end{equation}
where $X_{h_\mu}$ is the Hamiltonian vector field of
the $R_p$-reduced Hamiltonian $h_\mu$, and $\textnormal{vlift}(f_\mu)=
\textnormal{vlift}(f_\mu)X_{h_\mu}$, $\textnormal{vlift}(u_\mu)=
\textnormal{vlift}(u_\mu)X_{h_\mu}$ are the changes of $X_{h_\mu}$
under the actions of the $R_p$-reduced external force $f_\mu$
and the $R_p$-reduced control law $u_\mu$.
The dynamical vector fields of
the regular point reducible RCH system $(T^\ast Q,G,\omega,
H,F,u)$ and the $R_p$-reduced RCH system $((T^\ast Q)_\mu,
\omega_\mu,h_\mu,f_\mu, u_\mu)$ satisfy the condition that
\begin{equation}X_{((T^\ast Q)_\mu, \omega_\mu, h_\mu, f_\mu,
u_\mu)}\cdot \pi_\mu=T\pi_\mu\cdot X_{(T^\ast
Q,G,\omega,H,F,u)}\cdot i_\mu. \label{4.4}
\end{equation}

Thus, we can prove
the following regular point reduction theorem for the RCH
system, which explains the relationship between the RpCH-equivalence
for the regular point reducible RCH systems with symmetries and momentum maps
and the RCH-equivalence for the associated $R_p$-reduced RCH systems, its
proof was given in Marsden \textit{et al.} \cite{mawazh10}
and Wang \cite{wa18}. This theorem can
be regarded as an extension of the regular point reduction theorem
of Hamiltonian system under the regular controlled Hamiltonian
equivalence conditions.
\begin{theo}
Two regular point reducible RCH systems $(T^\ast Q_i, G_i, \omega_i,
H_i, F_i,W_i)$, $i=1,2,$ are RpCH-equivalent if and only
if the associated $R_p$-reduced RCH systems $((T^\ast
Q_i)_{\mu_i},\omega_{i\mu_i},h_{i\mu_i},f_{i\mu_i},\\
W_{i\mu_i}),$ $i=1,2,$ are RCH-equivalent.
\end{theo}

It is worthy of noting that when the external force and the control of a
regular point reducible RCH system $(T^*Q,G,\omega,H,F,W)$ are both zeros;
that is, $F=0 $ and $W=\emptyset$, in this case the RCH system
is just a regular point reducible Hamiltonian system $(T^*Q,G,\omega,H)$ itself.
Then the following theorem explains the relationship between the equivalence
for the regular point reducible Hamiltonian systems with symmetries and the
equivalence for the associated $R_p$-reduced Hamiltonian systems.

\begin{theo}
Two regular point reducible Hamiltonian systems $(T^\ast Q_i, G_i, \omega_i,
H_i )$, $i=1,2,$ are equivalent if and only
if the associated $R_p$-reduced Hamiltonian systems $((T^\ast
Q_i)_{\mu_i},\omega_{i\mu_i},h_{i\mu_i}),$ $i=1,2,$ are equivalent.
\end{theo}
See the proof and the more details in Wang \cite{wa18}.

\subsection{Regular Orbit Reduction of the RCH System}

Since the orbit reduction of a Hamiltonian system is an alternative
approach to symplectic reduction, which is different from the
Marsden-Weinstein reduction, in the subsection,
we consider the RCH system with symmetry and a momentum map,
and give the regular orbit reduced RCH system and
the RoCH-equivalence for the regular orbit reducible RCH system,
and prove the regular orbit reduction theorem.\\

First, we consider the regular orbit reducible RCH system.
Assume that the cotangent lifted left action
$\Phi^{T^\ast}:G\times T^\ast Q\rightarrow T^\ast Q$ is
symplectic, free and proper, and that the action admits an
$\operatorname{Ad}^\ast$-equivariant momentum map $\mathbf{J}:T^\ast
Q\rightarrow \mathfrak{g}^\ast$. Let $\mu\in \mathfrak{g}^\ast$ be a
regular value of the momentum map $\mathbf{J}$ and
$\mathcal{O}_\mu=G\cdot \mu\subset \mathfrak{g}^\ast$ be the
$G$-orbit of the coadjoint $G$-action through the point $\mu$. Since
$G$ acts freely, properly and symplectically on $T^\ast Q$, the
quotient space $(T^\ast Q)_{\mathcal{O}_\mu}=
\mathbf{J}^{-1}(\mathcal{O}_\mu)/G$ is a regular quotient symplectic
manifold with the symplectic form $\omega_{\mathcal{O}_\mu}$
uniquely characterized by the relation that
\begin{equation}i_{\mathcal{O}_\mu}^\ast \omega=\pi_{\mathcal{O}_{\mu}}^\ast
\omega_{\mathcal{O}
_\mu}+\mathbf{J}_{\mathcal{O}_\mu}^\ast\omega_{\mathcal{O}_\mu}^+,
\label{4.5}
\end{equation}
where $\mathbf{J}_{\mathcal{O}_\mu}$ is
the restriction of the momentum map $\mathbf{J}$ to
$\mathbf{J}^{-1}(\mathcal{O}_\mu)$; that is,
$\mathbf{J}_{\mathcal{O}_\mu}=\mathbf{J}\cdot i_{\mathcal{O}_\mu}$
and $\omega_{\mathcal{O}_\mu}^+$ is the $+$-symplectic structure on
the orbit $\mathcal{O}_\mu$ given by
\begin{equation}\omega_{\mathcal{O}_\mu}^
+(\nu)(\xi_{\mathfrak{g}^\ast}(\nu),\eta_{\mathfrak{g}^\ast}(\nu))
=<\nu,[\xi,\eta]>,\;\; \forall\;\nu\in\mathcal{O}_\mu, \;
\xi,\eta\in \mathfrak{g}. \label{4.6}
\end{equation}
The maps
$i_{\mathcal{O}_\mu}:\mathbf{J}^{-1}(\mathcal{O}_\mu)\rightarrow
T^\ast Q$ and
$\pi_{\mathcal{O}_\mu}:\mathbf{J}^{-1}(\mathcal{O}_\mu)\rightarrow
(T^\ast Q)_{\mathcal{O}_\mu}$ are the natural injection and the
projection, respectively. The pair $((T^\ast
Q)_{\mathcal{O}_\mu},\omega_{\mathcal{O}_\mu})$ is called the
symplectic orbit reduced space of $(T^\ast Q,\omega)$ at $\mu$.\\

Assume that $H:T^\ast Q\rightarrow \mathbb{R}$ is a $G$-invariant
Hamiltonian, the flow $F_t$ of the Hamiltonian vector field $X_H$
leaves the connected components of
$\mathbf{J}^{-1}(\mathcal{O}_\mu)$ invariant and commutes with the
$G$-action, so it induces a flow $f_t^{\mathcal{O}_\mu}$ on $(T^\ast
Q)_{\mathcal{O}_\mu}$, defined by $f_t^{\mathcal{O}_\mu}\cdot
\pi_{\mathcal{O}_\mu}=\pi_{\mathcal{O}_\mu} \cdot F_t\cdot
i_{\mathcal{O}_\mu}$, and the vector field $X_{h_{\mathcal{O}_\mu}}$
generated by the flow $f_t^{\mathcal{O}_\mu}$ on $((T^\ast
Q)_{\mathcal{O}_\mu},\omega_{\mathcal{O}_\mu})$ is Hamiltonian with
the associated regular orbit reduced Hamiltonian function
$h_{\mathcal{O}_\mu}:(T^\ast Q)_{\mathcal{O}_\mu}\rightarrow
\mathbb{R}$ defined by $h_{\mathcal{O}_\mu}\cdot
\pi_{\mathcal{O}_\mu}= H\cdot i_{\mathcal{O}_\mu}$, and the
Hamiltonian vector fields $X_H$ and $X_{h_{\mathcal{O}_\mu}}$ are
$\pi_{\mathcal{O}_\mu}$-related. In general, we
maybe thought that the structure of the symplectic orbit reduced
space $((T^\ast Q)_{\mathcal{O}_\mu},\omega_{\mathcal{O}_\mu})$ is
more complex than that of the symplectic point reduced space
$((T^\ast Q)_\mu,\omega_\mu)$, but, from the regular reduction
diagram ( see Ortega and Ratiu \cite{orra04}),
we know that the regular orbit reduced space $((T^\ast
Q)_{\mathcal{O}_\mu},\omega_{\mathcal{O}_\mu})$ is symplectically
diffeomorphic to the regular point reduced space $((T^*Q)_\mu,
\omega_\mu)$, and hence, is also symplectically diffeomorphic to a
symplectic fiber bundle. Moreover, assume that the fiber-preserving map
$F:T^\ast Q\rightarrow T^\ast Q$ and the control subset
$W$ of\; $T^\ast Q$ are both $G$-invariant.
In order to get the $R_o$-reduced RCH system, we also assume that
$F(\mathbf{J}^{-1}(\mathcal{O}_\mu ))\subset \mathbf{J}^{-1}(\mathcal{O}_\mu)£¬$
and that $W \cap \mathbf{J}^{-1}(\mathcal{O}_\mu)\neq \emptyset $.
Thus, we can introduce a regular orbit
reducible RCH system as follows (see Marsden \textit{et al.} \cite{mawazh10},
Wang \cite{wa18}):

\begin{defi}
(Regular Orbit Reducible RCH System) A 6-tuple $(T^\ast Q, G,
\omega,H,F,W)$, where the Hamiltonian $H: T^\ast Q\rightarrow
\mathbb{R}$, the fiber-preserving map $F: T^\ast Q\rightarrow T^\ast
Q$ and the fiber submanifold $W$ of $T^\ast Q$ are all
$G$-invariant, is called a regular orbit reducible RCH system, if
there exists an orbit $\mathcal{O}_\mu, \; \mu\in\mathfrak{g}^\ast$,
where $\mu$ is a regular value of the momentum map $\mathbf{J}$,
such that the regular orbit reduced system; that is, the 5-tuple
$((T^\ast Q)_{\mathcal{O}_\mu},
\omega_{\mathcal{O}_\mu},h_{\mathcal{O}_\mu},f_{\mathcal{O}_\mu},
W_{\mathcal{O}_\mu})$, where $(T^\ast
Q)_{\mathcal{O}_\mu}=\mathbf{J}^{-1}(\mathcal{O}_\mu)/G$,
$\pi_{\mathcal{O}_\mu}^\ast \omega_{\mathcal{O}_\mu}
=i_{\mathcal{O}_\mu}^\ast\omega-\mathbf{J}_{\mathcal{O}_\mu}^\ast\omega_{\mathcal{O}_\mu}^+$,
$h_{\mathcal{O}_\mu}\cdot \pi_{\mathcal{O}_\mu} =H\cdot
i_{\mathcal{O}_\mu}$, $F(\mathbf{J}^{-1}(\mathcal{O}_\mu))\subset
\mathbf{J}^{-1}(\mathcal{O}_\mu)$, $f_{\mathcal{O}_\mu}\cdot
\pi_{\mathcal{O}_\mu}=\pi_{\mathcal{O}_\mu}\cdot F\cdot
i_{\mathcal{O}_\mu}$, and $W \cap
\mathbf{J}^{-1}(\mathcal{O}_\mu)\neq \emptyset $,
$W_{\mathcal{O}_\mu}=\pi_{\mathcal{O}_\mu}(W \cap
\mathbf{J}^{-1}(\mathcal{O}_\mu))$, is an RCH system, that is simply
written as the $R_o$-reduced RCH system. Where $((T^\ast
Q)_{\mathcal{O}_\mu},\omega_{\mathcal{O}_\mu})$ is the $R_o$-reduced
space, the function $h_{\mathcal{O}_\mu}:(T^\ast
Q)_{\mathcal{O}_\mu}\rightarrow \mathbb{R}$ is called the $R_o$-reduced
Hamiltonian, the fiber-preserving map $f_{\mathcal{O}_\mu}:(T^\ast
Q)_{\mathcal{O}_\mu} \rightarrow (T^\ast Q)_{\mathcal{O}_\mu}$ is
called the $R_o$-reduced (external) force map, and $W_{\mathcal{O}_\mu}$ is a
fiber submanifold of $(T^\ast Q)_{\mathcal{O}_\mu}$, and is called
the $R_o$-reduced control subset.
\end{defi}

It is worthy of noting that for the regular orbit reducible RCH system
$(T^\ast Q,G,\omega,H,F,W)$, the $G$-invariant external force map
$F: T^*Q \rightarrow T^*Q $ has to satisfy the conditions that
$F(\mathbf{J}^{-1}(\mathcal{O}_\mu))\subset
\mathbf{J}^{-1}(\mathcal{O}_\mu), $ and $f_{\mathcal{O}_\mu}\cdot
\pi_{\mathcal{O}_\mu}=\pi_{\mathcal{O}_\mu}\cdot F\cdot
i_{\mathcal{O}_\mu}$, such that we can define the reduced external
force map $f_{\mathcal{O}_\mu}:(T^\ast Q)_{\mathcal{O}_\mu}
\rightarrow (T^\ast Q)_{\mathcal{O}_\mu}. $ The condition that $W \cap
\mathbf{J}^{-1}(\mathcal{O}_\mu)\neq \emptyset $ in above definition
makes that the  $G$-invariant control subset $W \cap
\mathbf{J}^{-1}(\mathcal{O}_\mu)$ can be reduced and that the reduced
control subset is $W_{\mathcal{O}_\mu}= \pi_{\mathcal{O}_\mu}(W \cap
\mathbf{J}^{-1}(\mathcal{O}_\mu))$. \\

In the following we consider the RCH system with symmetry and a momentum map,
and give the RoCH-equivalence for the regular orbit reducible RCH systems,
and prove the regular orbit reduction theorem.
Denote by $X_{(T^\ast Q,G,\omega,H,F,u)}$ the dynamical vector field of the
regular orbit reducible RCH system $(T^\ast Q, G,\omega, H,F,W)$
with a control law $u$. Assume that it can be expressed by
\begin{equation}X_{(T^\ast Q,G,\omega,H,F,u)}
=X_H+\textnormal{vlift}(F)+\textnormal{vlift}(u).\label{4.7}
\end{equation}

By using the above expression (4.7),
for the regular orbit reducible RCH system we can
introduce the regular orbit reducible controlled Hamiltonian
equivalence (RoCH-equivalence) as follows.
\begin{defi}
(RoCH-equivalence) Suppose that we have two regular orbit reducible
RCH systems $(T^\ast Q_i, G_i, \omega_i, H_i, F_i, W_i)$, $i=1,2$,
we say that they are RoCH-equivalent, or simply, that $(T^\ast Q_1, G_1,
\omega_1, H_1, F_1, W_1)\stackrel{RoCH}{\sim}(T^\ast Q_2, G_2,
\omega_2, H_2, F_2, W_2)$, if there exists a diffeomorphism
$\varphi:Q_1\rightarrow Q_2$ such that the following controlled Hamiltonian
matching conditions hold:

\noindent {\bf RoCH-1:} For $\mathcal{O}_{\mu_i},\; \mu_i\in
\mathfrak{g}^\ast_i$, the regular reducible orbits of the RCH systems
$(T^\ast Q_i, G_i, \omega_i, H_i, F_i, W_i)$, $i=1,2$, the map
$\varphi^\ast_{\mathcal{O}_\mu}=i_{\mathcal{O}_{\mu_1}}^{-1}\cdot\varphi^\ast\cdot
i_{\mathcal{O}_{\mu_2}}:\mathbf{J}_2^{-1}(\mathcal{O}_{\mu_2})\rightarrow
\mathbf{J}_1^{-1}(\mathcal{O}_{\mu_1})$ is $(G_2,G_1)$-equivariant,
 and $W_1\cap \mathbf{J}_1^{-1}(\mathcal{O}_{\mu_1})=\varphi_{\mathcal{O}_\mu}^\ast
(W_2\cap \mathbf{J}_2^{-1}(\mathcal{O}_{\mu_2}))$, and
$\mathbf{J}_{2\mathcal{O}_{\mu_2}}^\ast
\omega_{2\mathcal{O}_{\mu_2}}^{+}=(\varphi_{\mathcal{O}_\mu}^\ast)^\ast
\cdot\mathbf{J}_{1\mathcal{O}_{\mu_1}}^\ast\omega_{1\mathcal{O}_{\mu_1}}^{+},$
where $\mu=(\mu_1, \mu_2)$, and denote by
$i_{\mathcal{O}_{\mu_1}}^{-1}(S)$ the pre-image of a subset
$S\subset T^\ast Q_1$ for the map
$i_{\mathcal{O}_{\mu_1}}:\mathbf{J}_1^{-1}(\mathcal{O}_{\mu_1})\rightarrow
T^\ast Q_1$.

\noindent {\bf RoCH-2:}
For each control law $u_1:
T^\ast Q_1 \rightarrow W_1, $ there exists the control law $u_2:
T^\ast Q_2 \rightarrow W_2 $  such that the two
closed-loop dynamical systems produce the same dynamical vector fields; that is,
$X_{(T^\ast Q_1,G_1,\omega_1,H_1,F_1,u_1)}\cdot \varphi^\ast
=T(\varphi^\ast) X_{(T^\ast Q_2,G_2,\omega_2,H_2,F_2,u_2)}$.
\end{defi}

It is worthy of noting that for the regular orbit reducible RCH
system, the induced equivalent map $\varphi^*$ not only keeps the restriction
of the $(+)$-symplectic structure on the regular reducible orbit to
$\mathbf{J}^{-1}(\mathcal{O}_\mu)$, but also keeps the equivariance
of $G$-action on the regular reducible orbit. If an $R_o$-reduced feedback control law
$u_{\mathcal{O}_\mu}:(T^\ast Q)_{\mathcal{O}_\mu}\rightarrow
W_{\mathcal{O}_\mu}$ is chosen, the $R_o$-reduced RCH system
$((T^\ast Q)_{\mathcal{O}_\mu},
\omega_{\mathcal{O}_\mu},h_{\mathcal{O}_\mu},f_{\mathcal{O}_\mu},u_{\mathcal{O}_\mu})$
is a closed-loop regular dynamic system with the control law
$u_{\mathcal{O}_\mu}$. Assume that its vector field $X_{((T^\ast
Q)_{\mathcal{O}_\mu}, \omega_{\mathcal{O}_\mu},
h_{\mathcal{O}_\mu},f_{\mathcal{O}_\mu},u_{\mathcal{O}_\mu})}$ can
be expressed by
\begin{equation}X_{((T^\ast Q)_{\mathcal{O}_\mu},
\omega_{\mathcal{O}_\mu},h_{\mathcal{O}_\mu},f_{\mathcal{O}_\mu},u_{\mathcal{O}_\mu})}=
X_{h_{\mathcal{O}_\mu}}+\textnormal{vlift}(f_{\mathcal{O}_\mu})
+\textnormal{vlift}(u_{\mathcal{O}_\mu}), \label{4.8}
\end{equation}
where $X_{h_{\mathcal{O}_\mu}}$ is the dynamical vector field of
the $R_o$-reduced Hamiltonian $h_{\mathcal{O}_\mu}$,
and $\textnormal{vlift}(f_{\mathcal{O}_\mu})=
\textnormal{vlift}(f_{\mathcal{O}_\mu})X_{h_{\mathcal{O}_\mu}}$,
$\textnormal{vlift}(u_{\mathcal{O}_\mu})=
\textnormal{vlift}(u_{\mathcal{O}_\mu})X_{h_{\mathcal{O}_\mu}}$,
are the changes of $X_{h_{\mathcal{O}_\mu}}$
under the actions of the $R_o$-reduced external force $f_{\mathcal{O}_\mu}$
and the $R_o$-reduced control law $u_{\mathcal{O}_\mu}$.
The dynamical vector fields of
the regular orbit reducible RCH system $(T^\ast Q,G,\omega,
H,F,u)$ and the $R_o$-reduced RCH system $((T^\ast Q)_{\mathcal{O}_\mu},
\omega_{\mathcal{O}_\mu}, h_{\mathcal{O}_\mu},
f_{\mathcal{O}_\mu}, u_{\mathcal{O}_\mu})$ satisfy the condition that
\begin{equation}X_{((T^\ast Q)_{\mathcal{O}_\mu},
\omega_{\mathcal{O}_\mu},h_{\mathcal{O}_\mu},f_{\mathcal{O}_\mu},u_{\mathcal{O}_\mu})}\cdot
\pi_{\mathcal{O}_\mu} =T\pi_{\mathcal{O}_\mu} \cdot X_{(T^\ast
Q,G,\omega,H,F,u)}\cdot
i_{\mathcal{O}_\mu}.
\label{4.9}\end{equation}
Thus, we can prove the
following regular orbit reduction theorem for the RCH system, which
explains the relationship between the RoCH-equivalence for the regular
orbit reducible RCH systems with symmetries and momentum maps
and the RCH-equivalence
for the associated $R_o$-reduced RCH systems,
its proof was given in Marsden \textit{et al.} \cite{mawazh10} and
Wang \cite{wa18}.
This theorem can be
regarded as an extension of the regular orbit reduction theorem of
Hamiltonian system under the regular controlled Hamiltonian equivalence
conditions.
\begin{theo}
Two regular orbit reducible RCH systems $(T^\ast Q_i, G_i,
\omega_i, H_i, F_i,W_i)$, $i=1,2,$ are RoCH-equivalent if and only
if the associated $R_o$-reduced RCH systems $((T^\ast
Q)_{\mathcal{O}_{\mu_i}}, \omega_{i\mathcal{O}_{\mu_i}},
h_{i\mathcal{O}_{\mu_i}}, f_{i\mathcal{O}_{\mu_i}},\\
W_{i\mathcal{O}_{\mu_i}})$, $i=1,2,$ are RCH-equivalent.
\end{theo}

It is worthy of noting that when the external force and the control of a
regular orbit reducible RCH system $(T^*Q,G,\omega,H,F,W)$ are both zeros;
that is, $F=0 $ and $W=\emptyset$, in this case the RCH system
is just a regular orbit reducible Hamiltonian system $(T^*Q,G,\omega,H)$ itself.
Then the following theorem explains the relationship between the equivalence
for the regular orbit reducible Hamiltonian systems with symmetries and the
equivalence for the associated $R_o$-reduced Hamiltonian systems.

\begin{theo}
If two regular orbit reducible Hamiltonian systems $(T^\ast Q_i, G_i,
\omega_i, H_i)$, $i=1,2,$ are equivalent, then their
associated $R_o$-reduced Hamiltonian systems $((T^\ast
Q)_{\mathcal{O}_{\mu_i}}, \omega_{i\mathcal{O}_{\mu_i}},
h_{i\mathcal{O}_{\mu_i}})$, $i=1,2,$ must be equivalent.
Conversely, if the $R_o$-reduced Hamiltonian systems $((T^\ast
Q)_{\mathcal{O}_{\mu_i}}, \omega_{i\mathcal{O}_{\mu_i}},
h_{i\mathcal{O}_{\mu_i}})$, $i=1,2,$ are equivalent and the
induced map
$\varphi^\ast_{\mathcal{O}_\mu}:\mathbf{J}_2^{-1}(\mathcal{O}_{\mu_2})\rightarrow
\mathbf{J}_1^{-1}(\mathcal{O}_{\mu_1})$ such that
$\mathbf{J}_{2\mathcal{O}_{\mu_2}}^\ast
\omega_{2\mathcal{O}_{\mu_2}}^{+}=(\varphi_{\mathcal{O}_\mu}^\ast)^\ast
\cdot\mathbf{J}_{1\mathcal{O}_{\mu_1}}^\ast\omega_{1\mathcal{O}_{\mu_1}}^{+},$
then the regular orbit reducible Hamiltonian systems $(T^\ast Q_i, G_i,
\omega_i, H_i)$, $i=1,2,$ are equivalent.
\end{theo}
See the proof and the more details in Wang \cite{wa18}.

\begin{rema}
If $(T^\ast Q, \omega)$ is a connected symplectic manifold, and
$\mathbf{J}:T^\ast Q\rightarrow \mathfrak{g}^\ast$ is a
non-equivariant momentum map with a non-equivariance group
one-cocycle $\sigma: G\rightarrow \mathfrak{g}^\ast$, which is
defined by $\sigma(g):=\mathbf{J}(g\cdot
z)-\operatorname{Ad}^\ast_{g^{-1}}\mathbf{J}(z)$, where $g\in G$ and
$z\in T^\ast Q$. Then we know that $\sigma$ produces a new affine
action $\Theta: G\times \mathfrak{g}^\ast \rightarrow
\mathfrak{g}^\ast $ defined by
$\Theta(g,\mu):=\operatorname{Ad}^\ast_{g^{-1}}\mu + \sigma(g)$,
where $\mu \in \mathfrak{g}^\ast$, with respect to which the given
momentum map $\mathbf{J}$ is equivariant. Assume that the $G$ acts
freely and properly on $T^\ast Q$, and that $\tilde{G}_\mu$ is the
isotropy subgroup of $\mu \in \mathfrak{g}^\ast$ relative to this
affine action $\Theta$, and that $\mathcal{O}_\mu= G\cdot \mu
\subset \mathfrak{g}^\ast$ is the G-orbit of the point $\mu$
with respect to the action $\Theta$,
and that $\mu$ is a regular value of $\mathbf{J}$.
Then the quotient space $(T^\ast
Q)_\mu=\mathbf{J}^{-1}(\mu)/\tilde{G}_\mu$ is a symplectic
manifold with the symplectic form $\omega_\mu$ uniquely characterized by
$(4.1)$; and the quotient space
$(T^\ast Q)_{\mathcal{O}_\mu}=\mathbf{J}^{-1}(\mathcal{O}_\mu)/ G $
is also a symplectic manifold with the symplectic form
$\omega_{\mathcal{O}_\mu}$ uniquely characterized by $(4.5)$;
see Ortega and Ratiu \cite{orra04}.
Moreover, in this case,
for the given regular point or regular orbit reducible RCH system
$(T^*Q,G,\omega,H,F,W)$, we can also prove the regular point reduction theorem
or the regular orbit reduction theorem, by using the above similar ways.
\end{rema}

\section{Hamilton-Jacobi Equations for the RCH System and Its Reduced Systems}

We know that under the actions of the external force $F$
and the control $u$, in general, the RCH system
$(T^*Q,\omega,H,F,u)$, the $R_p$-reduced RCH
system $((T^\ast Q)_\mu, \omega_\mu, h_\mu, f_\mu, u_\mu)$
and the $R_o$-reduced RCH system $((T^\ast Q)_{\mathcal{O}_\mu},
\omega_{\mathcal{O}_\mu}, h_{\mathcal{O}_\mu},
f_{\mathcal{O}_\mu}, u_{\mathcal{O}_\mu})$, all of them
are not Hamiltonian systems, and hence, we cannot describe
the Hamilton-Jacobi equations for the
RCH system, the $R_p$-reduced RCH
system and the $R_o$-reduced RCH
system, from the viewpoint of a generating
function the same as in Theorem 2.1
(given by Abraham and Marsden in \cite{abma78}).
However, for a given RCH system $(T^*Q,\omega,H,F,W)$ in which $\omega$ is the
canonical symplectic form on $T^*Q$, by using
Lemma 2.3 and the expression (3.1) of the dynamical
vector field of the RCH system,
we can derive precisely two types of geometric constraint conditions of
the canonical symplectic form for the dynamical vector field of the RCH system;
that is, the two types of Hamilton-Jacobi equations for the RCH system.
Moreover, we generalize the above results for a regular
reducible RCH system with symmetry and a momentum map,
and derive precisely the two types of Hamilton-Jacobi equations
for the $R_p$-reduced RCH system and the $R_o$-reduced RCH system,
and we prove that the
RCH-equivalence for the RCH system, and the RpCH-equivalence and
the RoCH-equivalence for the regular reducible RCH systems with symmetries
and momentum maps, leave the
solutions of the corresponding Hamilton-Jacobi equations invariant.
These research reveal the deeply internal
relationships of the geometrical structures of phase spaces, the dynamical
vector fields and the controls of the RCH system and the reduced
RCH system.

\subsection{Hamilton-Jacobi Equations for an RCH System }

For a given RCH system $(T^*Q,\omega,H,F,W)$ in which $\omega$ is the
canonical symplectic form on $T^*Q$, by using
Lemma 2.3 and the expression (3.1) of the dynamical vector field $X_{(T^\ast Q,\omega,H,F,u)}$,
we can derive precisely two types of geometric constraint conditions of
the canonical symplectic form for the dynamical vector field of the RCH system;
that is, the two types of Hamilton-Jacobi equations for the RCH system.
For convenience, the maps involved in the theorem
and its proof are shown in Diagram-2.

\begin{center}
\hskip 0cm \xymatrix{ T^* Q \ar[r]^\varepsilon
& T^* Q \ar[d]_{X_{H\cdot \varepsilon}} \ar[dr]^{\tilde{X}^\varepsilon} \ar[r]^{\pi_Q}
& Q \ar[d]^{\tilde{X}^\gamma} \ar[r]^{\gamma} & T^*Q \ar[d]^{\tilde{X}} \\
  & T(T^*Q) & TQ \ar[l]^{T\gamma} & T(T^* Q)\ar[l]^{T\pi_Q}}
\end{center}
$$\mbox{Diagram-2}$$

\begin{theo}
(Hamilton-Jacobi Theorem for an RCH system)
For the RCH system $(T^*Q,\omega,H,F,W)$ with the
canonical symplectic form $\omega$ on $T^*Q$, assume that $\gamma: Q
\rightarrow T^*Q$ is a one-form on $Q$, and that
$\lambda=\gamma \cdot \pi_{Q}: T^* Q \rightarrow T^* Q $, and that for any
symplectic map $\varepsilon: T^* Q \rightarrow T^* Q $, denote that
$\tilde{X}^\gamma = T\pi_{Q}\cdot \tilde{X} \cdot \gamma$ and
$\tilde{X}^\varepsilon = T\pi_{Q}\cdot \tilde{X} \cdot \varepsilon$,
where $\tilde{X}=X_{(T^\ast Q,\omega,H,F,u)}$ is the dynamical vector field of the RCH system
$(T^*Q,\omega,H,F,W)$ with a control law $u$.
Then the following two assertions hold:\\
\noindent $(\mathbf{i})$
If the one-form $\gamma: Q \rightarrow T^*Q $ is closed with respect to
$T\pi_Q: TT^* Q \rightarrow TQ, $ then $\gamma$ is a solution of the Type I
Hamilton-Jacobi equation
$T\gamma\cdot \tilde{X}^\gamma= X_H\cdot \gamma ,$
where $X_H$ is the Hamiltonian vector field
of the corresponding Hamiltonian system $(T^*Q,\omega,H).$ \\
\noindent $(\mathbf{ii})$
The $\varepsilon$ is a solution of the Type II
Hamilton-Jacobi equation $T\gamma\cdot \tilde{X}^\varepsilon= X_H\cdot
\varepsilon $ if and only if it is a solution of the equation
$T\varepsilon\cdot X_{H\cdot\varepsilon}= T\lambda \cdot \tilde{X} \cdot \varepsilon,$
where $X_H$ and $ X_{H\cdot\varepsilon} \in
TT^*Q $ are the Hamiltonian vector fields of the functions $H$ and $H\cdot\varepsilon:
T^*Q\rightarrow \mathbb{R}, $ respectively.
\end{theo}
See the proof and the more details in Wang \cite{wa13d}.
It is worth noting that
the Type I Hamilton-Jacobi equation
$T\gamma\cdot \tilde{X}^\gamma= X_H \cdot \gamma ,$
is the equation of the differential one-form $\gamma$, and that
the Type II Hamilton-Jacobi equation $T\gamma\cdot \tilde{X}^\varepsilon
= X_H \cdot \varepsilon ,$ is the equation of
the symplectic diffeomorphism map $\varepsilon$.\\

Moreover, considering the RCH-equivalence of the RCH systems,
we can prove the following theorem, which states that
the solutions of two types of Hamilton-Jacobi equations for the RCH systems remain
invariant under the conditions of RCH-equivalence if the corresponding Hamiltonian
systems are equivalent:

\begin{theo}
Suppose that two RCH systems, $(T^\ast Q_i,\omega_i,H_i,F_i,W_i)$,
$i=1,2 $ are RCH-equivalent with an equivalent map $\varphi: Q_1
\rightarrow Q_2 $, and that the corresponding Hamiltonian systems, $(T^\ast
Q_i,\omega_i,H_i),$ $ i= 1,2 $ are also equivalent.
Under the hypotheses and notations of Theorems 5.1,
we have that the following two assertions hold:\\
\noindent $(\mathrm{i})$ If the one-form $\gamma_2: Q_2 \rightarrow T^* Q_2$ is closed with
respect to $T\pi_{Q_2}: TT^* Q_2 \rightarrow TQ_2, $ then $\gamma_1=
\varphi^* \cdot \gamma_2\cdot \varphi: Q_1 \rightarrow T^* Q_1 $ is also closed with respect to $T\pi_{Q_1}:
TT^* Q_1 \rightarrow TQ_1, $ and hence, $\gamma_1$ is
a solution of the Type I Hamilton-Jacobi equation for the RCH system
$(T^*Q_1,\omega_1,H_1,F_1,W_1)$.

\noindent $(\mathrm{ii})$ If the symplectic map $\varepsilon_2: T^*Q_2\rightarrow T^* Q_2$ is
a solution of the Type II Hamilton-Jacobi equation for the RCH system
$(T^*Q_2,\omega_2,H_2, F_2,W_2)$, then $\varepsilon_1=
\varphi^* \cdot \varepsilon_2\cdot \varphi_*: T^*Q_1 \rightarrow
T^* Q_1 $ is a symplectic map, and hence, $\varepsilon_1$ is a solution
of the Type II Hamilton-Jacobi equation for the RCH
system $(T^*Q_1,\omega_1,H_1,F_1,W_1). $
\end{theo}
See the proof and the more details in Wang \cite{wa13d}.

\subsection{Hamilton-Jacobi Equations for an $R_p$-reduced RCH System }

For a given regular point reducible RCH system
$(T^*Q,G,\omega,H,F,W)$ with the
canonical symplectic form $\omega$ on $T^*Q$ and
with an $R_p$-reduced RCH system $((T^\ast
Q)_\mu, \omega_\mu,h_\mu, f_\mu, u_\mu )$, by using Lemma 2.3,
we can give precisely the geometric constraint
conditions of the $R_p$-reduced symplectic form for the
dynamical vector field of the regular point reducible RCH system;
 that is, the Type I and Type II
Hamilton-Jacobi equations for the $R_p$-reduced RCH system $((T^\ast
Q)_\mu, \omega_\mu,h_\mu, f_\mu, u_\mu )$.
First, by using the fact that the one-form $\gamma: Q
\rightarrow T^*Q $ is closed with respect to
$T\pi_Q: TT^* Q \rightarrow TQ, $ and that $\textmd{Im}(\gamma)\subset
\mathbf{J}^{-1}(\mu), $ and that $\gamma$ is $G_\mu$-invariant,
we can prove the Type I
Hamilton-Jacobi theorem for the $R_p$-reduced RCH system $((T^\ast
Q)_\mu, \omega_\mu,h_\mu, f_\mu, u_\mu )$.
For convenience, the maps involved in
the theorem and its proof are shown in Diagram-3.
\begin{center}
\hskip 0cm \xymatrix{ \mathbf{J}^{-1}(\mu) \ar[r]^{i_\mu} & T^* Q \ar[d]_{X_H} \ar[r]^{\pi_Q}
& Q \ar[d]_{\tilde{X}^\gamma} \ar[r]^{\gamma}
& T^*Q \ar[d]_{\tilde{X}} \ar[r]^{\pi_\mu}
& (T^* Q)_\mu \ar[d]_{X_{h_\mu}} \\
& T(T^*Q)  & TQ \ar[l]^{T\gamma}
& T(T^*Q) \ar[l]^{T\pi_Q} \ar[r]_{T\pi_\mu} & T(T^* Q)_\mu }
\end{center}
$$\mbox{Diagram-3}$$

\begin{theo}
(Type I Hamilton-Jacobi Theorem for an $R_p$-reduced RCH system)
For the regular point reducible RCH system
$(T^*Q,G,\omega,H,F,W)$ with the
canonical symplectic form $\omega$ on $T^*Q$
and with an $R_p$-reduced RCH system
$((T^\ast Q)_\mu, \omega_\mu,h_\mu, f_\mu, u_\mu )$,
assume that $\gamma: Q \rightarrow T^*Q$ is a one-form
on $Q$, and that $\tilde{X}^\gamma = T\pi_{Q}\cdot \tilde{X} \cdot \gamma$,
where $\tilde{X}=X_{(T^\ast Q,G,\omega,H,F,u)}$ is the dynamical vector
field of the regular point reducible RCH system
$(T^*Q,G,\omega,H,F,W)$ with a control law $u$. Moreover,
assume that $\mu \in \mathfrak{g}^\ast $ is a regular value of the momentum
map $\mathbf{J}$, and that $\textmd{Im}(\gamma)\subset
\mathbf{J}^{-1}(\mu), $ and that $\gamma$ is $G_\mu$-invariant, and
that $\bar{\gamma}=\pi_\mu(\gamma): Q \rightarrow (T^* Q)_\mu. $
If the one-form $\gamma: Q \rightarrow T^*Q $ is closed with respect to
$T\pi_Q: TT^* Q \rightarrow TQ, $
then $\bar{\gamma}$ is a solution of the equation
$T\bar{\gamma}\cdot \tilde{X}^\gamma= X_{h_\mu}\cdot \bar{\gamma}, $
where $X_{h_\mu}$ is the Hamiltonian vector field of the $R_p$-reduced
Hamiltonian function $h_\mu:(T^\ast Q)_\mu\rightarrow \mathbb{R}, $
and the equation is called the Type I Hamilton-Jacobi equation for
the $R_p$-reduced RCH System $((T^\ast Q)_\mu, \omega_\mu,h_\mu, f_\mu, u_\mu)$.
\end{theo}
See the proof and the more details in Wang \cite{wa13d}. Moreover,
for any $G_\mu$-invariant symplectic map $\varepsilon: T^* Q \rightarrow T^* Q $,
we can prove the Type II
Hamilton-Jacobi theorem for the $R_p$-reduced RCH system $((T^\ast
Q)_\mu, \omega_\mu,h_\mu, f_\mu, u_\mu )$; see Wang \cite{wa13d}.\\

\begin{rema}
It is worth noting that
the Type I Hamilton-Jacobi equation
$T\bar{\gamma}\cdot \tilde{X}^\gamma= X_{h_\mu}\cdot
\bar{\gamma} $ is the equation of the
$R_p$-reduced differential one-form $\bar{\gamma}$, and that
the Type II Hamilton-Jacobi equation $T\bar{\gamma}\cdot \tilde{X}^\varepsilon
= X_{h_\mu}\cdot \bar{\varepsilon} $ is the equation of the symplectic
diffeomorphism map $\varepsilon$ and the $R_p$-reduced
symplectic diffeomorphism map $\bar{\varepsilon}. $
If both the external force and the control of the regular
point reducible RCH system $(T^*Q,G,\omega,H,F,W)$ are zero;
that is, $F=0 $ and $W=\emptyset$, in this
case, the RCH system is just a regular point reducible Hamiltonian
system $(T^*Q,G,\omega,H)$ itself. From the proofs of
the Theorem 5.3 above, we can also get the Theorem 2.4; that is,
the Type I Hamilton-Jacobi equation
for the associated Marsden-Weinstein reduced Hamiltonian system ( given in
Wang \cite{wa17}). Thus, Theorem 5.3 can be regarded as an extension
of the Type I Hamilton-Jacobi equation for the Marsden-Weinstein reduced Hamiltonian
system to that for the $R_p$-reduced RCH system.
\end{rema}

Moreover, considering the RpCH-equivalence of the regular point
reducible RCH systems and using
Theorems 4.3, 4.4 and 5.2, we can obtain the Theorem 5.5,
which states that the solutions of the two types of Hamilton-Jacobi equations for
the regular point reducible RCH systems remain invariant
under the conditions of RpCH-equivalence if the corresponding Hamiltonian
systems are equivalent.
\begin{theo}
Suppose that two regular point reducible RCH systems,
$(T^\ast Q_i,G_i,\omega_i, H_i, F_i, W_i)$,
$i=1,2 $ are RpCH-equivalent with an equivalent map $\varphi: Q_1
\rightarrow Q_2 $, and that the associated $R_p$-reduced RCH systems
are $((T^\ast Q_i)_{\mu_i},\omega_{i\mu_i},h_{i\mu_i},f_{i\mu_i}, u_{i\mu_i}),$
$i=1,2.$ Assume that the corresponding Hamiltonian systems, $(T^\ast
Q_i,G_i,\omega_i,H_i),$ $ i= 1,2 $ are also equivalent. Then,
under the hypotheses and the notations of Theorems 4.3,
4.4 and 5.3, we have that the following two assertions hold:\\
\noindent $(\mathrm{i})$ If the one-form
$\gamma_2: Q_2 \rightarrow T^* Q_2$ is closed with
respect to $T\pi_{Q_2}: TT^* Q_2 \rightarrow TQ_2, $ and
$\bar{\gamma}_2=\pi_{2\mu_2}(\gamma_2): Q_2 \rightarrow (T^* Q_2)_{\mu_2} $
is a solution of the Type I Hamilton-Jacobi equation for
the $R_p$-reduced RCH system $((T^\ast Q_2)_{\mu_2},
\omega_{2\mu_2},h_{2\mu_2}, f_{2\mu_2}, u_{2\mu_2})$,
then $\gamma_1=\varphi^* \cdot \gamma_2\cdot \varphi: Q_1 \rightarrow T^* Q_1 $ is
a solution of the Type I Hamilton-Jacobi equation for the
RCH system $(T^*Q_1,G_1,\omega_1,H_1,F_1,W_1) $,
and $\bar{\gamma}_1=\pi_{1\mu_1}(\gamma_1): Q_1 \rightarrow (T^* Q_1)_{\mu_1}$ is
a solution of the Type I Hamilton-Jacobi equation for the $R_p$-reduced
RCH system $((T^\ast Q_1)_{\mu_1}, \omega_{1\mu_1},
h_{1\mu_1}, f_{1\mu_1}, u_{1\mu_1})$.

\noindent $(\mathrm{ii})$ If the $G_{2\mu_2}$-invariant
symplectic map $\varepsilon_2: T^*Q_2\rightarrow T^* Q_2$ and
$\bar{\varepsilon}_2=\pi_{2\mu_2}(\varepsilon_2):
\mathbf{J_2}^{-1}(\mu_2) (\subset T^*Q_2) \rightarrow (T^* Q_2)_{\mu_2} $
satisfy the Type II Hamilton-Jacobi equation for the $R_p$-reduced RCH system
$((T^\ast Q_2)_{\mu_2}, \omega_{2\mu_2},h_{2\mu_2}, f_{2\mu_2}, u_{2\mu_2})$,
then $\varepsilon_1=
\varphi^* \cdot \varepsilon_2\cdot \varphi_*: T^*Q_1 \rightarrow T^* Q_1 $ and
$\bar{\varepsilon}_1=\pi_{1\mu_1}(\varepsilon_1):
\mathbf{J_1}^{-1}(\mu_1) (\subset T^*Q_1) \rightarrow (T^* Q_1)_{\mu_1} $
satisfy the Type II Hamilton-Jacobi equation for the $R_p$-reduced
RCH system $((T^\ast Q_1)_{\mu_1}, \omega_{1\mu_1},h_{1\mu_1}, f_{1\mu_1}, u_{1\mu_1}). $
\end{theo}
See the proof and the more details in Wang \cite{wa13d}.

\subsection{Hamilton-Jacobi Equations for an $R_o$-reduced RCH System }

For a given regular orbit reducible RCH system
$(T^*Q,G,\omega,H,F,W)$ with the
canonical symplectic form $\omega$ on $T^*Q$
and with an $R_o$-reduced RCH system $((T^\ast
Q)_{\mathcal{O}_\mu}, \omega_{\mathcal{O}_\mu},h_{\mathcal{O}_\mu},
f_{\mathcal{O}_\mu},u_{\mathcal{O}_\mu})$,
by using Lemma 2.3, we can give precisely the geometric constraint
conditions of the $R_o$-reduced symplectic form for the
dynamical vector field of the regular orbit reducible RCH system;
that is, the Type I and Type II
Hamilton-Jacobi equations for the $R_o$-reduced RCH system
$((T^\ast Q)_{\mathcal{O}_\mu}, \omega_{\mathcal{O}_\mu},h_{\mathcal{O}_\mu},
f_{\mathcal{O}_\mu},u_{\mathcal{O}_\mu})$.\\

In the following, by using Lemma 2.3 and the $R_o$-reduced symplectic
form $\omega_{\mathcal{O}_\mu}$,
for any $G$-invariant symplectic map
$\varepsilon: T^* Q \rightarrow T^* Q $,
we can derive precisely the Type II
Hamilton-Jacobi equation for the $R_o$-reduced
RCH system $((T^\ast Q)_{\mathcal{O}_\mu}, \omega_{\mathcal{O}_\mu},h_{\mathcal{O}_\mu},
f_{\mathcal{O}_\mu},u_{\mathcal{O}_\mu})$ as follows:

\begin{theo}
(Type II Hamilton-Jacobi Theorem for an $R_o$-reduced RCH system)
For a given regular orbit reducible RCH system
$(T^*Q,G,\omega,H,F,W)$ with the
canonical symplectic form $\omega$ on $T^*Q$
and with an $R_o$-reduced RCH system
$((T^\ast Q)_{\mathcal{O}_\mu}, \omega_{\mathcal{O}_\mu},h_{\mathcal{O}_\mu},
f_{\mathcal{O}_\mu},u_{\mathcal{O}_\mu})$,
assume that $\gamma: Q \rightarrow T^*Q$ is a one-form on $Q$,
and that $\lambda=\gamma \cdot \pi_{Q}: T^* Q
\rightarrow T^* Q $. For any $G$-invariant symplectic map
$\varepsilon: T^* Q \rightarrow T^* Q $,
denote that $\tilde{X}^\varepsilon = T\pi_{Q}\cdot
\tilde{X} \cdot \varepsilon$, where
$\tilde{X}=X_{(T^\ast Q,G,\omega,H,F,u)}$ is the dynamical vector
field of the regular orbit reducible RCH system
$(T^*Q,G,\omega,H,F,W)$ with a control law $u$.
Moreover, assume that $\mu \in \mathfrak{g}^\ast $ is a regular value of the momentum
map $\mathbf{J}$, and that $\mathcal{O}_\mu, \; (\mu \in \mathfrak{g}^\ast) $ is the
regular reducible orbit of the corresponding Hamiltonian system $(T^*Q,G,\omega,H)$, and that
$\textmd{Im}(\gamma)\subset \mathbf{J}^{-1}(\mu), $ and that $\gamma$ is
$G$-invariant, and that $\varepsilon(\mathbf{J}^{-1}(\mathcal{O}_\mu))
\subset \mathbf{J}^{-1}(\mathcal{O}_\mu). $
Denote that $\bar{\gamma}=\pi_{\mathcal{O}_\mu}(\gamma): Q
\rightarrow (T^* Q)_{\mathcal{O}_\mu} $,
$\bar{\lambda}=\pi_{\mathcal{O}_\mu}(\lambda):
\mathbf{J}^{-1}(\mathcal{O}_\mu) (\subset T^*Q) \rightarrow (T^*
Q)_{\mathcal{O}_\mu} $, and $\bar{\varepsilon}=\pi_{\mathcal{O}_\mu}(\varepsilon):
\mathbf{J}^{-1}(\mathcal{O}_\mu) (\subset T^*Q) \rightarrow (T^*
Q)_{\mathcal{O}_\mu} $. Then $\varepsilon$ and $\bar{\varepsilon}$ satisfy the equation
$T\bar{\varepsilon}\cdot X_{h_{\mathcal{O}_\mu}\cdot\bar{\varepsilon}}
= T\bar{\lambda} \cdot \tilde{X}\cdot\varepsilon $
if and only if they satisfy the equation $T\bar{\gamma}\cdot \tilde{X}^\varepsilon=
X_{h_{\mathcal{O}_\mu}}\cdot \bar{\varepsilon}$,
where $X_{h_{\mathcal{O}_\mu}}$ and
$ X_{h_{\mathcal{O}_\mu} \cdot \bar{\varepsilon}} \in TT^*Q $
are the Hamiltonian vector fields of the $R_o$-reduced Hamiltonian
functions $h_{\mathcal{O}_\mu}$ and
$h_{\mathcal{O}_\mu} \cdot \bar{\varepsilon}: T^*Q\rightarrow
\mathbb{R}, $ respectively.
The equation $T\bar{\gamma}\cdot \tilde{X}^\varepsilon=
X_{h_{\mathcal{O}_\mu}}\cdot \bar{\varepsilon}$
is called the Type II Hamilton-Jacobi
equation for the $R_o$-reduced RCH system $((T^\ast
Q)_{\mathcal{O}_\mu}, \omega_{\mathcal{O}_\mu},h_{\mathcal{O}_\mu},
f_{\mathcal{O}_\mu},u_{\mathcal{O}_\mu})$.
Here the maps involved in the theorem are shown in Diagram-4.

\begin{center}
\hskip 0cm \xymatrix{ \mathbf{J}^{-1}(\mathcal{O}_\mu)
\ar[r]^{i_{\mathcal{O}_\mu}} & T^* Q \ar[d]_{X_{H\cdot \varepsilon}}
\ar[dr]^{\tilde{X}^\varepsilon} \ar[r]^{\pi_Q}
& Q \ar[r]^{\gamma} & T^*Q \ar[d]_{\tilde{X}}
\ar[dr]^{X_{h_{\mathcal{O}_\mu} \cdot\bar{\varepsilon}}} \ar[r]^{\pi_{\mathcal{O}_\mu}}
& (T^* Q)_{\mathcal{O}_\mu} \ar[d]^{X_{h_{\mathcal{O}_\mu}}} \\
  & T(T^*Q)  & TQ \ar[l]^{T\gamma} & T(T^*Q) \ar[l]^{T\pi_Q}
   \ar[r]_{T\pi_{\mathcal{O}_\mu}} & T(T^* Q)_{\mathcal{O}_\mu} }
\end{center}
$$\mbox{Diagram-4}$$
\end{theo}
See the proof and the more details in Wang \cite{wa13d}.\\

Moreover, for the regular orbit reducible RCH system
$(T^*Q,G,\omega,H,F,W)$ with an $R_o$-reduced RCH system $((T^\ast
Q)_{\mathcal{O}_\mu}, \omega_{\mathcal{O}_\mu},h_{\mathcal{O}_\mu},f_{\mathcal{O}_\mu},u_{\mathcal{O}_\mu})$,
we know that the Hamiltonian vector fields
$X_{H}$ and $X_{h_{\mathcal{O}_\mu}}$ for
the corresponding Hamiltonian system $(T^*Q,G,\omega,H)$
and its $R_o$-reduced system $((T^\ast Q)_{\mathcal{O}_\mu}, \\
\omega_{\mathcal{O}_\mu},h_{\mathcal{O}_\mu})$ are
$\pi_{\mathcal{O}_\mu}$-related; that is, $X_{
h_{\mathcal{O}_\mu}}\cdot \pi_{\mathcal{O}_\mu} =
T\pi_{\mathcal{O}_\mu}\cdot X_{H}\cdot i_{\mathcal{O}_\mu}. $ Then
we can prove the Theorem 5.7, which states the
relationship between the solutions of Type II Hamilton-Jacobi
equations and the regular orbit reduction:

\begin{theo}
For the regular orbit reducible RCH system
$(T^*Q,G,\omega,H,F,W)$ with the
canonical symplectic form $\omega$ on $T^*Q$
and with an $R_o$-reduced RCH system $((T^\ast
Q)_{\mathcal{O}_\mu}, \omega_{\mathcal{O}_\mu},
h_{\mathcal{O}_\mu},f_{\mathcal{O}_\mu},u_{\mathcal{O}_\mu})$,
assume that $\gamma: Q \rightarrow T^*Q$ is a
one-form on $Q$, and that $\varepsilon: T^* Q \rightarrow T^* Q $ is
a $G$-invariant symplectic map,
$\bar{\varepsilon}=\pi_{\mathcal{O}_\mu}(\varepsilon):
\mathbf{J}^{-1}(\mathcal{O}_\mu) (\subset T^*Q) \rightarrow (T^*
Q)_{\mathcal{O}_\mu} $.
Under the hypotheses and notations of Theorem 5.6,
then we have that $\varepsilon$ is a solution of the Type
II Hamilton-Jacobi equation $T\gamma\cdot \tilde{X}^\varepsilon=
X_H\cdot \varepsilon $ for the regular orbit reducible RCH
system $(T^*Q,G,\omega,H,F,W) $ if and only if $\varepsilon$ and
$\bar{\varepsilon}$ satisfy the Type II Hamilton-Jacobi equation
$T\bar{\gamma}\cdot \tilde{X}^\varepsilon=
X_{h_{\mathcal{O}_\mu}\cdot\bar{\varepsilon}} $ for the
$R_o$-reduced RCH system $((T^\ast Q)_{\mathcal{O}_\mu},
\omega_{\mathcal{O}_\mu},h_{\mathcal{O}_\mu},
f_{\mathcal{O}_\mu},u_{\mathcal{O}_\mu})$.
\end{theo}
See the proof and the more details in Wang \cite{wa13d}.\\

It is worth noting that the different symplectic forms determine
the different regular reduced RCH systems. From (4.5) we know that,
for the regular orbit reduced symplectic space
$(T^\ast Q)_{\mathcal{O}_\mu}= \mathbf{J}^{-1}(\mathcal{O}_\mu)/G
\cong \mathbf{J}^{-1}(\mu)/G \times \mathcal{O}_\mu, $ if we give a stronger
assumption condition; that is, for the one-form $\gamma: Q \rightarrow T^*Q$ on $Q,$
we assume that $\textmd{Im}(\gamma)\subset
\mathbf{J}^{-1}(\mu), $ (note that it is not $\textmd{Im}(\gamma)\subset
\mathbf{J}^{-1}(\mathcal{O}_\mu) $),
and that $\gamma$ is $G$-invariant, then for any $V\in TQ $ and $w\in TT^*Q, $
we have that
$\mathbf{J}_{\mathcal{O}_\mu}^\ast\omega_{\mathcal{O}_\mu}^+(T\gamma
\cdot V, \; w)=0, $ and hence, from (4.5), $i_{\mathcal{O}_\mu}^\ast
\omega=\pi_{\mathcal{O}_{\mu}}^\ast \omega_{\mathcal{O}
_\mu}+\mathbf{J}_{\mathcal{O}_\mu}^\ast\omega_{\mathcal{O}_\mu}^+, $
we have that $\pi_{\mathcal{O}_\mu}^*\omega_{\mathcal{O}_\mu}=
i_{\mathcal{O}_\mu}^*\omega= \omega $ along $\textmd{Im}(\gamma)$.
Thus, we can use Lemma 2.3 for the regular orbit reduced symplectic form $\omega_{\mathcal{O}_\mu}$
in the proofs of Theorems 5.6 and 5.7. It is easy to give the wrong results
without the precise analysis for the regular orbit reduction case.

\section{Controlled Hamiltonian Systems with Nonholonomic Constraints}

In mechanics, it happens very often that systems have constraints.
A nonholonomic Hamiltonian system is a Hamiltonian system with nonholonomic constraint,
and a nonholonomic RCH system is also a RCH system with nonholonomic constraint.
Usually, under the restrictions given by the nonholonomic constraints,
in general, the dynamical vector fields of the nonholonomic Hamiltonian system
and the nonholonomic RCH system
may not be Hamiltonian. Thus, we cannot
describe the Hamilton-Jacobi equations for
the nonholonomic Hamiltonian system and the nonholonomic RCH system
from the viewpoint of a generating function as in Theorem 2.1.
Since the Hamilton-Jacobi theory is developed based on the
Hamiltonian picture of the dynamics, it is a natural idea to extend the
Hamilton-Jacobi theory to the nonholonomic Hamiltonian system
and the nonholonomic RCH system,
and to do so with symmetries and momentum maps;
( see Le\'{o}n and Wang \cite{lewa15}, and Wang \cite{wa22a}).\\

In order to describe the nonholonomic Hamiltonian
system and the nonholonomic RCH system,
in the following we first give the completeness and the regularity
conditions for nonholonomic constraint of a mechanical system;
see Le\'{o}n and Wang in \cite{lewa15}. In fact,
in order to describe the dynamics of a nonholonomic mechanical system,
we need some restriction conditions for the nonholonomic constraints of
the system. First, we note that the set of Hamiltonian vector fields
forms a Lie algebra with respect to the Lie bracket, since
$X_{\{f,g\}}=-[X_f, X_g]. $ However, the Lie bracket operator, in
general, may not be closed on the restriction of a nonholonomic
constraint. Thus, we have to give the completeness
condition for the nonholonomic constraints of a system.\\

{\bf $\mathcal{D}$-completeness } Let $Q$ be a smooth manifold and
$TQ$ its tangent bundle. A distribution $\mathcal{D} \subset TQ$ is
said to be {\bf completely nonholonomic} (or bracket-generating) if
$\mathcal{D}$ along with all of its iterated Lie brackets that
$[\mathcal{D},\mathcal{D}], [\mathcal{D}, [\mathcal{D},\mathcal{D}]],
\cdots ,$ spans the tangent bundle $TQ$. Moreover, we consider a
nonholonomic mechanical system on $Q$, which is
given by a Lagrangian function $L: TQ \rightarrow \mathbb{R}$
subjects to constraint determined by a nonholonomic
distribution $\mathcal{D}\subset TQ$ on the configuration manifold $Q$.
Then the nonholonomic system is said to be {\bf completely nonholonomic},
if the distribution $\mathcal{D} \subset TQ$ determined by the nonholonomic
constraint is completely nonholonomic.\\

{\bf $\mathcal{D}$-regularity } In the following we always assume
that $Q$ is an $n$-dimensional smooth manifold with coordinates $(q^i)$,
and that $TQ$ its tangent bundle with coordinates $(q^i,\dot{q}^i)$,
and that $T^\ast Q$ its cotangent bundle with coordinates $(q^i,p_j)$,
which are the canonical cotangent coordinates of $T^\ast Q$,
and that $\omega= dq^{i}\wedge dp_{i}$ is canonical symplectic form on $T^{\ast}Q$.
If the Lagrangian $L: TQ \rightarrow \mathbb{R}$ is hyperregular;
that is, the Hessian matrix
$(\partial^2L/\partial\dot{q}^i\partial\dot{q}^j)$ is nondegenerate
everywhere, then the Legendre transformation $FL: TQ \rightarrow T^*
Q$ is a diffeomorphism. In this case the Hamiltonian $H: T^* Q
\rightarrow \mathbb{R}$ is given by $H(q,p)=\dot{q}\cdot
p-L(q,\dot{q}) $ with Hamiltonian vector field $X_H$,
which is defined by the Hamilton's equation
$\mathbf{i}_{X_H}\omega=\mathbf{d}H$, and
$\mathcal{M}=\mathcal{F}L(\mathcal{D})$ is a constraint submanifold
in $T^* Q$. In particular, for the nonholonomic constraint
$\mathcal{D}\subset TQ$, the Lagrangian $L$ is said to be {\bf
$\mathcal{D}$-regular}, if the restriction of Hessian matrix
$(\partial^2L/\partial\dot{q}^i\partial\dot{q}^j)$ on $\mathcal{D}$
is nondegenerate everywhere. Moreover, a nonholonomic system is said
to be {\bf $\mathcal{D}$-regular}, if its Lagrangian $L$ is {\bf
$\mathcal{D}$-regular}. Note that the restriction of a positive
definite symmetric bilinear form to a subspace is also positive
definite, and hence nondegenerate. Thus, for a simple nonholonomic
mechanical system; that is, whose Lagrangian is the total kinetic
energy minus potential energy, it is {\bf $\mathcal{D}$-regular}
automatically.

\subsection{Nonholonomic Hamiltonian System and Hamilton-Jacobi Equations}

A nonholonomic Hamiltonian system is a 4-tuple
$(T^\ast Q,\omega,\mathcal{D},H)$, which is a Hamiltonian system with a
$\mathcal{D}$-completely and $\mathcal{D}$-regularly nonholonomic
constraint $\mathcal{D} \subset TQ$. Under
the restriction given by constraints, in general, the dynamical
vector field of a nonholonomic Hamiltonian system may not be
Hamiltonian. However, the nonholonomic Hamiltonian system is
a dynamical system closely related to a Hamiltonian system.
In the following we can derive a distributional Hamiltonian system of the nonholonomic
Hamiltonian system $(T^*Q,\omega,\mathcal{D},H )$,
by analyzing carefully the structure for the nonholonomic
dynamical vector field and by using the method similar to
that used in Le\'{o}n and Wang \cite{lewa15} and
Bates and $\acute{S}$niatycki \cite{basn93}.
The distributional Hamiltonian system is
very important and it is also called a semi-Hamiltonian
system in Patrick \cite{pa07}.\\

Assume that $L:TQ \rightarrow \mathbb{R}$ is a
hyperregular Lagrangian, and that the Legendre transformation
$\mathcal{F}L: TQ \rightarrow T^*Q$ is a diffeomorphism.
We consider that the constraint submanifold
$\mathcal{M}=\mathcal{F}L(\mathcal{D})\subset T^*Q$ and
$i_{\mathcal{M}}: \mathcal{M}\rightarrow T^*Q $ is the inclusion,
and that the symplectic form $\omega_{\mathcal{M}}= i_{\mathcal{M}}^* \omega $,
is induced from the canonical symplectic form $\omega$ on $T^* Q$.
We define the distribution $\mathcal{F}$ as the pre-image of the nonholonomic
constraint $\mathcal{D}$ for the map $T\pi_Q: TT^* Q \rightarrow TQ$;
that is, $\mathcal{F}=(T\pi_Q)^{-1}(\mathcal{D})\subset TT^*Q,
$ which is a distribution along $\mathcal{M}$, and
$\mathcal{F}^\circ:=\{\alpha \in T^*T^*Q | <\alpha,v>=0, \; \forall
v\in TT^*Q \}$ is the annihilator of $\mathcal{F}$ in
$T^*T^*Q_{|\mathcal{M}}$. We consider the following nonholonomic
constraints condition that
\begin{align} (\mathbf{i}_X \omega -\mathbf{d}H) \in \mathcal{F}^\circ,
\;\;\;\;\;\; X \in T \mathcal{M}.
\label{6.1} \end{align}
From Cantrijn et al.
\cite{calemama99}, we know that there exists an unique nonholonomic
vector field $X_n$ satisfying the above condition $(6.1)$ if the
admissibility condition $\mathrm{dim}\mathcal{M}=
\mathrm{rank}\mathcal{F}$ and the compatibility condition
$T\mathcal{M}\cap \mathcal{F}^\bot= \{0\}$ hold, where
$\mathcal{F}^\bot$ denotes the symplectic orthogonal of
$\mathcal{F}$ with respect to the canonical symplectic form
$\omega$ on $T^*Q$. In particular, when we consider the Whitney sum
decomposition $T(T^*Q)_{|\mathcal{M}}=T\mathcal{M}\oplus
\mathcal{F}^\bot$ and the canonical projection $P:
T(T^*Q)_{|\mathcal{M}} \rightarrow T\mathcal{M}$,
then we have that $X_n= P(X_H)$.
See \cite{gone79} for the more details.\\

From the condition (6.1) we know that the nonholonomic vector field,
in general, may not be Hamiltonian, because of the restriction
of nonholonomic constraint. However, we hope to study the dynamical
vector field of nonholonomic Hamiltonian system by using the similar
method of studying Hamiltonian vector field.
From Le\'{o}n and Wang \cite{lewa15} and
Bates and $\acute{S}$niatycki \cite{basn93}, we can define the
distribution $ \mathcal{K}=\mathcal {F}\cap T\mathcal{M}.$ and
$\mathcal{K}^\bot=\mathcal {F}^\bot\cap T\mathcal{M}, $ where
$\mathcal{K}^\bot$ denotes the symplectic orthogonal of
$\mathcal{K}$ with respect to the canonical symplectic form
$\omega$. If the admissibility condition $\mathrm{dim}\mathcal{M}=
\mathrm{rank}\mathcal{F}$ and the compatibility condition
$T\mathcal{M}\cap \mathcal{F}^\bot= \{0\}$ hold, then we know that the
restriction of the symplectic form $\omega_{\mathcal{M}}$ on
$T^*\mathcal{M}$ fibrewise to the distribution $\mathcal{K}$; that
is, $\omega_\mathcal{K}= \tau_{\mathcal{K}}\cdot
\omega_{\mathcal{M}}$ is non-degenerate, where $\tau_{\mathcal{K}}$
is the restriction map to distribution $\mathcal{K}$. It is worthy
of noting that $\omega_\mathcal{K}$ is not a true two-form on a
manifold, so it does not make sense to speak about it being closed.
We call $\omega_\mathcal{K}$ as a distributional two-form to avoid
any confusion. Because $\omega_\mathcal{K}$ is non-degenerate as a
bilinear form on each fibre of $\mathcal{K}$, there exists a vector
field $X_{\mathcal{K}}$ on $\mathcal{M}$ which takes values in the
constraint distribution $\mathcal{K}$,
such that the distributional Hamiltonian equation holds; that is,
\begin{align}
\mathbf{i}_{X_\mathcal{K}}\omega_{\mathcal{K}}=\mathbf{d}H_\mathcal {K}
\label{6.2} \end{align}
where $\mathbf{d}H_\mathcal{K}$ is the restriction of
$\mathbf{d}H_\mathcal{M}$ to $\mathcal{K}$.
and the function $H_{\mathcal{K}}$ satisfies
$\mathbf{d}H_{\mathcal{K}}= \tau_{\mathcal{K}}\cdot \mathbf{d}H_{\mathcal {M}}$,
and $H_\mathcal{M}= \tau_{\mathcal{M}}\cdot H$ is the restriction of $H$ to
$\mathcal{M}$. Moreover, from the distributional Hamiltonian equation (6.2),
we have that $X_{\mathcal{K}}= \tau_{\mathcal{K}}\cdot X_H.$
Thus, the geometric formulation of a distributional Hamiltonian system may
be summarized as follows.

\begin{defi} (Distributional Hamiltonian System)
Assume that the 4-tuple $(T^*Q,\omega,\mathcal{D},H)$ is a nonholonomic
Hamiltonian system, where $\omega$ is the canonical
symplectic form on $T^* Q$, and $\mathcal{D}\subset TQ$ is a
$\mathcal{D}$-completely and $\mathcal{D}$-regularly nonholonomic
constraint of the system. If there exist a distribution
$\mathcal{K}$, an associated non-degenerate distributional two-form
$\omega_\mathcal{K}$ induced by the canonical symplectic form
and a vector field $X_\mathcal {K}$ on the
constraint submanifold $\mathcal{M}=\mathcal{F}L(\mathcal{D})\subset
T^*Q$, such that the distributional Hamiltonian equation holds; that is,
$\mathbf{i}_{X_\mathcal{K}}\omega_\mathcal{K}=\mathbf{d}H_\mathcal
{K}$, where $\mathbf{d}H_\mathcal{K}$ is the restriction of
$\mathbf{d}H_\mathcal{M}$ to $\mathcal{K}$, and
the function $H_{\mathcal{K}}$ satisfies
$\mathbf{d}H_{\mathcal{K}}= \tau_{\mathcal{K}}\cdot \mathbf{d}H_{\mathcal {M}}$
as defined above,
then the triple $(\mathcal{K},\omega_{\mathcal{K}},H_{\mathcal{K}})$
is called a distributional Hamiltonian system of the nonholonomic
Hamiltonian system $(T^*Q,\omega,\mathcal{D},H)$, and $X_\mathcal
{K}$ is called a nonholonomic dynamical
vector field of the distributional Hamiltonian system
$(\mathcal{K},\omega_{\mathcal {K}},H_{\mathcal{K}})$. Under the above circumstances, we refer to
$(T^*Q,\omega,\mathcal{D},H)$ as a nonholonomic Hamiltonian system
with an associated distributional Hamiltonian system
$(\mathcal{K},\omega_{\mathcal {K}},H_{\mathcal{K}})$.
\end{defi}

Since the non-degenerate
distributional two-form $\omega_{\mathcal{K}}$ is not symplectic,
and the distributional Hamiltonian system
$(\mathcal{K},\omega_{\mathcal {K}},H_{\mathcal{K}})$ is not yet a Hamiltonian system,
and hence, we can not describe the Hamilton-Jacobi equation for a
distributional Hamiltonian system the same as in Theorem 2.1.
However, for a given nonholonomic Hamiltonian system
$(T^*Q,\omega,\mathcal{D},H)$ with an associated distributional Hamiltonian
system $(\mathcal{K},\omega_{\mathcal {K}},H_{\mathcal{K}})$,
we can derive precisely the geometric constraint conditions of
the non-degenerate distributional two-form $\omega_\mathcal{K}$
for the nonholonomic dynamical vector field $X_\mathcal {K}$;
that is, the two types of Hamilton-Jacobi equations for the distributional Hamiltonian
system $(\mathcal{K},\omega_{\mathcal {K}},H_{\mathcal{K}})$. In order to do this, we need
first give an important notions and a key lemma, (see also Le\'{o}n and Wang \cite{lewa15}).
This lemma and the Lemma 2.3 offer an important tool for the proofs of the two types of Hamilton-Jacobi
theorems for the distributional Hamiltonian system and the nonholonomic
reduced distributional Hamiltonian system.\\

Assume that $\omega$ is the canonical symplectic form on $T^*Q$,
and $\mathcal{D}\subset TQ$ is a $\mathcal{D}$-regularly nonholonomic
constraint, and the projection $\pi_Q: T^* Q \rightarrow Q $ induces
the map $T\pi_{Q}: TT^* Q \rightarrow TQ. $
If the one-form $\gamma$ is closed,
then $\mathbf{d}\gamma(x,y)=0, \; \forall\; x, y \in TQ$.
Note that for any $v, w \in TT^* Q, $ we have that
$\mathbf{d}\gamma(T\pi_{Q}(v),T\pi_{Q}(w))=\pi^*(\mathbf{d}\gamma )(v, w)$
is a two-form on the cotangent bundle $T^*Q$, where
$\pi^*: T^*Q \rightarrow T^*T^*Q.$ Thus,
in the following we can introduce a weaker notion
which is an extension of the Definition 2.2 to the nonholonomic context.

\begin{defi}
The one-form $\gamma$ is called to be closed
on $\mathcal{D}$ with respect to $T\pi_{Q}:
TT^* Q \rightarrow TQ $ if, for any $v, w \in TT^* Q $
and $T\pi_{Q}(v), \; T\pi_{Q}(w) \in \mathcal{D},$  we have that
$\mathbf{d}\gamma(T\pi_{Q}(v),T\pi_{Q}(w))=0. $
\end{defi}

The notion that $\gamma$ is closed on $\mathcal{D}$ with respect to $T\pi_{Q}:
TT^* Q \rightarrow TQ, $ is weaker than that $\gamma$ is closed on $\mathcal{D}$;
that is, $\mathbf{d}\gamma(x,y)=0, \; \forall\; x, y \in \mathcal{D}$.
In fact, if $\gamma$ is a closed one-form on $\mathcal{D}$,
then it must be closed on $\mathcal{D}$ with respect to
$T\pi_{Q}: TT^* Q \rightarrow TQ. $
Conversely, if $\gamma$ is closed on $\mathcal{D}$ with respect to
$T\pi_{Q}: TT^* Q \rightarrow TQ, $ then it may not be closed on $\mathcal{D}$.
Now, we give the important lemma as follows:
\begin{lemm}
Assume that $\gamma: Q \rightarrow T^*Q$ is a one-form on $Q$, and
$\lambda=\gamma \cdot \pi_{Q}: T^* Q \rightarrow T^* Q .$
If the Lagrangian $L$ is $\mathcal{D}$-regular, and
$\textmd{Im}(\gamma)\subset \mathcal{M}=\mathcal{F}L(\mathcal{D}), $
then we have that $ X_{H}\cdot \gamma \in \mathcal{F}$ along
$\gamma$ and $ X_{H}\cdot \lambda \in \mathcal{F}$ along
$\lambda$; that is, $T\pi_{Q}(X_H\cdot\gamma(q))\in
\mathcal{D}_{q}, \; \forall q \in Q $, and $T\pi_{Q}(X_H\cdot\lambda(q,p))\in
\mathcal{D}_{q}, \; \forall q \in Q, \; (q,p) \in T^* Q. $
Moreover, if the symplectic map $\varepsilon: T^* Q \rightarrow T^* Q $
with respect to the canonical symplectic form $\omega$ on $T^*Q$,
satisfies the condition $\varepsilon(\mathcal{M})\subset \mathcal{M},$ then
we have that $ X_{H}\cdot \varepsilon \in \mathcal{F}$ along
$\varepsilon. $
\end{lemm}
See the proof and the more details in
Le\'{o}n and Wang \cite{lewa15}.\\

By using the Lemma2.3 and the above Lemma 6.3,
we can derive precisely the geometric constraint conditions of
the non-degenerate distributional two-form $\omega_{\mathcal{K}}$
for the nonholonomic dynamical vector field $X_\mathcal {K}$;
that is, the two types of
Hamilton-Jacobi equations for the distributional Hamiltonian
system $(\mathcal{K},\omega_{\mathcal {K}},H_{\mathcal{K}})$.
For convenience, the
maps involved in the theorem and its proof are shown in Diagram-5.

\begin{center}
\hskip 0cm \xymatrix{& \mathcal{M} \ar[d]_{X_{\mathcal{K}}}
\ar[r]^{i_{\mathcal{M}}} & T^* Q \ar[d]_{X_{H}}
\ar[dr]^{X^\varepsilon} \ar[r]^{\pi_Q}
& Q \ar[d]^{X^\gamma} \ar[r]^{\gamma} & T^*Q \ar[d]^{X_{H\cdot\varepsilon}} \\
& \mathcal{K}  & T(T^*Q) \ar[l]^{\tau_{\mathcal{K}}} & TQ
\ar[l]^{T\gamma} & T(T^* Q) \ar[l]^{T\pi_Q}}
\end{center}
$$\mbox{Diagram-5}$$

\begin{theo} ( Hamilton-Jacobi Theorem for a Distributional Hamiltonian System)
For the nonholonomic Hamiltonian system $(T^*Q,\omega,\mathcal{D},H)$
with an associated distributional Hamiltonian system
$(\mathcal{K},\omega_{\mathcal {K}},H_{\mathcal{K}})$, assume that $\gamma: Q
\rightarrow T^*Q$ is a one-form on $Q$, and that $\lambda=\gamma \cdot
\pi_{Q}: T^* Q \rightarrow T^* Q, $ and that for any
symplectic map $\varepsilon: T^* Q \rightarrow T^* Q $, denote that
$X^\gamma = T\pi_{Q}\cdot X_H \cdot \gamma$ and
$X^\varepsilon = T\pi_{Q}\cdot X_H \cdot \varepsilon$,
where $X_{H}$ is the dynamical
vector field of the corresponding unconstrained Hamiltonian system
$(T^*Q,\omega,H)$. Moreover, assume that $\textmd{Im}(\gamma)\subset
\mathcal{M}=\mathcal{F}L(\mathcal{D}), $
and that $\varepsilon(\mathcal{M})\subset \mathcal{M},$
and that $ \textmd{Im}(T\gamma)\subset \mathcal{K}. $
Then the following two assertions hold:\\
\noindent $(\mathbf{i})$ If the
one-form $\gamma: Q \rightarrow T^*Q $ is closed on $\mathcal{D}$ with respect to
$T\pi_Q: TT^* Q \rightarrow TQ, $ then $\gamma$ is a
solution of the Type I Hamilton-Jacobi equation $T\gamma \cdot
X^\gamma= X_{\mathcal{K}} \cdot \gamma $ for the distributional Hamiltonian system
$(\mathcal{K},\omega_{\mathcal {K}},H_{\mathcal{K}})$, where
$X_{\mathcal{K}}$
is the nonholonomic dynamical vector field of the distributional Hamiltonian system.\\
\noindent $(\mathbf{ii})$ The $\varepsilon$ is a solution of the
Type II Hamilton-Jacobi equation $T\gamma \cdot
X^\varepsilon= X_{\mathcal{K}} \cdot \varepsilon $
if and only if it is a solution of the equation
$\tau_{\mathcal{K}}\cdot T\varepsilon(X_{H\cdot\varepsilon})= T\lambda \cdot X_H\cdot\varepsilon,$
where $X_{H\cdot\varepsilon}$ is the Hamiltonian vector field of the function
$H\cdot\varepsilon: T^* Q\rightarrow \mathbb{R} $.
\end{theo}
See the proof and the more details in
Le\'{o}n and Wang \cite{lewa15}.

\subsection{Nonholonomic Controlled Hamiltonian System and Its Reductions}

A nonholonomic RCH system is a 6-tuple $(T^\ast Q,\omega,\mathcal{D},H,F,W)$,
which is an RCH system with a
$\mathcal{D}$-completely and $\mathcal{D}$-regularly nonholonomic
constraint $\mathcal{D} \subset TQ$. Under
the restriction given by the constraints, in general, the dynamical
vector field of a nonholonomic RCH system may not be
Hamiltonian vector field. However the nonholonomic RCH system is
a dynamical system closely related to a Hamiltonian system.
In the following we can derive a distributional RCH system of the nonholonomic
RCH system $(T^*Q,\omega,\mathcal{D},H,F,W)$,
by analyzing carefully the structure for the nonholonomic
dynamical vector field and by using the method similar
to that used in Le\'{o}n and Wang \cite{lewa15}.\\

From the above discussion in \S 6.1, we know that
there exist a distribution
$\mathcal{K}$, an associated non-degenerate distributional two-form
$\omega_\mathcal{K}$ induced by the canonical symplectic form
and a vector field $X_\mathcal {K}$ on the
constraint submanifold $\mathcal{M}=\mathcal{F}L(\mathcal{D})\subset
T^*Q$, such that the distributional Hamiltonian equation (6.2) holds,
and we have that $X_{\mathcal{K}}= \tau_{\mathcal{K}}\cdot X_H.$\\

Moreover, if considering the external force $F$ and the control subset $W$,
and define that $F_\mathcal{K}=\tau_{\mathcal{K}}\cdot \textnormal{vlift}(F_{\mathcal{M}})X_H,$
and that, for the control law $u\in W$,
$u_\mathcal{K}= \tau_{\mathcal{K}}\cdot  \textnormal{vlift}(u_{\mathcal{M}})X_H,$
where $F_\mathcal{M}= \tau_{\mathcal{M}}\cdot F$ and
$u_\mathcal{M}= \tau_{\mathcal{M}}\cdot u$ are the restrictions of
$F$ and $u$ to $\mathcal{M}$; that is, $F_\mathcal{K}$ and $u_\mathcal{K}$
are the restrictions of the changes of the Hamiltonian vector field $X_H$
under the actions of $F_\mathcal{M}$ and $u_\mathcal{M}$ to $\mathcal{K}$.
Then the 5-tuple $(\mathcal{K},\omega_{\mathcal{K}},
H_\mathcal{K}, F_\mathcal{K}, u_\mathcal{K})$
is a distributional RCH system of the nonholonomic
RCH system $(T^*Q,\omega,\mathcal{D},H,F,W)$ with the control law $u\in W$.
Thus, the geometric formulation of the distributional RCH
system may be summarized as follows:

\begin{defi} (Distributional RCH System)
Assume that the 6-tuple $(T^*Q,\omega,\mathcal{D},H,F,W)$ is a nonholonomic
RCH system, where $\omega$ is the canonical
symplectic form on $T^* Q$, and $\mathcal{D}\subset TQ$ is a
$\mathcal{D}$-completely and $\mathcal{D}$-regularly nonholonomic
constraint of the system, and the external force
$F: T^*Q\rightarrow T^*Q$ is the fiber-preserving map,
and the control subset $W\subset T^*Q$ is a fiber submanifold of $T^*Q$.
For a control law $u\in W,$ if there exist a distribution
$\mathcal{K}$, an associated non-degenerate distributional two-form
$\omega_\mathcal{K}$ induced by the canonical symplectic form
and a vector field $X_\mathcal {K}$ on the
constraint submanifold $\mathcal{M}=\mathcal{F}L(\mathcal{D})\subset
T^*Q$, such that the distributional Hamiltonian equation holds; that is,
$\mathbf{i}_{X_\mathcal{K}}\omega_\mathcal{K}=\mathbf{d}H_\mathcal
{K}$, where $\mathbf{d}H_\mathcal{K}$ is the restriction of
$\mathbf{d}H_\mathcal{M}$ to $\mathcal{K}$, and
the function $H_{\mathcal{K}}$ satisfies
$\mathbf{d}H_{\mathcal{K}}= \tau_{\mathcal{K}}\cdot \mathbf{d}H_{\mathcal {M}},$
and $F_{\mathcal{K}}=\tau_{\mathcal{K}}\cdot \textnormal{vlift}(F_{\mathcal{M}})X_H$,
and $u_{\mathcal{K}}=\tau_{\mathcal{K}}\cdot \textnormal{vlift}(u_{\mathcal{M}})X_H$
as defined above,
then the 5-tuple $(\mathcal{K},\omega_{\mathcal{K}},
H_{\mathcal{K}}, F_{\mathcal{K}}, u_{\mathcal{K}})$
is called a distributional RCH system of the nonholonomic
RCH system $(T^*Q,\omega,\mathcal{D},H,F,u)$, and $X_\mathcal
{K}$ is called a nonholonomic dynamical vector field. Denote that
\begin{align}
\tilde{X}=X_{(\mathcal{K},\omega_{\mathcal{K}},
H_{\mathcal{K}}, F_{\mathcal{K}}, u_{\mathcal{K}})}
=X_\mathcal {K}+ F_{\mathcal{K}}+u_{\mathcal{K}}
\label{6.3} \end{align}
is the dynamical vector field of the distributional RCH system
$(\mathcal{K},\omega_{\mathcal {K}}, H_{\mathcal{K}}, F_{\mathcal{K}}, u_{\mathcal{K}})$,
which is the synthesis
of the nonholonomic dynamical vector field $X_{\mathcal{K}}$ and
the vector fields $F_{\mathcal{K}}$ and $u_{\mathcal{K}}$.
Under the above circumstances, we refer to
$(T^*Q,\omega,\mathcal{D}, H, F, u)$ as a nonholonomic RCH system
with an associated distributional RCH system
$(\mathcal{K},\omega_{\mathcal {K}}, H_{\mathcal{K}}, F_{\mathcal{K}}, u_{\mathcal{K}})$.
\end{defi}

It is worthy of noting that
if the external force and control of a distributional RCH
system $(\mathcal{K},\omega_{\mathcal {K}},H_{\mathcal {K}},
F_{\mathcal {K}}, u_{\mathcal {K}})$ are both zeros;
that is, $F_{\mathcal {K}}=0 $ and $u_{\mathcal {K}}=0$, in this case,
the distributional RCH system is just a distributional Hamiltonian system
$(\mathcal{K},\omega_{\mathcal {K}},H_{\mathcal {K}})$;
see Le\'{o}n and Wang \cite{lewa15} for more details.
Thus, the distributional RCH system
can be regarded as an extension of the distributional Hamiltonian system to
the system with the external force and the control. \\

We know that the reduction of the nonholonomically constrained mechanical systems
is also very important subject in the study of the geometric mechanics, which is
regarded as a useful tool for simplifying and studying
concrete nonholonomic systems; see
Le\'{o}n and Wang \cite{lewa15},
Bates and $\acute{S}$niatycki \cite{basn93},  Cantrijn \textit{et al.}
\cite{calemama99}, Cendra \textit{et al.} \cite{cemara01},
Cushman \textit{et al.} \cite{cudusn10} and \cite{cukesnba95}, Koiller \cite{ko92},
Le\'{o}n and Rodrigues \cite{lero89}, and so on.
In the following we consider the nonholonomic RCH system
with symmetry, as well as with a momentum map.
By using the similar method for the nonholonomic reduction
(given in Bates and $\acute{S}$niatycki \cite{basn93},
Le\'{o}n and Wang \cite{lewa15} and Wang \cite{wa21c}),
and analyzing carefully the structures for the nonholonomic reduced
dynamical vector fields, we also give the geometric formulations of
the nonholonomic reduced distributional RCH systems.\\

In the following we consider that a nonholonomic RCH system
with symmetry and a momentum map
is 8-tuple $(T^*Q,G, \omega,\mathbf{J},\mathcal{D},H,F,W)$,
where $\omega$ is the canonical symplectic form on $T^* Q$,
and the Lie group $G$, which may not be Abelian, acts smoothly by the left on $Q$,
its tangent lifted action on $TQ$ and its cotangent lifted action on $T^\ast Q$,
and $\mathcal{D}\subset TQ$ is a
$\mathcal{D}$-completely and $\mathcal{D}$-regularly nonholonomic
constraint of the system, and $\mathcal{D}$, $H, F$ and $W$ are all
$G$-invariant. Thus, the nonholonomic RCH system with symmetry and a momentum map
is a regular point reducible RCH system with $G$-invariant
nonholonomic constraint $\mathcal{D}$.
Moreover, in the following we shall give
carefully a geometric formulation of the $\mathbf{J}$-nonholonomic
$R_p$-reduced distributional RCH system, by using the momentum map and the
nonholonomic reduction compatible with regular point reduction.\\

Note that the Legendre transformation $\mathcal{F}L: TQ \rightarrow T^*Q $
is a fiber-preserving map, and that $\mathcal{D}\subset TQ$ is $G$-invariant
for the tangent lifted left action $\Phi^{T}: G\times TQ\rightarrow TQ, $
then the constraint submanifold
$\mathcal{M}=\mathcal{F}L(\mathcal{D})\subset T^*Q$ is
$G$-invariant for the cotangent lifted left action $\Phi^{T^\ast}:
G\times T^\ast Q\rightarrow T^\ast Q$.
For the nonholonomic RCH system with symmetry
and a momentum map  $(T^*Q,G, \omega,\mathbf{J},\mathcal{D},H,F,W)$,
in the same way, we define the distribution $\mathcal{F}$, which is the pre-image of the
nonholonomic constraints $\mathcal{D}$ for the map $T\pi_Q: TT^* Q
\rightarrow TQ$; that is, $\mathcal{F}=(T\pi_Q)^{-1}(\mathcal{D})$,
and the distribution $\mathcal{K}=\mathcal{F} \cap T\mathcal{M}$.
Moreover, we can also define the distributional two-form $\omega_\mathcal{K}$,
which is induced from the canonical symplectic form $\omega$ on $T^* Q$; that is,
$\omega_\mathcal{K}= \tau_{\mathcal{K}}\cdot \omega_{\mathcal{M}},$ and
$\omega_{\mathcal{M}}= i_{\mathcal{M}}^* \omega $.
If the admissibility condition $\mathrm{dim}\mathcal{M}=
\mathrm{rank}\mathcal{F}$ and the compatibility condition
$T\mathcal{M}\cap \mathcal{F}^\bot= \{0\}$ hold, then
$\omega_\mathcal{K}$ is non-degenerate as a
bilinear form on each fibre of $\mathcal{K}$, there exists a vector
field $X_\mathcal{K}$ on $\mathcal{M}$ which takes values in the
constraint distribution $\mathcal{K}$, such that for the function $H_\mathcal{K}$,
the following distributional Hamiltonian equation holds, that is,
\begin{align}
\mathbf{i}_{X_\mathcal{K}}\omega_\mathcal{K}
=\mathbf{d}H_\mathcal{K},
\label{6.4} \end{align}
where the function $H_{\mathcal{K}}$ satisfies
$\mathbf{d}H_{\mathcal{K}}= \tau_{\mathcal{K}}\cdot \mathbf{d}H_{\mathcal {M}}$,
and $H_\mathcal{M}= \tau_{\mathcal{M}}\cdot H$
is the restriction of $H$ to $\mathcal{M}$, and
from the equation (6.4), we have that
$X_{\mathcal{K}}=\tau_{\mathcal{K}}\cdot X_H $.\\

Since the nonholonomic RCH system with symmetry and a momentum map
is a regular point reducible RCH system with $G$-invariant
nonholonomic constraint $\mathcal{D}$,
for a regular value $\mu\in\mathfrak{g}^\ast$ of the
momentum map $\mathbf{J}:T^\ast Q\rightarrow \mathfrak{g}^\ast$,
we assume that the constraint submanifold $\mathcal{M}$
is clean intersection with $\mathbf{J}^{-1}(\mu)$; that is,
$\mathcal{M} \cap \mathbf{J}^{-1}(\mu)\neq \emptyset$.
Note that $\mathcal{M}$ is also $G_\mu (\subset G)$ action
invariant, and so is $\mathbf{J}^{-1}(\mu)$, because $\mathbf{J}$ is
$\operatorname{Ad}^\ast$-equivariant. It follows that the quotient
space $\mathcal{M}_\mu =(\mathcal{M}\cap \mathbf{J}^{-1}(\mu))
/G_\mu \subset (T^\ast Q)_\mu$ of the $G_\mu$-orbit in
$\mathcal{M}\cap \mathbf{J}^{-1}(\mu)$, is a smooth manifold with
the projection $\pi_\mu: \mathcal{M}\cap \mathbf{J}^{-1}(\mu)
\rightarrow \mathcal{M}_\mu$, which is a surjective submersion.
Denote that $i_{\mathcal{M}_\mu}: \mathcal{M}_\mu\rightarrow
(T^*Q)_\mu, $ and that $\omega_{\mathcal{M}_\mu}= i_{\mathcal{M}_\mu}^*
\omega_\mu $; that is, the symplectic form
$\omega_{\mathcal{M}_\mu}$ is induced from the $R_p$-reduced symplectic
form $\omega_\mu$ on $(T^* Q)_\mu$ given in (4.1), where $i_{\mathcal{M}_\mu}^*:
T^*(T^*Q)_\mu \rightarrow T^*\mathcal{M}_\mu. $ Moreover, the
distribution $\mathcal{F}$ is pushed down to a distribution
$\mathcal{F}_\mu= T\pi_\mu\cdot \mathcal{F}$ on $(T^\ast Q)_\mu$,
and we define $\mathcal{K}_\mu=\mathcal{F}_\mu \cap
T\mathcal{M}_\mu$. Assume that $\omega_{\mathcal{K}_\mu}=
\tau_{\mathcal{K}_\mu}\cdot \omega_{\mathcal{M}_\mu}$ is the
restriction of the symplectic form $\omega_{\mathcal{M}_\mu}$ on
$T^*\mathcal{M}_\mu$ fibrewise to the distribution $\mathcal{K}_\mu$,
where $\tau_{\mathcal{K}_\mu}$ is the restriction map to distribution
$\mathcal{K}_\mu$. The distributional two-form $\omega_{\mathcal{K}_\mu}$ is not
a ''true two-form'' on a manifold, which is called
as a $\mathbf{J}$-nonholonomic $R_p$-reduced
distributional two-form to avoid any confusion.\\

From the above construction we know that,
if the admissibility condition $\mathrm{dim}\mathcal{M}_\mu=
\mathrm{rank}\mathcal{F}_\mu$ and the compatibility condition
$T\mathcal{M}_\mu\cap \mathcal{F}_\mu^\bot= \{0\}$ hold, where
$\mathcal{F}_\mu^\bot$ denotes the symplectic orthogonal of
$\mathcal{F}_\mu$ with respect to the $R_p$-reduced symplectic form
$\omega_\mu$, then $\omega_{\mathcal{K}_\mu}$ is non-degenerate
as a bilinear form on each fibre of $\mathcal{K}_\mu$,
and hence, there exists a vector field $X_{\mathcal{K}_\mu}$
on $\mathcal{M}_\mu$, which takes values in the constraint
distribution $\mathcal{K}_\mu$, such that for the function
$h_{\mathcal {K}_\mu}$, the $\mathbf{J}$-nonholonomic
$R_p$-reduced distributional Hamiltonian equation holds, that is,
\begin{align}
\mathbf{i}_{X_{\mathcal{K}_\mu}}\omega_{\mathcal{K}_\mu}
=\mathbf{d}h_{\mathcal{K}_\mu},
\label{6.5} \end{align}
where $\mathbf{d}h_{\mathcal{K}_\mu}$ is the restriction
of $\mathbf{d}h_{\mathcal{M}_\mu}$ to $\mathcal{K}_\mu$,
and the function $h_{\mathcal {K}_\mu}$ satisfies
$\mathbf{d}h_{\mathcal{K}_\mu}= \tau_{\mathcal{K}_\mu}\cdot \mathbf{d}h_{\mathcal{M}_\mu} $,
and $h_{\mathcal{M}_\mu}= \tau_{\mathcal{M}_\mu}\cdot h_\mu$ is the
restriction of $h_\mu$ to $\mathcal{M}_\mu$,
and $h_\mu$ is the $R_p$-reduced
Hamiltonian function $h_\mu: (T^* Q)_\mu \rightarrow \mathbb{R}$ defined
by $h_\mu\cdot \pi_\mu= H\cdot i_\mu$.
In addition, from the distributional Hamiltonian equation (6.4),
$\mathbf{i}_{X_\mathcal{K}}\omega_\mathcal{K}=\mathbf{d}H_\mathcal{K},$
we have that $X_{\mathcal{K}}=\tau_{\mathcal{K}}\cdot X_H, $
and from the $\mathbf{J}$-nonholonomic $R_p$-reduced
distributional Hamiltonian equation (6.5),
$\mathbf{i}_{X_{\mathcal{K_\mu}}}\omega_{\mathcal{K_\mu}}
=\mathbf{d}h_{\mathcal{K_\mu}}$, we have that
$X_{\mathcal{K_\mu}}
=\tau_{\mathcal{K_\mu}}\cdot X_{h_{\mathcal{K_\mu}}},$
where $ X_{h_{\mathcal{K_\mu}}}$ is the Hamiltonian vector field of
the function $h_{\mathcal{K_\mu}},$
and the vector fields $X_{\mathcal{K}}$
and $X_{\mathcal{K_\mu}}$ are $\pi_{\mu}$-related,
that is, $X_{\mathcal{K_\mu}}\cdot \pi_{\mu}=T\pi_{\mu}\cdot X_{\mathcal{K}}.$ \\

Moreover, if considering the external force $F$ and the control subset $W$,
and we define the vector fields $F_\mathcal{K}
=\tau_{\mathcal{K}}\cdot \textnormal{vlift}(F_{\mathcal{M}})X_H,$
and for a control law $u\in W$, define that
$u_\mathcal{K}= \tau_{\mathcal{K}}\cdot  \textnormal{vlift}(u_{\mathcal{M}})X_H,$
where $F_\mathcal{M}= \tau_{\mathcal{M}}\cdot F$ and
$u_\mathcal{M}= \tau_{\mathcal{M}}\cdot u$ are the restrictions of
$F$ and $u$ to $\mathcal{M}$; that is, $F_\mathcal{K}$ and $u_\mathcal{K}$
are the restrictions of the changes of Hamiltonian vector field $X_H$
under the actions of $F_\mathcal{M}$ and $u_\mathcal{M}$ to $\mathcal{K}$,
then the 5-tuple $(\mathcal{K},\omega_{\mathcal{K}},
H_\mathcal{K}, F_\mathcal{K}, u_\mathcal{K})$
is a distributional RCH system corresponding to the nonholonomic RCH system with symmetry
and a momentum map $(T^*Q,G, \omega,\mathbf{J},\mathcal{D},H,F,u)$,
and the dynamical vector field of the distributional RCH system
can be expressed by
\begin{align}
\tilde{X}=X_{(\mathcal{K},\omega_{\mathcal{K}},
H_{\mathcal{K}}, F_{\mathcal{K}}, u_{\mathcal{K}})}
=X_\mathcal {K}+ F_{\mathcal{K}}+u_{\mathcal{K}},
\label{6.6} \end{align}
which is the synthesis
of the nonholonomic dynamical vector field $X_{\mathcal{K}}$ and
the vector fields $F_{\mathcal{K}}$ and $u_{\mathcal{K}}$.
Assume that the vector fields $F_\mathcal{K}$ and $u_\mathcal{K}$
on $\mathcal{M}$ are pushed down to the vector fields
$f_{\mathcal{M}_\mu}=T\pi_\mu \cdot F_\mathcal{K}$ and
$u_{\mathcal{M}_\mu}=T\pi_\mu \cdot u_\mathcal{K}$ on $\mathcal{M}_\mu$.
Then we define that $f_{\mathcal{K}_\mu}
=T\pi_{\mathcal{K}_\mu}\cdot f_{\mathcal{M}_\mu}$, and that
$u_{\mathcal{K}_\mu}=T\pi_{\mathcal{K}_\mu}\cdot u_{\mathcal{M}_\mu};$
that is, $f_{\mathcal{K}_\mu}$ and
$u_{\mathcal{K}_\mu}$ are the restrictions of
$f_{\mathcal{M}_\mu}$ and $u_{\mathcal{M}_\mu}$ to $\mathcal{K}_\mu$.
As a consequence, the 5-tuple $(\mathcal{K}_\mu,\omega_{\mathcal{K}_\mu},
h_{\mathcal{K}_\mu}, f_{\mathcal{K}_\mu}, u_{\mathcal{K}_\mu})$
is a $\mathbf{J}$-nonholonomic $R_p$-reduced distributional RCH system of the nonholonomic
RCH system with symmetry and a momentum map
$(T^*Q,G,\omega,\mathbf{J},\mathcal{D},H,F,W)$,
as well as with a control law $u \in W$.
Thus, the geometrical formulation
of the $\mathbf{J}$-nonholonomic $R_p$-reduced distributional RCH
system may be summarized as follows.

\begin{defi} ($\mathbf{J}$-Nonholonomic $R_p$-reduced Distributional RCH System)
Assume that the 8-tuple $(T^*Q,G,\omega,\mathbf{J},\mathcal{D},H,F,W)$
is a nonholonomic RCH system with symmetry and a momentum map,
where $\omega$ is the canonical symplectic form on $T^* Q$,
and $\mathcal{D}\subset TQ$ is a $\mathcal{D}$-completely and
$\mathcal{D}$-regularly nonholonomic constraint of the system, and
$\mathcal{D}$, $H, F$ and $W$ are all $G$-invariant. For a regular value
$\mu \in \mathfrak{g}^\ast$ of the momentum map $\mathbf{J}:T^\ast
Q\rightarrow \mathfrak{g}^\ast$, assume that there exists a
$\mathbf{J}$-nonholonomic $R_p$-reduced distribution
$\mathcal{K}_\mu$, an associated non-degenerate and $\mathbf{J}$-nonholonomic
$R_p$-reduced distributional two-form
$\omega_{\mathcal{K}_\mu}$ and a vector field $X_{\mathcal {K}_\mu}$
on the $\mathbf{J}$-nonholonomic $R_p$-reduced constraint submanifold
$\mathcal{M}_\mu=(\mathcal{M}\cap \mathbf{J}^{-1}(\mu)) /G_\mu, $
where $\mathcal{M}=\mathcal{F}L(\mathcal{D}),$ and $\mathcal{M}\cap
\mathbf{J}^{-1}(\mu)\neq \emptyset, $ and $G_\mu=\{g\in G \;| \;
\operatorname{Ad}_g^\ast \mu=\mu \}$, such that the
$\mathbf{J}$-nonholonomic $R_p$-reduced distributional Hamiltonian
equation (6.5) holds, that is,
$ \mathbf{i}_{X_{\mathcal{K}_\mu}}\omega_{\mathcal{K}_\mu} =
\mathbf{d}h_{\mathcal{K}_\mu}, $
where $\mathbf{d}h_{\mathcal{K}_\mu}$ is the restriction of
$\mathbf{d}h_{\mathcal{M}_\mu}$ to $\mathcal{K}_\mu$,
and the function $h_{\mathcal {K}_\mu}$, and the vector fields $f_{\mathcal {K}_\mu}$
and $u_{\mathcal {K}_\mu}$ are defined above. Then the 5-tuple
$(\mathcal{K}_\mu,\omega_{\mathcal {K}_\mu},h_{\mathcal {K}_\mu},
f_{\mathcal {K}_\mu}, u_{\mathcal {K}_\mu})$ is called a
$\mathbf{J}$-nonholonomic $R_p$-reduced distributional RCH system
of the nonholonomic RCH system with symmetry and a momentum map
$(T^*Q,G,\omega,\mathbf{J},\mathcal{D},H, F,u)$
with a control law $u \in W$, and $X_{\mathcal{K}_\mu}$
is called the $\mathbf{J}$-nonholonomic $R_p$-reduced
dynamical vector field. Denote that
\begin{align}
\hat{X}_\mu=X_{(\mathcal{K}_\mu,\omega_{\mathcal{K}_\mu},h_{\mathcal {K}_\mu},
f_{\mathcal {K}_\mu}, u_{\mathcal {K}_\mu})}
=X_{\mathcal{K}_\mu}+ f_{\mathcal{K}_\mu}+u_{\mathcal{K}_\mu}
\label{6.7} \end{align}
is the dynamical vector field of the
$\mathbf{J}$-nonholonomic $R_p$-reduced distributional RCH system
$(\mathcal{K}_\mu,\omega_{\mathcal{K}_\mu}, \\ h_{\mathcal {K}_\mu},
f_{\mathcal {K}_\mu}, u_{\mathcal {K}_\mu})$,
which is the synthesis
of the $\mathbf{J}$-nonholonomic $R_p$-reduced
dynamical vector field $X_{\mathcal{K}_\mu}$ and
the vector fields $f_{\mathcal{K}_\mu}$ and $u_{\mathcal{K}_\mu}$.
Under the above circumstances, we refer to
$(T^*Q,G,\omega,\mathbf{J},\mathcal{D},H, F,u)$ as
a $\mathbf{J}$-nonholonomic point reducible RCH system
with the associated distributional RCH system
$(\mathcal{K},\omega_{\mathcal {K}},H_{\mathcal{K}}, F_{\mathcal{K}}, u_{\mathcal{K}})$
and the $\mathbf{J}$-nonholonomic $R_p$-reduced
distributional RCH system
$(\mathcal{K}_\mu,\omega_{\mathcal{K}_\mu},h_{\mathcal {K}_\mu},\\
f_{\mathcal {K}_\mu}, u_{\mathcal {K}_\mu})$.
\end{defi}

Since the non-degenerate and $\mathbf{J}$-nonholonomic $R_p$-reduced
distributional two-form $\omega_{\mathcal{K}_\mu}$ is not a "true two-form"
on a manifold, and it is not symplectic, and hence
the $\mathbf{J}$-nonholonomic $R_p$-reduced distributional RCH system
$(\mathcal{K}_\mu,\omega_{\mathcal {K}_\mu},h_{\mathcal {K}_\mu},
f_{\mathcal {K}_\mu}, u_{\mathcal {K}_\mu})$
is not a Hamiltonian system, and has no yet generating function,
and hence, we can not describe the Hamilton-Jacobi equation for a
$\mathbf{J}$-nonholonomic $R_p$-reduced
distributional RCH system the same as in Theorem 2.1.
However, for a given $\mathbf{J}$-nonholonomic regular point reducible RCH system $(T^*Q,G,\omega,\mathbf{J},\mathcal{D},H, F, u)$ with the associated distributional RCH system
$(\mathcal{K},\omega_{\mathcal {K}},H_{\mathcal{K}}, F_{\mathcal{K}}, u_{\mathcal{K}})$
and the $\mathbf{J}$-nonholonomic $R_p$-reduced distributional RCH system
$(\mathcal{K}_\mu,\omega_{\mathcal {K}_\mu},h_{\mathcal {K}_\mu},
f_{\mathcal {K}_\mu}, u_{\mathcal {K}_\mu})$, by using Lemma 2.3 and 6.3,
we can derive precisely
the geometric constraint conditions of the $\mathbf{J}$-nonholonomic $R_p$-reduced distributional two-form
$\omega_{\mathcal{K}_\mu}$ for the $\mathbf{J}$-nonholonomic regular point reducible dynamical vector field;
that is, the two types of Hamilton-Jacobi equations for the
$\mathbf{J}$-nonholonomic $R_p$-reduced distributional RCH system
$(\mathcal{K}_\mu,\omega_{\mathcal {K}_\mu},h_{\mathcal {K}_\mu},
f_{\mathcal {K}_\mu}, u_{\mathcal {K}_\mu})$.
At first, by using the fact that the one-form $\gamma: Q
\rightarrow T^*Q $ is closed on $\mathcal{D}$ with respect to
$T\pi_Q: TT^* Q \rightarrow TQ, $
and that $\textmd{Im}(\gamma)\subset \mathcal{M} \cap
\mathbf{J}^{-1}(\mu), $ and that $\gamma$ is $G_\mu$-invariant,
as well as that $ \textmd{Im}(T\bar{\gamma}_\mu)\subset \mathcal{K}_\mu, $
we can prove the Type I
Hamilton-Jacobi theorem for the $\mathbf{J}$-nonholonomic $R_p$-reduced distributional
RCH system. For convenience, the maps involved in the
following theorem and its proof are shown in Diagram-6.
\begin{center}
\hskip 0cm \xymatrix{ \mathbf{J}^{-1}(\mu) \ar[r]^{i_\mu}
& T^* Q  \ar[r]^{\pi_Q}
& Q \ar[d]_{\tilde{X}^\gamma} \ar[r]^{\gamma}
& T^*Q \ar[d]_{\tilde{X}} \ar[r]^{\pi_\mu} & (T^* Q)_\mu \ar[d]_{\hat{X}_\mu}
& \mathcal{M}_\mu  \ar[l]_{i_{\mathcal{M}_\mu}} \ar[d]_{X_{\mathcal{K}_\mu}}\\
& T(T^*Q)  & TQ \ar[l]_{T\gamma} & T(T^*Q) \ar[l]^{T\pi_Q} \ar[r]_{T\pi_\mu}
& T(T^* Q)_\mu \ar[r]^{\tau_{\mathcal{K}_\mu}} & \mathcal{K}_\mu }
\end{center}
$$\mbox{Diagram-6}$$

\begin{theo} (Type I Hamilton-Jacobi Theorem for a $\mathbf{J}$-Nonholonomic
$R_p$-reduced Distributional RCH System) For a given
$\mathbf{J}$-nonholonomic regular point reducible RCH system
$(T^*Q,G,\omega,\mathbf{J}, \mathcal{D}, H, \\ F, u)$
with the associated distributional RCH system
$(\mathcal{K},\omega_{\mathcal {K}},H_{\mathcal{K}}, F_{\mathcal{K}}, u_{\mathcal{K}})$
and the $\mathbf{J}$-nonholonomic $R_p$-reduced distributional RCH system
$(\mathcal{K}_\mu,\omega_{\mathcal{K}_\mu}, h_{\mathcal {K}_\mu},
f_{\mathcal {K}_\mu}, u_{\mathcal {K}_\mu})$, assume that $\gamma:
Q \rightarrow T^*Q$ is a one-form on $Q$, and that
$\tilde{X}^\gamma = T\pi_{Q}\cdot \tilde{X}\cdot \gamma$,
where $\tilde{X}=X_{(\mathcal{K},\omega_{\mathcal{K}},
H_{\mathcal{K}}, F_{\mathcal{K}}, u_{\mathcal{K}})}
=X_\mathcal {K}+ F_{\mathcal{K}}+u_{\mathcal{K}}$
is the dynamical vector field of the distributional RCH system
$(\mathcal{K},\omega_{\mathcal{K}},
H_{\mathcal{K}}, F_{\mathcal{K}}, u_{\mathcal{K}})$
corresponding to the $\mathbf{J}$-nonholonomic regular point reducible RCH
system with symmetry and a momentum map $(T^*Q,G,\omega,\mathbf{J},\mathcal{D},H, F,u)$.
Moreover, assume that $\mu\in\mathfrak{g}^\ast$ is a regular value of the momentum
map $\mathbf{J}$, and that $\textmd{Im}(\gamma)\subset \mathcal{M} \cap
\mathbf{J}^{-1}(\mu), $ and that $\gamma$ is $G_\mu$-invariant, and denote that
$\bar{\gamma}_\mu=\pi_\mu(\gamma): Q \rightarrow \mathcal{M}_\mu $,
and that $ \textmd{Im}(T\bar{\gamma}_\mu)\subset \mathcal{K}_\mu. $
If the one-form $\gamma: Q \rightarrow T^*Q $ is closed on $\mathcal{D}$ with respect to
$T\pi_Q: TT^* Q \rightarrow TQ, $ then
$\bar{\gamma}_\mu$ is a solution of the
equation $T\bar{\gamma}_\mu\cdot
\tilde{X}^\gamma= X_{\mathcal{K}_\mu}\cdot \bar{\gamma}_\mu. $
Here $X_{\mathcal{K}_\mu}$ is the $\mathbf{J}$-nonholonomic $R_p$-reduced
dynamical vector field. The equation
$T\bar{\gamma}_\mu\cdot \tilde{X}^\gamma= X_{\mathcal{K}_\mu}\cdot
\bar{\gamma}_\mu,$ is called the Type I Hamilton-Jacobi equation for the
$\mathbf{J}$-nonholonomic $R_p$-reduced distributional RCH system
$(\mathcal{K}_\mu,\omega_{\mathcal{K}_\mu},h_{\mathcal {K}_\mu},
f_{\mathcal {K}_\mu}, u_{\mathcal {K}_\mu})$.
\end{theo}
See the proof and the more details in Wang \cite{wa22a}.
Moreover, for any $G_\mu$-invariant symplectic map $\varepsilon: T^* Q \rightarrow T^* Q $,
we can prove the Type II
Hamilton-Jacobi theorem for the $\mathbf{J}$-nonholonomic
$R_p$-reduced distributional RCH system.\\

In the following we consider that a nonholonomic RCH system
with symmetry and momentum map
is 8-tuple $(T^*Q,G, \omega,\mathbf{J},\mathcal{D},H,F,W)$,
which is a regular orbit reducible RCH system with $G$-invariant
nonholonomic constraint $\mathcal{D}$.
We can give the geometric formulation of the $\mathbf{J}$-nonholonomic
$R_o$-reduced distributional RCH system, by using the momentum map and the
nonholonomic reduction compatible with the regular orbit reduction.\\

For a regular value $\mu\in\mathfrak{g}^\ast$ of the
momentum map $\mathbf{J}:T^\ast Q\rightarrow \mathfrak{g}^\ast$,
$\mathcal{O}_\mu=G\cdot \mu\subset \mathfrak{g}^\ast$ is the
$G$-orbit of the coadjoint $G$-action through the point $\mu$,
we assume that the constraint submanifold $\mathcal{M}$
is clean intersection with $\mathbf{J}^{-1}(\mathcal{O}_\mu)$; that is,
$\mathcal{M} \cap \mathbf{J}^{-1}(\mathcal{O}_\mu)\neq \emptyset$.
It follows that the quotient
space $\mathcal{M}_{\mathcal{O}_\mu} =(\mathcal{M}\cap \mathbf{J}^{-1}(\mathcal{O}_\mu))
/G \subset (T^\ast Q)_{\mathcal{O}_\mu}$ of the $G$-orbit in
$\mathcal{M}\cap \mathbf{J}^{-1}(\mathcal{O}_\mu)$, is a smooth manifold with
projection $\pi_{\mathcal{O}_\mu}: \mathcal{M}\cap \mathbf{J}^{-1}(\mathcal{O}_\mu)
\rightarrow \mathcal{M}_{\mathcal{O}_\mu}$ which is a surjective submersion.
Denote that $i_{\mathcal{M}_{\mathcal{O}_\mu}}: \mathcal{M}_{\mathcal{O}_\mu}\rightarrow
(T^*Q)_{\mathcal{O}_\mu}, $ and that $\omega_{\mathcal{M}_{\mathcal{O}_\mu}}= i_{\mathcal{M}_{\mathcal{O}_\mu}}^*
\omega_{\mathcal{O}_\mu} $; that is, the symplectic form
$\omega_{\mathcal{M}_{\mathcal{O}_\mu}}$ is induced from the $R_o$-reduced symplectic
form $\omega_{\mathcal{O}_\mu}$ on $(T^* Q)_{\mathcal{O}_\mu}$ given in (4.5),
where $i_{\mathcal{M}_{\mathcal{O}_\mu}}^*:
T^*(T^*Q)_{\mathcal{O}_\mu} \rightarrow T^*\mathcal{M}_{\mathcal{O}_\mu}.$
Moreover, the distribution $\mathcal{F}$ is pushed down to a distribution
$\mathcal{F}_{\mathcal{O}_\mu}= T\pi_{\mathcal{O}_\mu}\cdot \mathcal{F}$
on $(T^\ast Q)_{\mathcal{O}_\mu}$,
and we define $\mathcal{K}_{\mathcal{O}_\mu}=\mathcal{F}_{\mathcal{O}_\mu} \cap
T\mathcal{M}_{\mathcal{O}_\mu}$. Assume that $\omega_{\mathcal{K}_{\mathcal{O}_\mu}}=
\tau_{\mathcal{K}_{\mathcal{O}_\mu}}\cdot \omega_{\mathcal{M}_{\mathcal{O}_\mu}}$ is the
restriction of the symplectic form $\omega_{\mathcal{M}_{\mathcal{O}_\mu}}$ on
$T^*\mathcal{M}_{\mathcal{O}_\mu}$ fibrewise to
the distribution $\mathcal{K}_{\mathcal{O}_\mu}$,
where $\tau_{\mathcal{K}_{\mathcal{O}_\mu}}$ is the restriction map to distribution
$\mathcal{K}_{\mathcal{O}_\mu}$. The distributional two-form
$\omega_{\mathcal{K}_{\mathcal{O}_\mu}}$ is not
a "true two-form" on a manifold, which is called
as a $\mathbf{J}$-nonholonomic $R_o$-reduced
distributional two-form to avoid any confusion.\\

From the above construction we know that,
if the admissibility condition $\mathrm{dim}\mathcal{M}_{\mathcal{O}_\mu}=
\mathrm{rank}\mathcal{F}_{\mathcal{O}_\mu}$ and the compatibility condition
$T\mathcal{M}_{\mathcal{O}_\mu}\cap \mathcal{F}_{\mathcal{O}_\mu}^\bot= \{0\}$ hold, where
$\mathcal{F}_{\mathcal{O}_\mu}^\bot$ denotes the symplectic orthogonal of
$\mathcal{F}_{\mathcal{O}_\mu}$ with respect to the $R_o$-reduced symplectic form
$\omega_{\mathcal{O}_\mu}$, then
$\omega_{\mathcal{K}_{\mathcal{O}_\mu}}$ is
non-degenerate as a bilinear form on each fibre of
$\mathcal{K}_{\mathcal{O}_\mu}$, and hence there exists a vector field
$X_{\mathcal{K}_{\mathcal{O}_\mu}}$
on $\mathcal{M}_{\mathcal{O}_\mu}$, which takes values in the constraint
distribution $\mathcal{K}_{\mathcal{O}_\mu}$,
such that for the function $h_{\mathcal{K}_{\mathcal{O}_\mu}}$,
the $\mathbf{J}$-nonholonomic $R_o$-reduced distributional
Hamiltonian equation holds, that is,
\begin{align}
\mathbf{i}_{X_{\mathcal{K}_{\mathcal{O}_\mu}}}\omega_{\mathcal{K}_{\mathcal{O}_\mu}}
=\mathbf{d}h_{\mathcal{K}_{\mathcal{O}_\mu}},
\label{6.8} \end{align}
where $\mathbf{d}h_{\mathcal{K}_{\mathcal{O}_\mu}}$ is the restriction
of $\mathbf{d}h_{\mathcal{M}_{\mathcal{O}_\mu}}$ to $\mathcal{K}_{\mathcal{O}_\mu}$, and
the function $h_{\mathcal{K}_{\mathcal{O}_\mu}}$ satisfies
$\mathbf{d}h_{\mathcal{K}_{\mathcal{O}_\mu}}= \tau_{\mathcal{K}_{\mathcal{O}_\mu}}\cdot \mathbf{d}h_{\mathcal{M}_{\mathcal{O}_\mu}} $,
and $h_{\mathcal{M}_{\mathcal{O}_\mu}}
= \tau_{\mathcal{M}_{\mathcal{O}_\mu}}\cdot h_{\mathcal{O}_\mu}$ is the
restriction of $h_{\mathcal{O}_\mu}$ to $\mathcal{M}_{\mathcal{O}_\mu}$,
and $h_{\mathcal{O}_\mu}$ is the $R_o$-reduced
Hamiltonian function $h_{\mathcal{O}_\mu}: (T^* Q)_{\mathcal{O}_\mu} \rightarrow \mathbb{R}$ defined
by $h_{\mathcal{O}_\mu}\cdot \pi_{\mathcal{O}_\mu}= H\cdot i_{\mathcal{O}_\mu}$.
In addition, from the distributional Hamiltonian equation; that is,
$\mathbf{i}_{X_\mathcal{K}}\omega_\mathcal{K}=\mathbf{d}H_\mathcal
{K},$ we have that $X_{\mathcal{K}}=\tau_{\mathcal{K}}\cdot X_H, $
and from the $\mathbf{J}$-nonholonomic $R_o$-reduced distributional Hamiltonian equation (6.8)
$\mathbf{i}_{X_{\mathcal{K_{\mathcal{O}_\mu}}}}\omega_{\mathcal{K_{\mathcal{O}_\mu}}}
=\mathbf{d}h_{\mathcal{K_{\mathcal{O}_\mu}}}$, we have that
$X_{\mathcal{K_{\mathcal{O}_\mu}}}
=\tau_{\mathcal{K_{\mathcal{O}_\mu}}}\cdot X_{h_{\mathcal{K_{\mathcal{O}_\mu}}}},$
where $ X_{h_{\mathcal{K_{\mathcal{O}_\mu}}}}$ is the Hamiltonian vector field of
the function $h_{\mathcal{K_{\mathcal{O}_\mu}}},$
and the vector fields $X_{\mathcal{K}}$
and $X_{\mathcal{K_{\mathcal{O}_\mu}}}$ are $\pi_{\mathcal{O}_\mu}$-related,
that is, $X_{\mathcal{K_{\mathcal{O}_\mu}}}\cdot \pi_{\mathcal{O}_\mu}
=T\pi_{\mathcal{O}_\mu}\cdot X_{\mathcal{K}}.$ \\

Moreover, assume that the vector fields $F_\mathcal{K}$ and $u_\mathcal{K}$
on $\mathcal{M}$ are pushed down to the vector fields
$f_{\mathcal{M}_{\mathcal{O}_\mu}}=T\pi_{\mathcal{O}_\mu} \cdot F_\mathcal{K}$ and
$u_{\mathcal{M}_{\mathcal{O}_\mu}}=T\pi_{\mathcal{O}_\mu} \cdot u_\mathcal{K}$ on $\mathcal{M}_{\mathcal{O}_\mu}$.
Then we can define that $f_{\mathcal{K}_{\mathcal{O}_\mu}}
=T\pi_{\mathcal{K}_{\mathcal{O}_\mu}}\cdot f_{\mathcal{M}_{\mathcal{O}_\mu}}$ and
$u_{\mathcal{K}_{\mathcal{O}_\mu}}
=T\pi_{\mathcal{K}_{\mathcal{O}_\mu}}\cdot u_{\mathcal{M}_{\mathcal{O}_\mu}};$
that is, $f_{\mathcal{K}_{\mathcal{O}_\mu}}$ and
$u_{\mathcal{K}_{\mathcal{O}_\mu}}$ are the restrictions of
$f_{\mathcal{M}_{\mathcal{O}_\mu}}$ and $u_{\mathcal{M}_{\mathcal{O}_\mu}}$ to $\mathcal{K}_{\mathcal{O}_\mu}$.
In consequence, the 5-tuple $(\mathcal{K}_{\mathcal{O}_\mu},\omega_{\mathcal{K}_{\mathcal{O}_\mu}},
h_{\mathcal{K}_{\mathcal{O}_\mu}}, f_{\mathcal{K}_{\mathcal{O}_\mu}},
u_{\mathcal{K}_{\mathcal{O}_\mu}})$
is a $\mathbf{J}$-nonholonomic $R_o$-reduced distributional RCH system of the nonholonomic
RCH system with symmetry and momentum map
$(T^*Q,G,\omega,\mathbf{J},\mathcal{D},H,F,W)$,
as well as with a control law $u \in W$.
Thus, the geometrical formulation
of the $\mathbf{J}$-nonholonomic $R_o$-reduced distributional RCH
system may be summarized as follows:

\begin{defi} ($\mathbf{J}$-Nonholonomic $R_o$-reduced Distributional RCH System)
Assume that the 8-tuple $(T^*Q,G,\omega,\mathbf{J},\mathcal{D},H, F, W)$
is a nonholonomic RCH system with symmetry and momentum map,
where $\omega$ is the canonical symplectic form on $T^* Q$,
and $\mathcal{D}\subset TQ$ is a $\mathcal{D}$-completely and
$\mathcal{D}$-regularly nonholonomic constraint of the system, and
$\mathcal{D}$, $H, F$ and $W$ are all $G$-invariant. For a regular value
$\mu\in\mathfrak{g}^\ast$ of the momentum map $\mathbf{J}$,
$\mathcal{O}_\mu=G\cdot \mu\subset \mathfrak{g}^\ast$ is the
$G$-orbit of the coadjoint $G$-action through the point $\mu$,
assume that there exists a
$\mathbf{J}$-nonholonomic $R_o$-reduced distribution
$\mathcal{K}_{\mathcal{O}_\mu}$, an associated non-degenerate
and $\mathbf{J}$-nonholonomic $R_o$-reduced distributional two-form
$\omega_{\mathcal{K}_{\mathcal{O}_\mu}}$
and a vector field $X_{\mathcal {K}_{\mathcal{O}_\mu}}$
on the $\mathbf{J}$-nonholonomic $R_o$-reduced constraint submanifold
$\mathcal{M}_{\mathcal{O}_\mu}=(\mathcal{M}\cap \mathbf{J}^{-1}(\mathcal{O}_\mu)) /G, $
where $\mathcal{M}=\mathcal{F}L(\mathcal{D}),$ and $\mathcal{M}\cap
\mathbf{J}^{-1}({\mathcal{O}_\mu})\neq \emptyset, $ such that the
$\mathbf{J}$-nonholonomic $R_o$-reduced distributional Hamiltonian
equation (6.8) holds, that is, $\mathbf{i}_{X_{\mathcal{K}_{\mathcal{O}_\mu}}}\omega_{\mathcal{K}_{\mathcal{O}_\mu}} =
\mathbf{d}h_{\mathcal{K}_{\mathcal{O}_\mu}}$, where
$\mathbf{d}h_{\mathcal{K}_{\mathcal{O}_\mu}}$ is the restriction of
$\mathbf{d}h_{\mathcal{M}_{\mathcal{O}_\mu}}$ to $\mathcal{K}_{\mathcal{O}_\mu}$,
and the function $h_{\mathcal{K}_{\mathcal{O}_\mu}}$,
and the vector fields $f_{\mathcal{K}_{\mathcal{O}_\mu}}$ and
$u_{\mathcal{K}_{\mathcal{O}_\mu}}$ are defined above.
Then the 5-tuple
$(\mathcal{K}_{\mathcal{O}_\mu},\omega_{\mathcal {K}_{\mathcal{O}_\mu}},
h_{\mathcal {K}_{\mathcal{O}_\mu}}, f_{\mathcal {K}_{\mathcal{O}_\mu}},
u_{\mathcal {K}_{\mathcal{O}_\mu}})$ is called a
$\mathbf{J}$-nonholonomic $R_o$-reduced distributional RCH system
of the nonholonomic RCH system with symmetry and momentum map
$(T^*Q,G,\omega,\mathbf{J},\mathcal{D},H, F, u)$
with a control law $u \in W$, and $X_{\mathcal
{K}_{\mathcal{O}_\mu}}$ is called the $\mathbf{J}$-nonholonomic
$R_o$-reduced dynamical vector field.
Denote that
\begin{align}
\hat{X}_{\mathcal{O}_\mu}=X_{(\mathcal{K}_{\mathcal{O}_\mu},
\omega_{\mathcal{K}_{\mathcal{O}_\mu}}, h_{\mathcal {K}_{\mathcal{O}_\mu}},
f_{\mathcal {K}_{\mathcal{O}_\mu}}, u_{\mathcal {K}_{\mathcal{O}_\mu}})}
=X_{\mathcal{K}_{\mathcal{O}_\mu}}+ f_{\mathcal{K}_{\mathcal{O}_\mu}}+u_{\mathcal{K}_{\mathcal{O}_\mu}}
\label{6.9} \end{align}
is the dynamical vector field of the
$\mathbf{J}$-nonholonomic $R_o$-reduced distributional RCH system
$(\mathcal{K}_{\mathcal{O}_\mu},\omega_{\mathcal{K}_{\mathcal{O}_\mu}},\\
h_{\mathcal {K}_{\mathcal{O}_\mu}},
f_{\mathcal {K}_{\mathcal{O}_\mu}}, u_{\mathcal {K}_{\mathcal{O}_\mu}})$,
which is the synthesis of the $\mathbf{J}$-nonholonomic
$R_o$-reduced dynamical vector field $X_{\mathcal {K}_{\mathcal{O}_\mu}}$ and
the vector fields $f_{\mathcal {K}_{\mathcal{O}_\mu}}$
and $u_{\mathcal {K}_{\mathcal{O}_\mu}}$.
Under the above circumstances, we refer to
$(T^*Q,G,\omega,\mathbf{J},\mathcal{D},H, F, u)$ as
a $\mathbf{J}$-nonholonomic regular orbit reducible RCH system
with the associated distributional RCH system
$(\mathcal{K},\omega_{\mathcal {K}},H_{\mathcal{K}}, F_{\mathcal{K}}, u_{\mathcal{K}})$
and the $\mathbf{J}$-nonholonomic $R_o$-reduced
distributional RCH system
$(\mathcal{K}_{\mathcal{O}_\mu},\omega_{\mathcal{K}_{\mathcal{O}_\mu}},
h_{\mathcal{K}_{\mathcal{O}_\mu}},
f_{\mathcal {K}_{\mathcal{O}_\mu}}, u_{\mathcal {K}_{\mathcal{O}_\mu}})$.
\end{defi}

 For a given $\mathbf{J}$-nonholonomic regular orbit reducible RCH system
 $(T^*Q,G,\omega,\mathbf{J},\mathcal{D},H, F, u)$ with the associated distributional RCH system
$(\mathcal{K},\omega_{\mathcal {K}},H_{\mathcal{K}}, F_{\mathcal{K}}, u_{\mathcal{K}})$
and the $\mathbf{J}$-nonholonomic $R_o$-reduced distributional RCH system
$(\mathcal{K}_{\mathcal{O}_\mu},\omega_{\mathcal {K}_{\mathcal{O}_\mu}},
h_{\mathcal{K}_{\mathcal{O}_\mu}}, f_{\mathcal {K}_{\mathcal{O}_\mu}},
u_{\mathcal {K}_{\mathcal{O}_\mu}})$,
by using Lemma 2.3 and 6.3, we can derive precisely
the geometric constraint conditions of the $\mathbf{J}$-nonholonomic $R_o$-reduced distributional two-form
$\omega_{\mathcal{K}_{\mathcal{O}_\mu}}$ for the nonholonomic reducible dynamical vector field;
that is, the two types of Hamilton-Jacobi equations for the
$\mathbf{J}$-nonholonomic $R_o$-reduced distributional RCH system
$(\mathcal{K}_{\mathcal{O}_\mu},\omega_{\mathcal
{K}_{\mathcal{O}_\mu}},h_{\mathcal{K}_{\mathcal{O}_\mu}},
f_{\mathcal {K}_{\mathcal{O}_\mu}}, u_{\mathcal {K}_{\mathcal{O}_\mu}})$.\\

In the following for any $G$-invariant symplectic map $\varepsilon: T^* Q \rightarrow T^* Q $,
we can give the Type II Hamilton-Jacobi theorem for the
$\mathbf{J}$-nonholonomic $R_o$-reduced distributional RCH system.
For convenience, the maps involved in the following
theorem are shown in Diagram-7.
\begin{center}
\hskip 0cm \xymatrix{ \mathbf{J}^{-1}(\mathcal{O}_\mu) \ar[r]^{i_{\mathcal{O}_\mu}} & T^* Q
\ar[d]_{X_{H\cdot \varepsilon}} \ar[dr]^{\tilde{X}^\varepsilon} \ar[r]^{\pi_Q}
& Q \ar[r]^{\gamma} & T^*Q \ar[d]_{\tilde{X}} \ar[dr]^{X_{h_{\mathcal{K}_{\mathcal{O}_\mu}}
\cdot\bar{\varepsilon}}} \ar[r]^{\pi_{\mathcal{O}_\mu}}
& (T^* Q)_{\mathcal{O}_\mu} \ar[d]^{X_{h_{\mathcal{K}_{\mathcal{O}_\mu}}}}
& \mathcal{M}_{\mathcal{O}_\mu} \ar[l]_{i_{\mathcal{M}_{\mathcal{O}_\mu}}}
\ar[d]_{X_{\mathcal{K}_{\mathcal{O}_\mu}}}\\
& T(T^*Q)  & TQ \ar[l]^{T\gamma} & T(T^*Q) \ar[l]^{T\pi_Q} \ar[r]_{T\pi_{\mathcal{O}_\mu}}
& T(T^* Q)_{\mathcal{O}_\mu} \ar[r]^{\tau_{\mathcal{K}_{\mathcal{O}_\mu}}} & \mathcal{K}_{\mathcal{O}_\mu} }
\end{center}
$$\mbox{Diagram-7}$$

\begin{theo} (Type II Hamilton-Jacobi Theorem for the $\mathbf{J}$-Nonholonomic
$R_o$-reduced Distributional RCH system) For a given
$\mathbf{J}$-nonholonomic regular orbit reducible RCH system
$(T^*Q,G,\omega,\mathbf{J},\\ \mathcal{D},H, F, u)$ with the associated distributional RCH system
$(\mathcal{K},\omega_{\mathcal {K}},H_{\mathcal{K}}, F_{\mathcal{K}}, u_{\mathcal{K}})$
and the $\mathbf{J}$-nonholonomic $R_o$-reduced distributional RCH system
$(\mathcal{K}_{\mathcal{O}_\mu},\omega_{\mathcal{K}_{\mathcal{O}_\mu}},
h_{\mathcal{K}_{\mathcal{O}_\mu}},
f_{\mathcal {K}_{\mathcal{O}_\mu}}, u_{\mathcal {K}_{\mathcal{O}_\mu}})$,
assume that $\gamma:
Q \rightarrow T^*Q$ is a one-form on $Q$, and that $\lambda=\gamma \cdot
\pi_{Q}: T^* Q \rightarrow T^* Q, $ and that for any
symplectic map $\varepsilon:T^* Q \rightarrow T^* Q, $
denote that $\tilde{X}^\varepsilon = T\pi_{Q}\cdot \tilde{X}\cdot \varepsilon$,
where $\tilde{X}=X_{(\mathcal{K},\omega_{\mathcal{K}},
H_{\mathcal{K}}, F_{\mathcal{K}}, u_{\mathcal{K}})}
=X_\mathcal {K}+ F_{\mathcal{K}}+u_{\mathcal{K}}$
is the dynamical vector field of the distributional RCH system
$(\mathcal{K},\omega_{\mathcal{K}},
H_{\mathcal{K}}, F_{\mathcal{K}}, u_{\mathcal{K}})$
corresponding to the $\mathbf{J}$-nonholonomic regular orbit reducible RCH
system with symmetry and a momentum map $(T^*Q,G,\omega,\mathbf{J},\mathcal{D},H, F,u)$.
Moreover, assume that $\mu\in\mathfrak{g}^\ast$ is a regular value of the momentum
map $\mathbf{J}$, and that $\textmd{Im}(\gamma)\subset \mathcal{M} \cap
\mathbf{J}^{-1}(\mu), $ and that $\gamma$
and $\varepsilon$ are $G$-invariant, and that
$\varepsilon(\mathcal{M}\cap \mathbf{J}^{-1}(\mathcal{O}_\mu))
\subset \mathcal{M}\cap \mathbf{J}^{-1}(\mathcal{O}_\mu). $ Denote that
$\bar{\gamma}_{\mathcal{O}_\mu}
=\pi_{\mathcal{O}_\mu}(\gamma): Q \rightarrow \mathcal{M}_{\mathcal{O}_\mu} $,
and that $ \textmd{Im}(T\bar{\gamma}_{\mathcal{O}_\mu})\subset \mathcal{K}_{\mathcal{O}_\mu}, $
and that
$\bar{\lambda}_{\mathcal{O}_\mu}=\pi_{\mathcal{O}_\mu}(\lambda): \mathcal{M}\cap \mathbf{J}^{-1}(\mathcal{O}_\mu)
(\subset T^* Q) \rightarrow \mathcal{M}_{\mathcal{O}_\mu}, $ and that
$\bar{\varepsilon}_{\mathcal{O}_\mu}=\pi_{\mathcal{O}_\mu}(\varepsilon):
\mathcal{M} \cap \mathbf{J}^{-1}(\mathcal{O}_\mu) (\subset T^* Q) \rightarrow
\mathcal{M}_{\mathcal{O}_\mu}. $
Then $\varepsilon$ and $\bar{\varepsilon}_{\mathcal{O}_\mu}$ satisfy the
equation $\tau_{\mathcal{K}_{\mathcal{O}_\mu}} \cdot T\bar{\varepsilon}(X_{h_{\mathcal{K}_{\mathcal{O}_\mu}}\cdot \bar{\varepsilon}_{\mathcal{O}_\mu}})
= T\bar{\lambda}_{\mathcal{O}_\mu} \cdot \tilde{X} \cdot\varepsilon $ if and only if
they satisfy the equation $T\bar{\gamma}_{\mathcal{O}_\mu}\cdot
\tilde{X}^\varepsilon= X_{\mathcal{K}_{\mathcal{O}_\mu}}\cdot \bar{\varepsilon}_{\mathcal{O}_\mu}. $
Here $X_{h_{\mathcal{K}_{\mathcal{O}_\mu}} \cdot\bar{\varepsilon}_{\mathcal{O}_\mu}}$
is the Hamiltonian vector field of the
function $h_{\mathcal{K}_{\mathcal{O}_\mu}}\cdot \bar{\varepsilon}_{\mathcal{O}_\mu}: T^* Q\rightarrow \mathbb{R}, $
and $X_{\mathcal{K}_{\mathcal{O}_\mu}}$ is the $\mathbf{J}$-nonholonomic
$R_o$-reduced dynamical vector field. The equation
$T\bar{\gamma}_{\mathcal{O}_\mu}\cdot \tilde{X}^\varepsilon= X_{\mathcal{K}_{\mathcal{O}_\mu}}\cdot
\bar{\varepsilon}_{\mathcal{O}_\mu},$ is called the Type II Hamilton-Jacobi equation for the
$\mathbf{J}$-nonholonomic $R_o$-reduced distributional RCH system
$(\mathcal{K}_{\mathcal{O}_\mu},\omega_{\mathcal{K}_{\mathcal{O}_\mu}},
h_{\mathcal{K}_{\mathcal{O}_\mu}},
f_{\mathcal {K}_{\mathcal{O}_\mu}}, u_{\mathcal {K}_{\mathcal{O}_\mu}})$.
\end{theo}
See the proof and the more details in Wang \cite{wa22a}.\\

For a given $\mathbf{J}$-nonholonomic regular orbit reducible RCH system
$(T^*Q,G,\omega,\mathbf{J},\mathcal{D},H, F, W)$ with the associated distributional RCH system
$(\mathcal{K},\omega_{\mathcal {K}},H_{\mathcal{K}}, F_{\mathcal{K}}, u_{\mathcal{K}})$
and the $\mathbf{J}$-nonholonomic $R_o$-reduced distributional RCH system
$(\mathcal{K}_{\mathcal{O}_\mu},\omega_{\mathcal{K}_{\mathcal{O}_\mu}},
h_{\mathcal{K}_{\mathcal{O}_\mu}},
f_{\mathcal {K}_{\mathcal{O}_\mu}}, u_{\mathcal {K}_{\mathcal{O}_\mu}})$,
we know that the nonholonomic dynamical vector field
$X_{\mathcal{K}}$ and the $\mathbf{J}$-nonholonomic $R_o$-reduced dynamical vector field
$X_{\mathcal{K}_{\mathcal{O}_\mu}}$ are $\pi_{\mathcal{O}_\mu}$-related;
that is, $X_{\mathcal{K}_{\mathcal{O}_\mu}}\cdot \pi_{\mathcal{O}_\mu}
=T\pi_{\mathcal{O}_\mu}\cdot X_{\mathcal{K}}\cdot i_{\mathcal{O}_\mu}.$ Then
we can prove the following Theorem 6.10, which states the
relationship between the solutions of the Type II Hamilton-Jacobi equations and
the $\mathbf{J}$-nonholonomic regular orbit reduction.

\begin{theo}
For a given $\mathbf{J}$-nonholonomic regular orbit reducible RCH system
$(T^*Q, G, \omega,\mathbf{J}, \mathcal{D}, \\ H, F, W)$ with the associated distributional RCH system
$(\mathcal{K},\omega_{\mathcal {K}},H_{\mathcal{K}}, F_{\mathcal{K}}, u_{\mathcal{K}})$
and the $\mathbf{J}$-nonholonomic $R_o$-reduced distributional RCH system
$(\mathcal{K}_{\mathcal{O}_\mu},\omega_{\mathcal{K}_{\mathcal{O}_\mu}},
h_{\mathcal{K}_{\mathcal{O}_\mu}},
f_{\mathcal {K}_{\mathcal{O}_\mu}}, u_{\mathcal {K}_{\mathcal{O}_\mu}})$,
assume that $\gamma: Q \rightarrow T^*Q$ is a one-form on $Q$, and that
$\varepsilon: T^* Q \rightarrow T^* Q $ is a symplectic map,
and that $\bar{\gamma}_{\mathcal{O}_\mu}=\pi_{\mathcal{O}_\mu}(\gamma): Q \rightarrow \mathcal{M}_{\mathcal{O}_\mu} $,
and that $\bar{\varepsilon}_{\mathcal{O}_\mu}=\pi_{\mathcal{O}_\mu}(\varepsilon):
\mathcal{M}\cap \mathbf{J}^{-1}(\mathcal{O}_\mu) (\subset T^*Q) \rightarrow (T^* Q)_{\mathcal{O}_\mu} $.
Under the hypotheses and the notations of Theorem 6.9, then we have that $\varepsilon$
is a solution of the Type II Hamilton-Jacobi equation $T\gamma\cdot
\tilde{X}^\varepsilon= X_{\mathcal{K}}\cdot \varepsilon $ for the distributional
RCH system $(\mathcal{K},\omega_{\mathcal {K}},H_{\mathcal {K}},
F_{\mathcal {K}}, u_{\mathcal {K}})$ if and
only if $\varepsilon$ and $\bar{\varepsilon}_{\mathcal{O}_\mu} $ satisfy the Type II Hamilton-Jacobi
equation $T\bar{\gamma}_{\mathcal{O}_\mu}\cdot \tilde{X}^\varepsilon=
X_{\mathcal{K}_{\mathcal{O}_\mu}}\cdot \bar{\varepsilon}_{\mathcal{O}_\mu} $ for the
$\mathbf{J}$-nonholonomic $R_o$-reduced distributional RCH system
$(\mathcal{K}_{\mathcal{O}_\mu},\omega_{\mathcal{K}_{\mathcal{O}_\mu}},
h_{\mathcal{K}_{\mathcal{O}_\mu}},
f_{\mathcal {K}_{\mathcal{O}_\mu}}, u_{\mathcal {K}_{\mathcal{O}_\mu}})$.
\end{theo}
See the proof and the more details in Wang \cite{wa22a}.\\

\begin{rema}
It is worthy of noting that
the Type I Hamilton-Jacobi equation
$T\bar{\gamma}_{\mathcal{O}_\mu}\cdot \tilde{X}^\gamma
= X_{\mathcal{K}_{\mathcal{O}_\mu}}\cdot \bar{\gamma}_{\mathcal{O}_\mu}$
is the equation of the $\mathbf{J}$-nonholonomic
$R_o$-reduced differential one-form $\bar{\gamma}_{\mathcal{O}_\mu}$; and that
the Type II Hamilton-Jacobi equation
$T\bar{\gamma}_{\mathcal{O}_\mu}\cdot \tilde{X}^\varepsilon
= X_{\mathcal{K}_{\mathcal{O}_\mu}}\cdot \bar{\varepsilon}_{\mathcal{O}_\mu} $
is the equation of the symplectic diffeomorphism map $\varepsilon$
and the $\mathbf{J}$-nonholonomic $R_o$-reduced
symplectic diffeomorphism map $\bar{\varepsilon}_{\mathcal{O}_\mu}. $
If a $\mathbf{J}$-nonholonomic regular point and regular orbit reducible RCH systems we considered
$(T^*Q,G,\omega, \mathbf{J}, \mathcal{D}, H, F, u)$ have not any constrains,
in this case, the $\mathbf{J}$-nonholonomic $R_p$-reduced and $R_o$-reduced distributional
RCH systems are just the $R_p$-reduced RCH system and the $R_o$-reduced RCH system.
From the Type I and Type II Hamilton-Jacobi theorems; that is,
Theorem 6.7 and Theorem 6.9, we can get the Theorem 5.3
and Theorem 5.6, which were given in Wang \cite{wa13d}.
It shows that Theorem 6.7 and Theorem 6.9 can be regarded as an extension of the two types of
Hamilton-Jacobi theorems for the $R_p$-reduced RCH system and the $R_o$-reduced RCH system
to that for the systems with the nonholonomic context.
If the $\mathbf{J}$-nonholonomic regular orbit reducible RCH system we considered
$(T^*Q,G,\omega, \mathbf{J}, \mathcal{D}, H, F, u)$
has not any the external force and the control; that is, $F=0 $ and $u=0$,
in this case, the $\mathbf{J}$-nonholonomic regular orbit reducible RCH system
is just the $\mathbf{J}$-nonholonomic regular orbit reducible
Hamiltonian system $(T^*Q,G,\omega,\mathbf{J},\mathcal{D},H)$.
and with the canonical symplectic form $\omega$ on $T^*Q$.
From the Type II Hamilton-Jacobi theorem; that is,
Theorem 6.9, we can get the Theorem 5.10 in Le\'{o}n and Wang \cite{lewa15}.
It shows that Theorem 6.9 can be regarded as an extension of the Type II
Hamilton-Jacobi theorem for the $\mathbf{J}$-nonholonomic regular orbit reducible Hamiltonian system
to that for the system with the external force and the control.
In particular, in this case,
if the $\mathbf{J}$-nonholonomic regular orbit reducible RCH system we considered has not any constrains;
that is, $F=0, \; u=0 $ and $\mathcal{D}=\emptyset$, then
the $\mathbf{J}$-nonholonomic regular orbit reducible RCH system
is just a regular orbit reducible Hamiltonian system $(T^*Q,G,\omega,H)$
with the canonical symplectic form $\omega$ on $T^*Q$,
we can obtain the Type II Hamilton-Jacobi
equation for the associated $R_o$-reduced Hamiltonian system,
which is given in Wang \cite{wa17}.
Thus, Theorem 6.9 can be regarded as an extension of the Type II Hamilton-Jacobi
theorem for a regular orbit reducible Hamiltonian system to that for the system
with the external force,
the control and the nonholonomic constrain.
\end{rema}

\section{Controlled Magnetic Hamiltonian System and Hamilton-Jacobi Equations}

A magnetic Hamiltonian system is a canonical Hamiltonian system
coupling the action of a magnetic field, so we can introduce a
magnetic symplectic form and use it to define a magnetic
Hamilton's equatins to describe the magnetic Hamiltonian system, and
we drive precisely the geometric constraint conditions of
the magnetic symplectic form for the magnetic Hamiltonian vector field;
that is, the Type I and Type II Hamilton-Jacobi equations.
Moreover, for the magnetic Hamiltonian system with a nonholonomic constraint,
we can give a distributional magnetic Hamiltonian system,
then derive its two types of Hamilton-Jacobi equations,
and we generalize the above results
to nonholonomic reducible magnetic Hamiltonian system with symmetry
(see Wang \cite{wa21c}).\\

A controlled magnetic Hamiltonian (CMH) system is a
regular controlled Hamiltonian (RCH) system
coupling the action of a magnetic field, which is regarded as
a magnetic Hamiltonian system with the external force and the control.
By using the notation of vertical lift map of a vector along a fiber
(see Marsden \textit{et al.} \cite{mawazh10} and Wang \cite{wa21b}),
we give a good expression of the dynamical vector field of the CMH system,
such that we can describe the magnetic vanishing condition
and the CMH-equivalence,
and derive precisely the geometric constraint conditions
of the magnetic symplectic form for the
dynamical vector field of the CMH system;
that is, the Type I and Type II Hamilton-Jacobi equations,
which are an extension of the two types of
Hamilton-Jacobi equations for the magnetic Hamiltonian system
to that for the system with the external force and the control,
and we prove that the
CMH-equivalence for the CMH systems leaves the
solutions of corresponding to Hamilton-Jacobi equations invariant,
if the associated magnetic Hamiltonian systems are equivalent.
Moreover, we consider the CMH system with the nonholonomic constraint,
and derive a distributional CMH system, which is
determined by a non-degenerate distributional two-form induced
from the magnetic symplectic form. Then
we generalize the above results for the nonholonomic
reducible CMH system with symmetry,
and prove two types of Hamilton-Jacobi theorems
for the nonholonomic reduced distributional CMH system
(see Wang \cite{wa21b}).
These research works reveal the deeply internal
relationships of the magnetic symplectic forms,
the nonholonomic constraints, the dynamical
vector fields and the controls of the CMH systems.

\subsection{Magnetic Hamiltonian System}

Let $Q$ be an $n$-dimensional smooth manifold and $TQ$
the tangent bundle, $T^* Q$ the cotangent bundle with a canonical
symplectic form $\omega$ and the projection $\pi_Q: T^* Q
\rightarrow Q $ induces the map $\pi^*_{Q}: T^* Q \rightarrow T^*T^*Q. $
We consider the magnetic symplectic form
$\omega^B= \omega- \pi^*_Q B,$ where $\omega$ is the canonical symplectic
form on $T^*Q$, and $B$ is the closed two-form on $Q$,
and the $\pi_Q^*B$ is called a magnetic term on $T^*Q$.
A magnetic Hamiltonian system is a triple $(T^\ast Q,\omega^B,H)$,
which is a Hamiltonian system defined by the
magnetic symplectic form $\omega^B$, that is,
a canonical Hamiltonian system
coupling the action of a magnetic field $B$. For a given Hamiltonian $H$,
the dynamical vector field $X^B_H$, which is called
the magnetic Hamiltonian vector field,
satisfies the magnetic Hamilton's equation, that is,
$\mathbf{i}_{X^B_{H} }\omega^B= \mathbf{d}H $.
In canonical cotangent bundle coordinates, for any $q \in Q
, \; (q,p)\in T^* Q, $ we have that
$$
\omega=\sum^n_{i=1} \mathbf{d}q^i \wedge \mathbf{d}p_i ,
\;\;\;\;\;\; B=\sum^n_{i,j=1}B_{ij}\mathbf{d}q^i \wedge
\mathbf{d}q^j ,\;\;\; \mathbf{d}B=0, $$
$$\omega^B= \omega -\pi_Q^*B=\sum^n_{i=1} \mathbf{d}q^i \wedge
\mathbf{d}p_i- \sum^n_{i,j=1}B_{ij}\mathbf{d}q^i \wedge
\mathbf{d}q^j,
$$
and the magnetic Hamiltonian vector field $X^B_H$ with respect to
the magnetic symplectic form $\omega^B$ can be expressed that
$$
X^B_H= \sum^n_{i=1} (\frac{\partial H}{\partial
p_i}\frac{\partial}{\partial q^i} - \frac{\partial H}{\partial
q^i}\frac{\partial}{\partial p_i})-
\sum^n_{i,j=1}B_{ij}\frac{\partial H}{\partial
p_j}\frac{\partial}{\partial p_i}.
$$
See Marsden \textit{et al.} \cite{mamiorpera07}.\\

In the following we can derive precisely the geometric constraint conditions of
the magnetic symplectic form for the dynamical vector field
of a magnetic Hamiltonian system; that is,
the Type I and Type II Hamilton-Jacobi equations for
the magnetic Hamiltonian system $(T^\ast Q,\omega^B,H)$.
In order to do this, in the following we first give
an important notion and prove a key lemma, which is an important
tool for the proofs of two types of
Hamilton-Jacobi theorem for the magnetic Hamiltonian system.\\

For the one-form $\gamma: Q \rightarrow T^*Q$, $\mathbf{d}\gamma$
is a two-form on $Q$. Assume that $B$ is a closed two-form on $Q$,
then we say that the $\gamma$ satisfies condition $\mathbf{d}\gamma=-B$
if, for any $ x, y \in TQ$, we have that $(\mathbf{d}\gamma +B)(x,y)=0.$
In the following we can give a new notion.
\begin{defi}
Assume that $\gamma: Q
\rightarrow T^*Q$ is a one-form on $Q$,
we say that the $\gamma$ satisfies condition that
$\mathbf{d}\gamma=-B$ with respect to $T\pi_{Q}:
TT^* Q \rightarrow TQ $ if, for any $v, w \in TT^* Q, $ we have
that $(\mathbf{d}\gamma +B)(T\pi_{Q}(v),T\pi_{Q}(w))=0. $
\end{defi}

From the above definition we know that if $\gamma$ satisfies condition
$\mathbf{d}\gamma=-B$, then it must satisfy condition
$\mathbf{d}\gamma=-B$ with respect to $T\pi_{Q}: TT^* Q \rightarrow
TQ. $ Conversely, if $\gamma$ satisfies condition
$\mathbf{d}\gamma=-B$ with respect to
$T\pi_{Q}: TT^* Q \rightarrow TQ, $ then it may not satisfy condition
$\mathbf{d}\gamma=-B$. \\

Now, we can give the following lemma,
which  can be regarded as an extension of the Lemma 2.3 ( given by
Wang \cite{wa17}) to it with the nonholonomic context,
and the lemma is a very important tool for our research.

\begin{lemm}
Assume that $\gamma: Q \rightarrow T^*Q$ is a one-form on $Q$, and
that $\lambda=\gamma \cdot \pi_{Q}: T^* Q \rightarrow T^* Q .$
For the magnetic symplectic form $\omega^B= \omega- \pi_Q^*B $ on $T^*Q$,
where $\omega$ is the canonical symplectic form on $T^*Q$,
and $B$ is a closed two-form on $Q$,
then we have that the following two assertions hold.\\
\noindent $(\mathrm{i})$ For any $v, w \in
TT^* Q, $ we have that $\lambda^*\omega^B(v,w)= -(\mathbf{d}\gamma+B)(T\pi_{Q}(v), \;
T\pi_{Q}(w))$. \\
\noindent $(\mathrm{ii})$ For any $v, w \in TT^* Q, $ we have that \\
$\omega^B(T\lambda \cdot v,w)= \omega^B(v, w-T\lambda \cdot
w)-(\mathbf{d}\gamma+B)(T\pi_{Q}(v), \; T\pi_{Q}(w)). $
\end{lemm}
See the proof and the more details in Wang \cite{wa21b, wa21c}.\\

Usually, under the impact of the
magnetic term $\pi_Q^*B$, the magnetic symplectic form
$\omega^B=\omega- \pi_Q^*B $,
in general, is not the canonical symplectic form $\omega$ on $T^*Q$,
we cannot prove the Hamilton-Jacobi theorem for a magnetic Hamiltonian system
the same as in Theorem 2.1. However, using Lemma 7.2
we can derive precisely the geometric constraint conditions of
the magnetic symplectic form for the dynamical vector field of
the magnetic Hamiltonian system; that is, the Type I and Type II
Hamilton-Jacobi equations for the magnetic Hamiltonian system.
For convenience, the maps involved in the
theorem are shown in Diagram-8.
\begin{center}
\hskip 0cm \xymatrix{ T^* Q \ar[r]^\varepsilon
& T^* Q  \ar[d]_{X^B_{H\cdot \varepsilon}} \ar[dr]^{X^\varepsilon} \ar[r]^{\pi_Q}
 & Q \ar[d]^{X^\gamma} \ar[r]^{\gamma} & T^*Q \ar[d]^{X^B_H} \\
 & T(T^*Q) & TQ \ar[l]^{T\gamma} & T(T^* Q)\ar[l]^{T\pi_Q}}
\end{center}
$$\mbox{Diagram-8}$$

\begin{theo}
(Hamilton-Jacobi Theorem for a Magnetic Hamiltonian System)
For a given magnetic Hamiltonian system $(T^*Q,\omega^B,H)$ with
a magnetic symplectic form $\omega^B= \omega- \pi_Q^*B $ on $T^*Q$,
where $\omega$ is the canonical symplectic form on $T^* Q$
and $B$ is a closed two-form on $Q$,
assume that $\gamma: Q
\rightarrow T^*Q$ is a one-form on $Q$, and
that $\lambda=\gamma \cdot \pi_{Q}: T^* Q \rightarrow T^* Q $, and that for any
symplectic map $\varepsilon: T^* Q \rightarrow T^* Q $ with respect to $\omega^B$,
denote that $X^\gamma = T\pi_{Q}\cdot X^B_H \cdot \gamma$
and $ X^\varepsilon = T\pi_{Q}\cdot X^B_H \cdot \varepsilon$,
where $X^B_H$ is the dynamical vector field
of the magnetic Hamiltonian system $(T^*Q,\omega^B,H)$;
that is, the magnetic Hamiltonian vector field.
Then the following two assertions hold:\\
\noindent $(\mathbf{i})$
If the one-form $\gamma: Q \rightarrow T^*Q $ satisfies the condition that
$\mathbf{d}\gamma=-B $ with respect to $T\pi_{Q}:
TT^* Q \rightarrow TQ, $ then $\gamma$ is a solution of the
Type I Hamilton-Jacobi equation
$T\gamma\cdot X^\gamma= X^B_H\cdot \gamma $
for the magnetic Hamiltonian system $(T^*Q,\omega^B,H)$. \\
\noindent $(\mathbf{ii})$
The $\varepsilon$ is a solution of the Type II
Hamilton-Jacobi equation $T\gamma \cdot X^\varepsilon= X^B_H\cdot
\varepsilon $ if and only if it is a solution of the equation
$T\varepsilon\cdot X^B_{H\cdot\varepsilon}= T\lambda \cdot X^B_H \cdot \varepsilon $,
where $ X^B_{H\cdot\varepsilon} \in
TT^*Q $ is the magnetic Hamiltonian vector field of the function $H\cdot\varepsilon:
T^*Q\rightarrow \mathbb{R} $.
\end{theo}
See the proof and the more details in Wang \cite{wa21c}.
It is worthy of noting that the Type I Hamilton-Jacobi equation
$T\gamma \cdot X^\gamma= X^B_H \cdot \gamma $
is the equation of the differential one-form $\gamma$, and that
the Type II Hamilton-Jacobi equation $T\gamma\cdot X^\varepsilon
= X^B_H \cdot \varepsilon $ is the equation of
the symplectic diffeomorphism map $\varepsilon$ with respect to $\omega^B$.
If $B=0$, in this case the magnetic symplectic form $\omega^B$
is just the canonical symplectic form $\omega$ on $T^*Q$, and the
condition that the one-form $\gamma: Q \rightarrow T^*Q $ satisfies the condition
$\mathbf{d}\gamma=-B $ with respect to $T\pi_{Q}:
TT^* Q \rightarrow TQ, $ becomes that  $\gamma $ is closed with respect to
$T\pi_Q: TT^* Q \rightarrow TQ.$ Thus, from above Theorem 7.3,
we can obtain Theorem 2.5 and Theorem 2.6 in Wang \cite{wa17}.
Thus, Theorem 7.3 can be regarded as an extension of the two types of Hamilton-Jacobi
equations for a canonical Hamiltonian system to that for the system
with the magnetic. \\

In order to describe the impact of the different geometric structures
and the constraints for the dynamics of a Hamiltonian system,
for the magnetic Hamiltonian system with a nonholonomic constraint,
we can define a distributional magnetic Hamiltonian system
by analyzing carefully the structure of the nonholonomic dynamical
vector field, and derive precisely its two types of Hamilton-Jacobi equations.
Moreover, we generalize the above results
to the nonholonomic reducible magnetic Hamiltonian system with symmetry.
We give the definition of a nonholonomic reduced distributional
magnetic Hamiltonian system, and prove the two types of
Hamilton-Jacobi theorems for the system,
which are the development of the Type I and Type II
Hamilton-Jacobi theorems for the nonholonomic reduced distributional
Hamiltonian system given in Le\'{o}n and Wang \cite{lewa15};
See Wang \cite{wa21c}for more details.

\subsection{Controlled Magnetic Hamiltonian System}

In order to describe the impact of the different geometric structures
for the dynamics and the Hamilton-Jacobi equations of an RCH system,
considering the external force and the
control, we can define a kind of controlled magnetic Hamiltonian
(CMH) system on $T^*Q$ as follows.
\begin{defi}
(CMH System) A controlled magnetic Hamiltonian (CMH) system
on $T^*Q$ is a 5-tuple $(T^*Q, \omega^B, H, F, W)$, which is a
magnetic Hamiltonian system $(T^\ast Q,\omega^B,H)$
with the external force $F$ and the control set $W$, where $F: T^*Q\rightarrow T^*Q$ is
the fiber-preserving map, and $W\subset T^*Q$ is a fiber submanifold,
which is called the control subset.
\end{defi}

From the above Definition 2.2 and Definition 7.4 we know that
a CMH system on $T^*Q$ is also an RCH system on $T^*Q$, however, its symplectic structure
is given by a magnetic symplectic form, and the set of
the CMH systems on $T^*Q$ is a subset of the set of the RCH systems on $T^*Q$.
However, the subset of CMH systems is not closed under the
actions of the external force and the control, because the
magnetic may be vanishing. \\

When a feedback control law $u: T^*Q\rightarrow W$ is
chosen, the 5-tuple $(T^*Q, \omega^B, H, F, u)$ is a regular
closed-loop dynamic system.
In order to describe the dynamics of the CMH system
$(T^*Q,\omega^B,H,F,W)$ with a control law $u$,
by using the notation of vertical lift map of a vector along a
fiber; (also see Marsden \textit{et al.} \cite{mawazh10}),
we can give a good expression of the dynamical vector field
of the CMH system as follows:
for a given CMH system $(T^\ast Q, \omega^B, H, F, W)$, the dynamical
vector field of the associated magnetic Hamiltonian system $(T^\ast Q,
\omega^B, H) $ is  $X^B_H$, which satisfies the equation
$\mathbf{i}_{X^B_H}\omega^B=\mathbf{d}H$. If considering the
external force $F: T^*Q \rightarrow T^*Q, $ by using the
notation of vertical lift map of a vector along a fiber, the change
of $X^B_H$ under the action of $F$ is such that
$$\textnormal{vlift}(F)X^B_H(\alpha_x)
= \textnormal{vlift}((TFX^B_H)(F(\alpha_x)), \alpha_x)
= (TFX^B_H)^v_\gamma(\alpha_x),$$
where $\alpha_x \in T^*_x Q, \; x\in Q $ and $\gamma$ is a straight
line in $T^*_x Q$ connecting $F_x(\alpha_x)$ and $\alpha_x$. In the
same way, when a feedback control law $u: T^\ast Q \rightarrow W$ is
chosen, the change of $X^B_H$ under the action of $u$ is such that
$$\textnormal{vlift}(u)X^B_H(\alpha_x)
= \textnormal{vlift}((TuX^B_H)(u(\alpha_x)), \alpha_x)
= (TuX^B_H)^v_\gamma(\alpha_x).$$
As a consequence, we can give an expression of the dynamical vector
field of the CMH system as follows:
\begin{theo}
The dynamical vector field of the CMH system $(T^\ast Q,\omega^B,H,F,W)$
with a control law $u$ is the synthesis of the magnetic Hamiltonian vector field
$X^B_H$ and its changes under the actions of the external force $F$
and the control law $u$; that is,
\begin{equation}
X_{(T^\ast Q,\omega^B,H,F,u)}(\alpha_x)
= X^B_H(\alpha_x)+ \textnormal{vlift}(F)X^B_H(\alpha_x)
+ \textnormal{vlift}(u)X^B_H(\alpha_x), \;\; \label{7.1}
\end{equation}
 for any $\alpha_x \in T^*_x
Q, \; x\in Q $. For convenience, that is simply written as
\begin{equation}X_{(T^\ast Q,\omega^B,H,F,u)}
=X^B_H +\textnormal{vlift}(F)^B +\textnormal{vlift}(u)^B. \;\; \label{7.2}
\end{equation}
\end{theo}
Where $\textnormal{vlift}(F)^B=\textnormal{vlift}(F)X^B_H$,
and $\textnormal{vlift}(u)^B=\textnormal{vlift}(u)X^B_H.$ are the
changes of $X^B_H$ under the actions of $F$ and $u$.
We also denote that $\textnormal{vlift}(W)^B= \bigcup\{\textnormal{vlift}(u)X^B_H |
\; u\in W\}$. \\

From the expression (7.2) of the dynamical vector
field of a CMH system, we know that under the actions of the external force $F$
and the control law $u$, in general, the dynamical vector
field may not be magnetic Hamiltonian, and hence the CMH system may not
be yet a magnetic Hamiltonian system. However,
it is a dynamical system closed relative to a
magnetic Hamiltonian system, and it can be explored and studied by the extending
methods for the external force and the control
in the study of the magnetic Hamiltonian system.\\

For the magnetic Hamiltonian system $(T^\ast Q, \omega^B, H)$, its
magnetic Hamiltonian vector field $X^B_H$ satisfies the equation
$\mathbf{i}_{X^B_H}\omega^B=\mathbf{d}H$, and for the
associated canonical Hamiltonian system $(T^\ast Q,
\omega, H) $, its canonical Hamiltonian vector field $X_H$
satisfies the equation $\mathbf{i}_{X_H}\omega=\mathbf{d}H$.
Denote that the vector field $X^0= X^B_H-X_H, $ and, from the
magnetic symplectic form $\omega^B= \omega- \pi_Q^*B$, we have that
$$
\mathbf{i}_{X^0}\omega=\mathbf{i}_{(X^B_H-X_H)}\omega
=\mathbf{i}_{X^B_H}\omega-\mathbf{i}_{X_H}\omega
=\mathbf{i}_{X^B_H}(\omega^B+ \pi_Q^*B)-\mathbf{i}_{X_H}\omega
=\mathbf{i}_{X^B_H}( \pi_Q^*B).
$$
Thus, $X^0$ is called the magnetic vector field and that
$\mathbf{i}_{X^0}\omega=\mathbf{i}_{X^B_H}( \pi_Q^*B)$
is called the magnetic equation, which is determined by
the magnetic term $\pi_Q^*B$ on $T^*Q$. When $B=0$,
then $X^0=0$, the magnetic equation holds trivially.
For the CMH system $(T^\ast Q, \omega^B, H, F, W)$,
from the expression (7.2) of its dynamical
vector field, we have that
\begin{equation}X_{(T^\ast Q,\omega^B,H,F,u)}
=X_H + X^0+\textnormal{vlift}(F)^B +\textnormal{vlift}(u)^B. \;\; \label{7.3}
\end{equation}
If we choose the external force $F$ and the control law $u$, such that
\begin{equation}
 X^0+\textnormal{vlift}(F)^B +\textnormal{vlift}(u)^B=0, \;\; \label{7.4}
\end{equation}
then from (7.3) we have that $X_{(T^\ast Q,\omega^B,H,F,u)}
=X_H; $ that is, in this case the dynamical vector
field of the CMH system is just the canonical Hamiltonian vector field itself,
and the motion of the CMH system is just the same as the motion of the canonical
Hamiltonian system without the actions of the magnetic, the external force and the control.
Thus, the condition (7.4) is called the magnetic vanishing condition for
the CMH system $(T^\ast Q, \omega^B, H, F, W)$, and we have the following theorem:
\begin{theo}
If the external force $F$ and the control law $u$ for the CMH system
$(T^\ast Q, \omega^B, H, F, u)$ satisfy the magnetic vanishing condition
(7.4), then the dynamical vector field $X_{(T^\ast Q,\omega^B,H,F,u)}$
of the CMH system is just the canonical Hamiltonian vector field $X_H$ for the
associated canonical Hamiltonian system $(T^\ast Q, \omega, H) $.
\end{theo}

For a given CMH system $(T^*Q,\omega^B,H,F,W)$ on $T^*Q$, by using
the above Lemma 7.2, we can derive precisely the geometric constraint
conditions of the magnetic symplectic form $\omega^B$ for the dynamical
vector field $X_{(T^\ast Q,\omega^B,H,F,u)}$ of the CMH system with a control law $u$;
that is, the Type I and Type II
Hamilton-Jacobi equations for the CMH system.
For convenience, the maps involved in the
theorem are shown in Diagram-9.
\begin{center}
\hskip 0cm \xymatrix{ T^* Q \ar[r]^\varepsilon
& T^* Q  \ar[d]_{X^B_{H\cdot \varepsilon}} \ar[dr]^{\tilde{X}^\varepsilon} \ar[r]^{\pi_Q}
 & Q \ar[d]^{\tilde{X}^\gamma} \ar[r]^{\gamma} & T^*Q \ar[d]^{\tilde{X}} \\
 & T(T^*Q) & TQ \ar[l]^{T\gamma} & T(T^* Q) \ar[l]^{T\pi_Q}}
\end{center}
$$\mbox{Diagram-9}$$
\begin{theo}
(Hamilton-Jacobi Theorem for a CMH System)
For a given CMH system $(T^*Q,\omega^B,H)$ with
a magnetic symplectic form $\omega^B= \omega- \pi_Q^*B $ on $T^*Q$,
where $\omega$ is the canonical symplectic form on $T^* Q$
and $B$ is a closed two-form on $Q$,
assume that $\gamma: Q
\rightarrow T^*Q$ is a one-form on $Q$, and
that $\lambda=\gamma \cdot \pi_{Q}: T^* Q \rightarrow T^* Q $, and that for any
symplectic map $\varepsilon: T^* Q \rightarrow T^* Q $ with respect to $\omega^B$,
denote that $\tilde{X}^\gamma = T\pi_{Q}\cdot \tilde{X}\cdot \gamma$
and $ \tilde{X}^\varepsilon = T\pi_{Q}\cdot \tilde{X}\cdot \varepsilon$,
where $\tilde{X}=X_{(T^\ast Q,\omega^B,H,F,u)}$ is the dynamical vector field
of the CMH system $(T^*Q,\omega^B,H,F,W)$ with a control law $u$.
Then the following two assertions hold:\\
\noindent $(\mathbf{i})$
If the one-form $\gamma: Q \rightarrow T^*Q $ satisfies the condition that
$\mathbf{d}\gamma=-B $ with respect to $T\pi_{Q}:
TT^* Q \rightarrow TQ, $ then $\gamma$ is a solution of the
Type I Hamilton-Jacobi equation
$T\gamma \cdot \tilde{X}^\gamma= X^B_H\cdot \gamma $,
where $X^B_H$ is the magnetic Hamiltonian vector field
of the associated magnetic Hamiltonian system $(T^*Q,\omega^B,H).$.\\
\noindent $(\mathbf{ii})$
The $\varepsilon$ is a solution of the Type II
Hamilton-Jacobi equation $T\gamma \cdot \tilde{X}^\varepsilon= X^B_H\cdot
\varepsilon $ if and only if it is a solution of the equation
$T\varepsilon\cdot X^B_{H\cdot\varepsilon}= T\lambda \cdot \tilde{X} \cdot \varepsilon $,
where $X^B_H$ and $ X^B_{H\cdot\varepsilon} \in
TT^*Q $ are the magnetic Hamiltonian vector fields of the functions $H$ and $H\cdot\varepsilon:
T^*Q\rightarrow \mathbb{R}, $ respectively.
\end{theo}
See the proof and the more details in Wang \cite{wa21b}.

\begin{rema}
It is worthy of noting that
the Type I Hamilton-Jacobi equation
$T\gamma\cdot \tilde{X}^\gamma= X^B_H \cdot \gamma $
is the equation of the differential one-form $\gamma$, and that
the Type II Hamilton-Jacobi equation $T\gamma\cdot \tilde{X}^\varepsilon
= X^B_H \cdot \varepsilon $ is the equation of
the symplectic diffeomorphism map $\varepsilon$ with respect to $\omega^B$.
If the external force and control of the CMH
system $(T^*Q,\omega^B,H,F,W)$ are both zeros, that is, $F=0 $
and $W=\emptyset$, in this case the CMH system
is just a magnetic Hamiltonian system $(T^*Q,\omega^B,H)$, and from the
Theorem 7.7, we can obtain the two types of Hamilton-Jacobi
equations for the associated magnetic Hamiltonian system; that is,  Theorem 7.3
(see Wang \cite{wa21c}).
Thus, Theorem 7.7 can be regarded as an extension of the two types
of Hamilton-Jacobi equations for a magnetic Hamiltonian system
to that for the system with the external force and the control.
If $B=0$, in this case the magnetic symplectic form $\omega^B$
is just the canonical symplectic form $\omega$ on $T^*Q$, and the
condition that the one-form $\gamma: Q \rightarrow T^*Q $ satisfies the condition
$\mathbf{d}\gamma=-B $ with respect to $T\pi_{Q}:
TT^* Q \rightarrow TQ $ becomes that  $\gamma $ is closed with respect to
$T\pi_Q: TT^* Q \rightarrow TQ.$ Thus, from the Theorem 7.7,
we can obtain Theorem 5.1 (see Wang \cite{wa13d}). Moreover,
if $B=0$, $F=0 $ and $W=\emptyset$, then, from the Theorem 7.7,
we can obtain Theorem 2.5 and Theorem 2.6 in Wang \cite{wa17}.
Thus, the Theorem 7.7 can be regarded as an extension of two types of Hamilton-Jacobi
equations for a canonical Hamiltonian system to that for the system
with the magnetic, the external force and the control.
\end{rema}

In the same way, for two given magnetic
Hamiltonian systems $(T^\ast Q_i,\omega^B_i,H_i),$ $ i= 1,2,$
we say them to be equivalent, if there exists a
diffeomorphism $\varphi: Q_1\rightarrow Q_2$, which
is symplectic with respect to their magnetic symplectic forms,
such that their magnetic Hamiltonian vector fields $X^B_{H_i}, \; i=1,2 $ satisfy
the condition $X^B_{H_1}\cdot \varphi^\ast =T(\varphi^\ast) X^B_{H_2}$. \\

For two given CMH systems $(T^\ast
Q_i,\omega^B_i,H_i,F_i,W_i),$ $ i= 1,2,$ we also want to define
their equivalence, that is, to look for a diffeomorphism
$\varphi: Q_1\rightarrow Q_2$, such that
$X_{(T^\ast Q_1,\omega^B_1,H_1,F_1,W_1)}\cdot \varphi^\ast
=T(\varphi^\ast) X_{(T^\ast Q_2,\omega^B_2,H_2,F_2,W_2)}$. However,
it is worthy of noting that because the
magnetic may be vanishing and
the set of CMH systems, as a subset of the set of RCH systems,
is not closed under the actions of the external force and the control,
so we cannot use the RCH-equivalence for the CMH systems.
On the other hand, when a CMH system is given,
the force map $F$ is determined,
but the feedback control law $u: T^\ast Q\rightarrow W$
could be chosen. In order to describe the feedback control law to
modify the structures of the CMH systems, the controlled magnetic Hamiltonian matching
conditions and CMH-equivalence are induced as follows:
\begin{defi}
(CMH-equivalence) Suppose that we have two CMH systems $(T^\ast
Q_i,\omega^B_i,H_i,F_i,W_i),$ $ i= 1,2,$ we say that they are
CMH-equivalent, or simply, that $(T^\ast
Q_1,\omega^B_1,H_1,F_1,W_1)\stackrel{CMH}{\sim}\\ (T^\ast
Q_2,\omega^B_2,H_2,F_2,W_2)$, if there exists a
diffeomorphism $\varphi: Q_1\rightarrow Q_2$, such that the
following controlled magnetic Hamiltonian matching conditions hold:

\noindent {\bf CMH-1:} The control subsets $W_i, \; i=1,2$ satisfy
the condition $W_1=\varphi^\ast (W_2),$ where the map
$\varphi^\ast= T^\ast \varphi:T^\ast Q_2\rightarrow T^\ast Q_1$
is cotangent lifted map of $\varphi$.

\noindent {\bf CMH-2:} For each control law $u_1:
T^\ast Q_1 \rightarrow W_1, $ there exists the control law $u_2:
T^\ast Q_2 \rightarrow W_2, $  such that the two
closed-loop dynamical systems produce the same dynamical vector fields; that is,
$X_{(T^\ast Q_1,\omega^B_1,H_1,F_1,u_1)}\cdot \varphi^\ast
=T(\varphi^\ast) X_{(T^\ast Q_2,\omega^B_2,H_2,F_2,u_2)}$,
where the map $T(\varphi^\ast):TT^\ast Q_2\rightarrow TT^\ast Q_1$
is the tangent map of $\varphi^\ast$.
\end{defi}

From the expression (7.1) of the dynamical vector field of the CMH system
and the condition $X_{(T^\ast Q_1,\omega^B_1,H_1,F_1,u_1)}\cdot \varphi^\ast
=T(\varphi^\ast) X_{(T^\ast Q_2,\omega^B_2,H_2,F_2,u_2)}$, we have that
$$
(X^B_{H_1}+\textnormal{vlift}(F_1)X^B_{H_1}+\textnormal{vlift}(u_1)X^B_{H_1})\cdot \varphi^\ast
= T(\varphi^\ast)[X^B_{H_2}+\textnormal{vlift}(F_2)X^B_{H_2}+\textnormal{vlift}(u_2)X^B_{H_2}].
$$
By using the notation of vertical lift map of a vector along a fiber,
for $\alpha_x \in T_x^\ast Q_2, \; x \in Q_2$, we have that
\begin{align*}
T(\varphi^\ast)\textnormal{vlift}(F_2)X^B_{H_2}(\alpha_x)
&=T(\varphi^\ast)\textnormal{vlift}((TF_2X^B_{H_2})(F_2(\alpha_x)), \alpha_x)\\
&=\textnormal{vlift}(T(\varphi^\ast)\cdot TF_2\cdot T(\varphi_\ast)X^B_{H_2}
(\varphi^\ast F_2 \varphi_\ast(\varphi^\ast \alpha_x)),\varphi^\ast \alpha)\\
&=\textnormal{vlift}(T(\varphi^\ast F_2 \varphi_\ast)X^B_{H_2}
(\varphi^\ast F_2 \varphi_\ast(\varphi^\ast \alpha_x)),\varphi^\ast \alpha)\\
&=\textnormal{vlift}(\varphi^\ast
F_2\varphi_\ast)X^B_{H_2}(\varphi^\ast \alpha_x),
\end{align*}
where the map $\varphi_\ast=(\varphi^{-1})^\ast: T^\ast Q_1\rightarrow T^\ast Q_2$.
In the same way, we have that
$T(\varphi^\ast)\textnormal{vlift}(u_2)X^B_{H_2}=\textnormal{vlift}(\varphi^\ast
u_2\varphi_\ast)X^B_{H_2}\cdot \varphi^\ast$. Note that
$\textnormal{vlift}(F)^B=\textnormal{vlift}(F)X^B_H$,
and $\textnormal{vlift}(u)^B=\textnormal{vlift}(u)X^B_H.$,
and hence we have that the explicit relation
between the two control laws $u_i \in W_i, \; i=1,2$ in {\bf CMH-2} is given by
\begin{align}
& (\textnormal{vlift}(u_1)^B -\textnormal{vlift}(\varphi^\ast u_2\varphi_\ast)^B)\cdot \varphi^\ast \nonumber \\
& = -X^B_{H_1}\cdot \varphi^\ast +T(\varphi^\ast) (X^B_{H_2})+
(-\textnormal{vlift}(F_1)^B+\textnormal{vlift}(\varphi^\ast F_2
\varphi_\ast)^B)\cdot \varphi^\ast.
\label{7.5}\end{align}
From the above relation (7.5) we know that when two CMH systems $(T^\ast
Q_i,\omega^B_i,H_i,F_i,W_i),$ $ i= 1,2,$ are CMH-equivalent with respect to $\varphi^*$, the
associated magnetic Hamiltonian systems $(T^\ast Q_i,\omega^B_i,H_i),$
$ i= 1,2,$ may not be equivalent with respect to $\varphi^*$.\\

On the other hand,
note that the magnetic vector field $X^0= X^B_H-X_H, $
from (7.3) and (7.5) we have that
\begin{align*}
& (\textnormal{vlift}(u_1)^B -\textnormal{vlift}(\varphi^\ast u_2\varphi_\ast)^B)\cdot \varphi^\ast \\
& = -(X_{H_1}+X^0_1)\cdot \varphi^\ast +T(\varphi^\ast) (X_{H_2}+X^0_2)+
(-\textnormal{vlift}(F_1)^B+\textnormal{vlift}(\varphi^\ast F_2
\varphi_\ast)^B)\cdot \varphi^\ast.
\end{align*}
and hence, we have that
\begin{align}
& (X^0_1+\textnormal{vlift}(F_1)^B+\textnormal{vlift}(u_1)^B )\cdot \varphi^\ast \nonumber \\
& = -X_{H_1}\cdot \varphi^\ast +T(\varphi^\ast) X_{H_2}+T(\varphi^\ast) (X^0_2
+\textnormal{vlift}(F_2)^B+\textnormal{vlift}(u_2)^B).
\label{7.6}\end{align}
If the associated canonical Hamiltonian systems $(T^\ast
Q_i,\omega_i,H_i),$ $ i= 1,2,$ are also equivalent with respect to $\varphi^*$; that is,
$T(\varphi^*) \cdot X_{H_2}= X_{H_1}\cdot \varphi^*$. In this case,
from (7.4) and (7.6), we have the following theorem.
\begin{theo}
Suppose that two CMH systems $(T^\ast Q_i,\omega^B_i,H_i,F_i,W_i)$,
$i=1,2,$ are CMH-equivalent with respect to $\varphi^*$,
and the associated canonical Hamiltonian systems $(T^\ast
Q_i,\omega_i,H_i),$ $ i= 1,2,$ are also equivalent with respect to $\varphi^*$.
Then we have the following fact that if one system satisfies the magnetic
vanishing condition, then another CMH-equivalent system must satisfy
the associated magnetic vanishing condition.
\end{theo}

Moreover, considering the CMH-equivalence of the CMH systems,
we can obtain the following Theorem 7.11, which states that
the solutions of the two types of Hamilton-Jacobi equations for the CMH systems leave
invariant under the conditions of CMH-equivalence, if the associated
magnetic Hamiltonian systems are equivalent.
\begin{theo}
Suppose that two CMH systems $(T^\ast Q_i,\omega^B_i,H_i,F_i,W_i)$,
$i=1,2,$ are CMH-equivalent with an equivalent map $\varphi: Q_1
\rightarrow Q_2 $, and the associated magnetic Hamiltonian systems $(T^\ast
Q_i,\omega^B_i,H_i),$ $ i= 1,2,$ are also equivalent with respect to $\varphi^*$,
under the hypotheses and the notations of the Theorem 7.7, we have that
the following two assertions hold:\\
\noindent $(\mathrm{i})$ If the one-form $\gamma_2: Q_2 \rightarrow T^* Q_2$
satisfies the condition that $\mathbf{d}\gamma_2=-B_2 $ with
respect to $T\pi_{Q_2}: TT^* Q_2 \rightarrow TQ_2, $ then $\gamma_1=
\varphi^* \cdot \gamma_2\cdot \varphi: Q_1 \rightarrow T^* Q_1 $ satisfies also the condition that
$\mathbf{d}\gamma_1=-B_1 $ with respect to $T\pi_{Q_1}:
TT^* Q_1 \rightarrow TQ_1, $ and hence, it is
a solution of the Type I Hamilton-Jacobi equation for the CMH system
$(T^*Q_1,\omega^B_1,H_1,F_1,W_1). $ Vice versa.

\noindent $(\mathrm{ii})$ If the symplectic map
$\varepsilon_2: T^*Q_2\rightarrow T^* Q_2$ with respect to $\omega^B_2$
is a solution of the Type II Hamilton-Jacobi equation for the CMH system
$(T^*Q_2,\omega^B_2,H_2, F_2,W_2)$, then $\varepsilon_1=
\varphi^* \cdot \varepsilon_2\cdot \varphi_*: T^*Q_1 \rightarrow
T^* Q_1 $ is a symplectic map with respect to $\omega^B_1$,
and it is a solution of the Type II Hamilton-Jacobi equation
for the CMH system $(T^*Q_1,\omega^B_1,H_1,F_1,W_1). $ Vice versa.
\end{theo}
See the proof and the more details in Wang \cite{wa21b}.\\

In order to describe the impact of the different group structure of symmetry
for the dynamics of a reducible RCH system,
we consider the regular point reduction of the CMH system
$(T^\ast Q,\mathcal{H},\omega^B,H,F, W)$ with symmetry of
the Heisenberg group $\mathcal{H}$. Here the configuration space $Q=
\mathcal{H}\times V, \; \mathcal{H}= \mathbb{R}^2\oplus \mathbb{R},
$ and $V$ is a $k$-dimensional vector space, and the cotangent
bundle $T^*Q$ with the magnetic symplectic form $\omega^B= \omega_0-
\pi^*_Q \bar{B},$ where $\omega_0$ is the usual canonical symplectic
form on $T^*Q$, and $\bar{B}= \pi_1^*B$ is the closed two-form on
$Q$, $B$ is a closed two-form on $\mathcal{H}$ and the projection $\pi_1:
Q=\mathcal{H}\times V \rightarrow \mathcal{H}$ induces the map
$\pi_1^*: T^* \mathcal{H}\rightarrow T^*Q$.
We note that there is a magnetic term on the cotangent
bundle of the Heisenberg group $\mathcal{H}$, which is related to a
curvature two-form of a mechanical connection determined by the
reduction of center action of the Heisenberg group $\mathcal{H} $
(see Marsden et al. \cite{mamipera98}),
such that we can define a controlled magnetic Hamiltonian system with
symmetry of the Heisenberg group $\mathcal{H}$, and study the regular point
reduction of this system. Since a CMH system is also a
RCH system, but its symplectic structure
is given by a magnetic symplectic form. Thus, the set of
the CMH systems with symmetries is a subset of the set of
the RCH systems with symmetries, and the subset is not complete
under the regular point reduction of the RCH system,
because the magnetic may be vanishing for the reduced RCH system.
It is worthy of noting that
it is different from the regular point reduction of an RCH system
defined on a cotangent bundle with the canonical structure,
the regular point reduction of a CMH system reveals
the deeper relationship of the intrinsic geometrical structures
of the RCH systems on a cotangent bundle.
See Wang \cite{wa15a} for more details.

\subsection{Nonholonomic Controlled Magnetic Hamiltonian System}

A nonholonomic CMH system is the 6-tuple
$(T^\ast Q,\omega^B,\mathcal{D},H,F,W)$,
which is a CMH system with a
$\mathcal{D}$-completely and $\mathcal{D}$-regularly nonholonomic
constraint $\mathcal{D} \subset TQ$.
Under the restriction given by constraint, in general, the dynamical
vector field of a nonholonomic CMH system may not be
magnetic Hamiltonian, however the system is a dynamical system
closely related to a magnetic Hamiltonian system.
By analyzing carefully the structure for the nonholonomic dynamical
vector field, we give a geometric formulation of the distributional CMH system,
which is determined by a non-degenerate distributional two-form induced
from the magnetic symplectic form. Moreover,
for the nonholonomic CMH system $(T^*Q,\omega^B,\mathcal{D},H,F,u)$
 with a control law $u$ and with an associated distributional CMH system
$(\mathcal{K},\omega^B_{\mathcal {K}},H_{\mathcal {K}}, F^B_{\mathcal {K}},u^B_{\mathcal {K}})$,
we can derive precisely the geometric constraint conditions of
the non-degenerate distributional two-form $\omega^B_{\mathcal{K}}$
for the dynamical vector field $X^B_{(\mathcal{K},\omega^B_{\mathcal{K}},
H_{\mathcal{K}}, F^B_{\mathcal{K}}, u^B_{\mathcal{K}})}$;
that is, the two types of Hamilton-Jacobi equation for the distributional
CMH system $(\mathcal{K},\omega^B_{\mathcal {K}}, H_{\mathcal {K}},
F^B_{\mathcal {K}}, u^B_{\mathcal {K}})$;
see Wang \cite{wa 21b} for more details.\\

In the following we consider the nonholonomic CMH
system with symmetry.
By analyzing carefully the structure of the reducible dynamical vector field
of the nonholonomic CMH system,
we give a geometric formulation of
the nonholonomic reduced distributional CMH system,
Moreover, we derive precisely the geometric constraint conditions of
the non-degenerate, and nonholonomic reduced distributional two-form
for the nonholonomic reducible dynamical vector field,
that is, the two types of Hamilton-Jacobi equations for the
nonholonomic reduced distributional CMH system,
which are an extension of the two types of Hamilton-Jacobi equations
for the nonholonomic reduced distributional Hamiltonian system
given in Le\'{o}n and Wang \cite{lewa15}.\\

Assume that the Lie group $G$ acts smoothly on the manifold $Q$ by the left,
and that we also consider the natural lifted actions on $TQ$ and $T^* Q$,
and assume that the cotangent lifted left action $\Phi^{T^\ast}:
G\times T^\ast Q\rightarrow T^\ast Q$ is free, proper and
symplectic with respect to the magnetic symplectic form
$\omega^B$ on $T^* Q$.
Then the orbit space $T^* Q/ G$ is a smooth manifold and the
canonical projection $\pi_{/G}: T^* Q \rightarrow T^* Q /G $ is
a surjective submersion. For the cotangent lifted left action
$\Phi^{T^\ast}: G\times T^\ast Q\rightarrow T^\ast Q$,
assume that $H: T^*Q \rightarrow \mathbb{R}$ is a
$G$-invariant Hamiltonian, and that the fiber-preserving map
$F:T^\ast Q\rightarrow T^\ast Q$ and the control subset
$W$ of\; $T^\ast Q$ are both $G$-invariant,
and that the $\mathcal{D}$-completely and
$\mathcal{D}$-regularly nonholonomic constraint $\mathcal{D}\subset
TQ$ is a $G$-invariant distribution for the tangent lifted left action $\Phi^{T}:
G\times TQ\rightarrow TQ$; that is, the tangent of group action maps
$\mathcal{D}_q$ to $\mathcal{D}_{gq}$ for any
$q\in Q $. A nonholonomic CMH system with symmetry
is 7-tuple $(T^*Q,G,\omega^B,\mathcal{D},H, F,W)$, which is an
CMH system with symmetry and $G$-invariant
nonholonomic constraint $\mathcal{D}$.\\

In the following we first consider the nonholonomic reduction
of a nonholonomic CMH system with symmetry
$(T^*Q,G,\omega^B,\mathcal{D},H, F,W)$.
Note that the Legendre transformation $\mathcal{F}L: TQ
\rightarrow T^*Q$ is a fiber-preserving map,
and that $\mathcal{D}\subset TQ$ is $G$-invariant
for the tangent lifted left action $\Phi^{T}: G\times TQ\rightarrow TQ, $
then the constraint submanifold
$\mathcal{M}=\mathcal{F}L(\mathcal{D})\subset T^*Q$ is
$G$-invariant for the cotangent lifted left action $\Phi^{T^\ast}:
G\times T^\ast Q\rightarrow T^\ast Q$,
For the nonholonomic CMH system with symmetry
$(T^*Q,G, \omega^B,\mathcal{D},H,F,W)$,
in the same way, we define the distribution $\mathcal{F}$, which is the pre-image of the
nonholonomic constraints $\mathcal{D}$ for the map $T\pi_Q: TT^* Q
\rightarrow TQ$; that is, $\mathcal{F}=(T\pi_Q)^{-1}(\mathcal{D})$,
and the distribution $\mathcal{K}=\mathcal{F} \cap T\mathcal{M}$.
Moreover, we can also define the distributional two-form $\omega^B_\mathcal{K}$,
which is induced from the magnetic symplectic form $\omega^B$ on $T^* Q$; that is,
$\omega^B_\mathcal{K}= \tau_{\mathcal{K}}\cdot \omega^B_{\mathcal{M}},$ and
$\omega^B_{\mathcal{M}}= i_{\mathcal{M}}^* \omega^B $.
If the admissibility condition $\mathrm{dim}\mathcal{M}=
\mathrm{rank}\mathcal{F}$ and the compatibility condition
$T\mathcal{M}\cap \mathcal{F}^\bot= \{0\}$ hold, then
$\omega^B_\mathcal{K}$ is non-degenerate as a
bilinear form on each fibre of $\mathcal{K}$, there exists a vector
field $X^B_\mathcal{K}$ on $\mathcal{M}$ which takes values in the
constraint distribution $\mathcal{K}$, such that for the function $H_\mathcal{K}$,
the following distributional magnetic Hamiltonian equation holds, that is,
\begin{align}
\mathbf{i}_{X^B_\mathcal{K}}\omega^B_\mathcal{K}
=\mathbf{d}H_\mathcal{K},
\label{7.7} \end{align}
where the function $H_{\mathcal{K}}$ satisfies
$\mathbf{d}H_{\mathcal{K}}= \tau_{\mathcal{K}}\cdot \mathbf{d}H_{\mathcal {M}}$,
and $H_\mathcal{M}= \tau_{\mathcal{M}}\cdot H$
is the restriction of $H$ to $\mathcal{M}$, and
from the equation (7.7), we have that
$X^B_{\mathcal{K}}=\tau_{\mathcal{K}}\cdot X^B_H $.\\

In the following we define that the quotient space
$\bar{\mathcal{M}}=\mathcal{M}/G$ of the $G$-orbit in $\mathcal{M}$
is a smooth manifold with projection $\pi_{/G}:
\mathcal{M}\rightarrow \bar{\mathcal{M}}( \subset T^* Q /G),$ which
is a surjective submersion. The reduced magnetic symplectic form
$\omega^B_{\bar{\mathcal{M}}}= \pi^*_{/G} \cdot \omega^B_{\mathcal{M}}$
on $\bar{\mathcal{M}}$ is induced from the magnetic symplectic form $\omega^B_{\mathcal{M}}
= i_{\mathcal{M}}^* \omega^B $ on $\mathcal{M}$.
Since $G$ is the symmetry group of the system
$(T^*Q,G,\omega^B,\mathcal{D},H, F,W)$, all intrinsically
defined vector fields and distributions are pushed down to
$\bar{\mathcal{M}}$. In particular, the vector field $X^B_\mathcal{M}$
on $\mathcal{M}$ is pushed down to a vector field
$X^B_{\bar{\mathcal{M}}}=T\pi_{/G}\cdot X^B_\mathcal{M}$, and the
distribution $\mathcal{K}$ is pushed down to a distribution
$T\pi_{/G}\cdot \mathcal{K}$ on $\bar{\mathcal{M}}$, and the
Hamiltonian $H$ is pushed down to $h_{\bar{\mathcal{M}}}$, such that
$h_{\bar{\mathcal{M}}}\cdot \pi_{/G}=
\tau_{\mathcal{M}}\cdot H$. However, $\omega^B_\mathcal{K}$ need not
to be pushed down to a distributional two-form defined on $T\pi_{/G}\cdot
\mathcal{K}$, despite of the fact that $\omega^B_\mathcal{K}$ is
$G$-invariant. This is because there may be infinitesimal symmetry
$\eta_{\mathcal{K}}$ that lies in $\mathcal{M}$, such that
$\mathbf{i}_{\eta_\mathcal{K}} \omega^B_\mathcal{K}\neq 0$. From Bates
and $\acute{S}$niatycki \cite{basn93}, we know that in order to eliminate
this difficulty, $\omega^B_\mathcal{K}$ is restricted to a
sub-distribution $\mathcal{U}$ of $\mathcal{K}$ defined by
$$\mathcal{U}=\{u\in\mathcal{K} \; | \; \omega^B_\mathcal{K}(u,v)
=0,\quad \forall \; v \in \mathcal{V}\cap \mathcal{K}\},$$ where
$\mathcal{V}$ is the distribution on $\mathcal{M}$ tangent to the
orbits of $G$ in $\mathcal{M}$ and it is spanned by the infinitesimal
symmetries. Clearly, $\mathcal{U}$ and $\mathcal{V}$ are both
$G$-invariant, project down to $\bar{\mathcal{M}}$ and
$T\pi_{/G}\cdot \mathcal{V}=0$, and define the distribution $\bar{\mathcal{K}}$ by
$\bar{\mathcal{K}}= T\pi_{/G}\cdot \mathcal{U}$. Moreover, we take
that $\omega^B_\mathcal{U}= \tau_{\mathcal{U}}\cdot
\omega^B_{\mathcal{M}}$ is the restriction of the induced magnetic symplectic form
$\omega^B_{\mathcal{M}}$ on $T^*\mathcal{M}$ fibrewise to the
distribution $\mathcal{U}$, where $\tau_{\mathcal{U}}$ is the
restriction map to distribution $\mathcal{U}$, and the
$\omega^B_{\mathcal{U}}$ is pushed down to a
distributional two-form $\omega^B_{\bar{\mathcal{K}}}$ on
$\bar{\mathcal{K}}$, such that $\pi_{/G}^*
\omega^B_{\bar{\mathcal{K}}}= \omega^B_{\mathcal{U}}$.
We know that distributional two-form
$\omega^B_{\bar{\mathcal{K}}}$ is not a "true two-form"
on a manifold, which is called the nonholonomic reduced
distributional two-form to avoid any confusion.\\

From the above construction we know that,
if the admissibility condition $\mathrm{dim}\bar{\mathcal{M}}=
\mathrm{rank}\bar{\mathcal{F}}$ and the compatibility condition
$T\bar{\mathcal{M}} \cap \bar{\mathcal{F}}^\bot= \{0\}$ hold, where
$\bar{\mathcal{F}}^\bot$ denotes the symplectic orthogonal of
$\bar{\mathcal{F}}$ with respect to the reduced magnetic symplectic form
$\omega^B_{\bar{\mathcal{M}}}$, then the nonholonomic reduced
distributional two-form
$\omega^B_{\bar{\mathcal{K}}}$ is non-degenerate as a bilinear form on
each fibre of $\bar{\mathcal{K}}$, and hence there exists a vector field
$X^B_{\bar{\mathcal{K}}}$ on $\bar{\mathcal{M}}$ which takes values in
the constraint distribution $\bar{\mathcal{K}}$, such that the
reduced distributional magnetic Hamiltonian equation holds, that is,
\begin{align}
\mathbf{i}_{X^B_{\bar{\mathcal{K}}}}\omega^B_{\bar{\mathcal{K}}}
=\mathbf{d}h_{\bar{\mathcal{K}}},
\label{7.8} \end{align}
where $\mathbf{d}h_{\bar{\mathcal{K}}}$ is the restriction of
$\mathbf{d}h_{\bar{\mathcal{M}}}$ to $\bar{\mathcal{K}}$ and
the function $h_{\bar{\mathcal{K}}}:\bar{M}(\subset T^* Q/G)\rightarrow \mathbb{R}$ satisfies
$\mathbf{d}h_{\bar{\mathcal{K}}}= \tau_{\bar{\mathcal{K}}}\cdot \mathbf{d}h_{\bar{\mathcal{M}}}$,
and $h_{\bar{\mathcal{M}}}\cdot \pi_{/G}= H_{\mathcal{M}}$ and
$H_{\mathcal{M}}$ is the restriction of the Hamiltonian function $H$
to $\mathcal{M}$, and the function
$h_{\bar{\mathcal{M}}}:\bar{M}(\subset T^* Q/G)\rightarrow \mathbb{R}$.
In addition, from the distributional magnetic Hamiltonian equation (7.7),
$\mathbf{i}_{X^B_\mathcal{K}}\omega^B_\mathcal{K}=\mathbf{d}H_\mathcal
{K},$ we have that $X^B_{\mathcal{K}}=\tau_{\mathcal{K}}\cdot X^B_H, $
and from the reduced distributional magnetic Hamiltonian equation (7.8),
$\mathbf{i}_{X^B_{\bar{\mathcal{K}}}}\omega^B_{\bar{\mathcal{K}}}
=\mathbf{d}h_{\bar{\mathcal{K}}}$, we have that
$X^B_{\bar{\mathcal{K}}}
=\tau_{\bar{\mathcal{K}}}\cdot X^B_{h_{\bar{\mathcal{K}}}},$
where $ X^B_{h_{\bar{\mathcal{K}}}}$ is the magnetic Hamiltonian vector field of
the function $h_{\bar{\mathcal{K}}}$ with respect to the reduced magnetic symplectic
form $\omega^B_{\bar{\mathcal{M}}}$,
and the vector fields $X^B_{\mathcal{K}}$
and $X^B_{\bar{\mathcal{K}}}$ are $\pi_{/G}$-related,
that is, $X^B_{\bar{\mathcal{K}}}\cdot \pi_{/G}=T\pi_{/G}\cdot X^B_{\mathcal{K}}.$ \\

Moreover, considering the external force $F$ and control subset $W$,
and we define the vector fields $F^B_\mathcal{K}
=\tau_{\mathcal{K}}\cdot \textnormal{vlift}(F_{\mathcal{M}})X^B_H,$
and for a control law $u\in W$,
$u^B_\mathcal{K}= \tau_{\mathcal{K}}\cdot  \textnormal{vlift}(u_{\mathcal{M}})X^B_H,$
where $F_\mathcal{M}= \tau_{\mathcal{M}}\cdot F$ and
$u_\mathcal{M}= \tau_{\mathcal{M}}\cdot u$ are the restrictions of
$F$ and $u$ to $\mathcal{M}$, that is, $F^B_\mathcal{K}$ and $u^B_\mathcal{K}$
are the restrictions of the changes of magnetic Hamiltonian vector field $X^B_H$
under the actions of $F_\mathcal{M}$ and $u_\mathcal{M}$ to $\mathcal{K}$,
then the 5-tuple $(\mathcal{K},\omega^B_{\mathcal{K}},
H_\mathcal{K}, F^B_\mathcal{K}, u^B_\mathcal{K})$
is a distributional CMH system corresponding to the nonholonomic CMH system with symmetry
$(T^*Q,G,\omega^B,\mathcal{D},H,F,u)$,
and the dynamical vector field of the distributional CMH system
can be expressed by
\begin{align}
\tilde{X}=X^B_{(\mathcal{K},\omega^B_{\mathcal{K}},
H_{\mathcal{K}}, F^B_{\mathcal{K}}, u^B_{\mathcal{K}})}
=X^B_\mathcal {K}+ F^B_{\mathcal{K}}+u^B_{\mathcal{K}},
\label{7.9} \end{align}
which is the synthesis
of the nonholonomic dynamical vector field $X^B_{\mathcal{K}}$ and
the vector fields $F^B_{\mathcal{K}}$ and $u^B_{\mathcal{K}}$.
Assume that the vector fields $F^B_\mathcal{K}$ and $u^B_\mathcal{K}$
on $\mathcal{M}$ are pushed down to the vector fields
$f^B_{\bar{\mathcal{M}}}= T\pi_{/G}\cdot F^B_\mathcal{K}$ and
$u^B_{\bar{\mathcal{M}}}=T\pi_{/G}\cdot u^B_\mathcal{K}$ on $\bar{\mathcal{M}}$.
Then we define that $f^B_{\bar{\mathcal{K}}}=T\tau_{\bar{\mathcal{K}}}\cdot f^B_{\bar{\mathcal{M}}}$ and
$u^B_{\bar{\mathcal{K}}}=T\tau_{\bar{\mathcal{K}}}\cdot u^B_{\bar{\mathcal{M}}};$
that is, $f^B_{\bar{\mathcal{K}}}$ and
$u^B_{\bar{\mathcal{K}}}$ are the restrictions of
$f^B_{\bar{\mathcal{M}}}$ and $u^B_{\bar{\mathcal{M}}}$ to $\bar{\mathcal{K}}$,
where $\tau_{\bar{\mathcal{K}}}$
is the restriction map to distribution $\bar{\mathcal{K}}$,
and $T\tau_{\bar{\mathcal{K}}}$ is the tangent map of $\tau_{\bar{\mathcal{K}}}$.
Then the 5-tuple $(\bar{\mathcal{K}},\omega^B_{\bar{\mathcal{K}}},
h_{\bar{\mathcal{K}}}, f^B_{\bar{\mathcal{K}}}, u^B_{\bar{\mathcal{K}}})$
is a nonholonomic reduced distributional CMH system of the nonholonomic
reducible CMH system with symmetry $(T^*Q,G,\omega^B,\mathcal{D},H,F,W)$,
as well as with a control law $u \in W$.
Thus, the geometrical formulation of a nonholonomic reduced distributional
CMH system may be summarized as follows.

\begin{defi} (Nonholonomic Reduced Distributional CMH System)
Assume that the 7-tuple \\ $(T^*Q,G,\omega^B,\mathcal{D},H,F,W)$ is a nonholonomic
reducible CMH system with symmetry, where $\omega^B$ is the magnetic
symplectic form on $T^* Q$, and $\mathcal{D}\subset TQ$ is a
$\mathcal{D}$-completely and $\mathcal{D}$-regularly nonholonomic
constraint of the system, and $\mathcal{D}$, $H, F$ and $W$ are all
$G$-invariant. If there exists a nonholonomic reduced distribution $\bar{\mathcal{K}}$,
an associated non-degenerate  and nonholonomic reduced
distributional two-form $\omega^B_{\bar{\mathcal{K}}}$
and a vector field $X^B_{\bar{\mathcal {K}}}$ on the reduced constraint
submanifold $\bar{\mathcal{M}}=\mathcal{M}/G, $ where
$\mathcal{M}=\mathcal{F}L(\mathcal{D})\subset T^*Q$, such that the
nonholonomic reduced distributional magnetic Hamiltonian equation holds: that is,
$ \mathbf{i}_{X^B_{\bar{\mathcal{K}}}}\omega^B_{\bar{\mathcal{K}}} =
\mathbf{d}h_{\bar{\mathcal{K}}}, $
where $\mathbf{d}h_{\bar{\mathcal{K}}}$ is the restriction of
$\mathbf{d}h_{\bar{\mathcal{M}}}$ to $\bar{\mathcal{K}}$ and
the function $h_{\bar{\mathcal{K}}}$ satisfies
$\mathbf{d}h_{\bar{\mathcal{K}}}= \tau_{\bar{\mathcal{K}}}\cdot \mathbf{d}h_{\bar{\mathcal{M}}}$
and $h_{\bar{\mathcal{M}}}\cdot \pi_{/G}= H_{\mathcal{M}}$,
and the vector fields $f^B_{\bar{\mathcal{K}}}=T\tau_{\bar{\mathcal{K}}}\cdot f^B_{\bar{\mathcal{M}}}$ and
$u^B_{\bar{\mathcal{K}}}=T\tau_{\bar{\mathcal{K}}}\cdot u^B_{\bar{\mathcal{M}}}$ as defined above.
Then the 5-tuple $(\bar{\mathcal{K}},\omega^B_{\bar{\mathcal {K}}},h_{\bar{\mathcal{K}}},
f^B_{\bar{\mathcal{K}}}, u^B_{\bar{\mathcal{K}}})$
is called a nonholonomic reduced distributional CMH system
of the nonholonomic reducible CMH system $(T^*Q,G,\omega^B,\mathcal{D},H,F,W)$
with a control law $u \in W$, and $X^B_{\bar{\mathcal {K}}}$ is
called a nonholonomic reduced dynamical vector field.
Denote that
\begin{align}
\hat{X}=X^B_{(\bar{\mathcal{K}},\omega^B_{\bar{\mathcal{K}}},
h_{\bar{\mathcal{K}}}, f^B_{\bar{\mathcal{K}}}, u^B_{\bar{\mathcal{K}}})}
=X^B_{\bar{\mathcal{K}}}+ f^B_{\bar{\mathcal{K}}}+u^B_{\bar{\mathcal{K}}}
\label{7.10} \end{align}
is the dynamical vector field of the
nonholonomic reduced distributional CMH system
$(\bar{\mathcal{K}},\omega^B_{\bar{\mathcal{K}}},h_{\bar{\mathcal{K}}},
f^B_{\bar{\mathcal{K}}}, \\ u^B_{\bar{\mathcal{K}}})$, which is the synthesis
of the nonholonomic reduced dynamical vector field $X^B_{\bar{\mathcal{K}}}$ and
the vector fields $F^B_{\bar{\mathcal{K}}}$ and $u^B_{\bar{\mathcal{K}}}$.
Under the above
circumstances, we refer to $(T^*Q,G,\omega^B,\mathcal{D},H,F,u)$ as a
nonholonomic reducible CMH system with the associated
distributional CMH system
$(\mathcal{K},\omega^B_{\mathcal {K}},H_{\mathcal{K}}, F^B_{\mathcal{K}}, u^B_{\mathcal{K}})$
and the nonholonomic reduced distributional CMH system
$(\bar{\mathcal{K}},\omega^B_{\bar{\mathcal{K}}},h_{\bar{\mathcal{K}}},
f^B_{\bar{\mathcal{K}}}, u^B_{\bar{\mathcal{K}}})$.
The dynamical vector fields
$\tilde{X}=X^B_{(\mathcal{K},\omega^B_{\mathcal{K}},
H_{\mathcal{K}}, F^B_{\mathcal{K}}, u^B_{\mathcal{K}})}$
and $\hat{X}= X^B_{(\bar{\mathcal{K}},\omega^B_{\bar{\mathcal{K}}},
h_{\bar{\mathcal{K}}}, f^B_{\bar{\mathcal{K}}}, u^B_{\bar{\mathcal{K}}})}$
are $\pi_{/G}$-related; that is,
$\hat{X}\cdot \pi_{/G}=T\pi_{/G}\cdot \tilde{X}.$
\end{defi}

For a given nonholonomic reducible CMH system
$(T^*Q,G,\omega^B,\mathcal{D},H,F,u)$ with the associated
distributional CMH system
$(\mathcal{K},\omega^B_{\mathcal {K}},H_{\mathcal{K}}, F^B_{\mathcal{K}}, u^B_{\mathcal{K}})$
and the nonholonomic reduced distributional CMH system
$(\bar{\mathcal{K}},\omega^B_{\bar{\mathcal{K}}},h_{\bar{\mathcal{K}}}, f^B_{\bar{\mathcal{K}}}, u^B_{\bar{\mathcal{K}}})$, the magnetic vector field $X^0= X^B_H-X_H, $
which is determined by the magnetic equation
$\mathbf{i}_{X^0}\omega=\mathbf{i}_{X^B_H}( \pi_Q^*B)$
on $T^*Q$. Note that the vector fields $X^0$, $ X^B_H$, $X_H $
and the distribution
$\mathcal{K}$ are pushed down to
$\bar{\mathcal{M}}$; that is,
$X^0_{\bar{\mathcal{M}}}=T\pi_{/G}\cdot X^0_\mathcal{M}$,
$X^B_{\bar{\mathcal{M}}}=T\pi_{/G}\cdot X^B_\mathcal{M}$,
and $X_{\bar{\mathcal{M}}}=T\pi_{/G}\cdot X_\mathcal{M}$,
where the vector field $X^0_\mathcal{M}$, $X^B_\mathcal{M}$
and $X_\mathcal{M}$ are the restrictions of $X^0$, $ X^B_H$ and $X_H$
on $\mathcal{M}$ and the
distribution $\mathcal{K}$ is pushed down to a distribution
$T\pi_{/G}\cdot \mathcal{K}$ on $\bar{\mathcal{M}}$,
and define the distribution $\bar{\mathcal{K}}$ by
$\bar{\mathcal{K}}= T\pi_{/G}\cdot \mathcal{U}$.
Denote that $X^0_{\bar{\mathcal {K}}}=\tau_{\bar{\mathcal {K}}}(X^0_{\bar{\mathcal{M}}})
= \tau_{\bar{\mathcal {K}}}(X^B_{\bar{\mathcal{M}}})- \tau_{\bar{\mathcal {K}}}(X_{\bar{\mathcal{M}}})
=X^B_{\bar{\mathcal {K}}}- X_{\bar{\mathcal {K}}},$
from the expression (7.10) of the dynamical
vector field of the nonholonomic reduced distributional CMH system
$(\bar{\mathcal{K}},\omega^B_{\bar{\mathcal{K}}},h_{\bar{\mathcal{K}}}, f^B_{\bar{\mathcal{K}}}, u^B_{\bar{\mathcal{K}}})$, we have that
\begin{align}\hat{X}
=X^B_{\bar{\mathcal{K}}}+ F^B_{\bar{\mathcal{K}}}+u^B_{\bar{\mathcal{K}}}
=X_{\bar{\mathcal{K}}}+ X^0_{\bar{\mathcal{K}}}+ F^B_{\bar{\mathcal{K}}}+u^B_{\bar{\mathcal{K}}}.
\label{7.11} \end{align}
If the vector fields $F^B_{\bar{\mathcal{K}}}$ and $u^B_{\bar{\mathcal{K}}}$ satisfy the following condition
\begin{equation}
 X^0_{\bar{\mathcal{K}}}+ F^B_{\bar{\mathcal{K}}}+u^B_{\bar{\mathcal{K}}}=0, \;\; \label{7.12}
\end{equation}
then from (7.11) we have that $X^B_{(\bar{\mathcal{K}},\omega^B_{\bar{\mathcal{K}}},
h_{\bar{\mathcal{K}}}, f^B_{\bar{\mathcal{K}}}, u^B_{\bar{\mathcal{K}}})}
=X_{\bar{\mathcal{K}}}, $ that is, in this case the dynamical vector
field of the nonholonomic reduced distributional CMH system is just the dynamical
vector field of the nonholonomic reduced canonical distributional Hamiltonian system
without the actions of the magnetic, the external force and the control.
Thus, the condition (7.12) is called the magnetic vanishing condition for
the nonholonomic reduced distributional CMH system $(\bar{\mathcal{K}},\omega^B_{\bar{\mathcal{K}}},h_{\bar{\mathcal{K}}}, f^B_{\bar{\mathcal{K}}}, u^B_{\bar{\mathcal{K}}})$.\\

Since the non-degenerate and nonholonomic reduced distributional two-form
$\omega^B_{\bar{\mathcal{K}}}$ is not a "true two-form"
on a manifold, and it is not symplectic, and hence
the nonholonomic reduced distributional CMH system
$(\bar{\mathcal{K}},\omega^B_{\bar{\mathcal{K}}},h_{\bar{\mathcal{K}}}, f^B_{\bar{\mathcal{K}}}, u^B_{\bar{\mathcal{K}}})$ is not a Hamiltonian system,
and it has no yet generating function,
and hence we can not describe the Hamilton-Jacobi equation for the nonholonomic reduced
distributional CMH system the same as in Theorem 2.1.
However, for a given nonholonomic reducible CMH system
$(T^*Q,G,\omega^B,\mathcal{D},H,F,u)$ with the associated
distributional CMH system
$(\mathcal{K},\omega^B_{\mathcal {K}},H_{\mathcal{K}}, F^B_{\mathcal{K}}, u^B_{\mathcal{K}})$
and the nonholonomic reduced distributional CMH system
$(\bar{\mathcal{K}},\omega^B_{\bar{\mathcal {K}}},h_{\bar{\mathcal{K}}}, f^B_{\bar{\mathcal{K}}}, u^B_{\bar{\mathcal{K}}})$, by using Lemma 7.2 and the following Lemma 7.13,
we can derive precisely
the geometric constraint conditions of the nonholonomic reduced distributional two-form
$\omega^B_{\bar{\mathcal{K}}}$ for the nonholonomic reducible dynamical vector field
$\tilde{X}=X^B_{(\mathcal{K},\omega^B_{\mathcal{K}},
H_{\mathcal{K}}, F^B_{\mathcal{K}}, u^B_{\mathcal{K}})}$;
that is, the two types of Hamilton-Jacobi equations for the
nonholonomic reduced distributional CMH system
$(\bar{\mathcal{K}},\omega^B_{\bar{\mathcal {K}}},h_{\bar{\mathcal{K}}}, f^B_{\bar{\mathcal{K}}}, u^B_{\bar{\mathcal{K}}})$.\\

The following Lemma 7.13
can be regarded as an extension of the Lemma 6.3 ( given in
Le\'{o}n and Wang \cite{lewa15}) to it with the nonholonomic magnetic context,
and the lemma is a very important tool for the proof of
the Hamilton-Jacobi theorem for the
nonholonomic reduced distributional CMH system.
\begin{lemm}
Assume that $\gamma: Q \rightarrow T^*Q$ is a one-form on $Q$, and that
$\lambda=\gamma \cdot \pi_{Q}: T^* Q \rightarrow T^* Q ,$ and that
$\omega$ is the canonical symplectic form on $T^*Q$, and that
$\omega^B= \omega- \pi_Q^*B $
is the magnetic symplectic form on $T^*Q$.
If the Lagrangian $L$ is $\mathcal{D}$-regular, and assume that
$\textmd{Im}(\gamma)\subset \mathcal{M},$
where $\mathcal{M}=\mathcal{F}L(\mathcal{D}), $
then we have that $ X^B_{H}\cdot \gamma \in \mathcal{F}$ along
$\gamma$, and that $ X^B_{H}\cdot \lambda \in \mathcal{F}$ along
$\lambda$; that is, $T\pi_{Q}(X^B_H\cdot\gamma(q))\in
\mathcal{D}_{q}, \; \forall q \in Q $ and $T\pi_{Q}(X^B_H\cdot\lambda(q,p))\in
\mathcal{D}_{q}, \; \forall q \in Q, \; (q,p) \in T^* Q. $
Moreover, if a symplectic map $\varepsilon: T^* Q \rightarrow T^* Q $
with respect to the magnetic symplectic form $\omega^B$ satisfies the
condition that $\varepsilon(\mathcal{M})\subset \mathcal{M},$ then
we have that $ X^B_{H}\cdot \varepsilon \in \mathcal{F}$ along
$\varepsilon. $
\end{lemm}
See the proof and the more details in Wang \cite{wa21b}.\\

For convenience, the maps involved in the following
theorem are shown in Diagram-10.
\begin{center}
\hskip 0cm \xymatrix{ & \mathcal{M} \ar[d]_{X^B_{\mathcal{K}}}
\ar[r]^{i_{\mathcal{M}}} & T^* Q \ar[d]_{X^B_{H}}
 \ar[dr]^{\tilde{X}^\varepsilon} \ar[r]^{\pi_Q}
  & Q \ar[d]^{\tilde{X}^\gamma} \ar[r]^{\gamma}
  & T^* Q \ar[d]_{\tilde{X}} \ar[r]^{\pi_{/G}}
  & T^* Q/G \ar[d]_{X^B_{h_{\bar{\mathcal{M}}}}}
  & \mathcal{\bar{M}} \ar[l]_{i_{\mathcal{\bar{M}}}} \ar[d]_{X^B_{\mathcal{\bar{K}}}}\\
  & \mathcal{K}
  & T(T^*Q) \ar[l]^{\tau_{\mathcal{K}}}
  & TQ \ar[l]^{T\gamma}
  & T(T^* Q) \ar[l]^{T\pi_Q} \ar[r]_{T\pi_{/G}}
  & T(T^* Q/G) \ar[r]_{\tau_{\mathcal{\bar{K}}}} & \mathcal{\bar{K}} }
\end{center}
$$\mbox{Diagram-10}$$
\begin{theo} (Hamilton-Jacobi Theorem for a Nonholonomic
Reduced Distributional CMH System)
For a given nonholonomic reducible CMH system
$(T^*Q,G,\omega^B,\mathcal{D},H,F,u)$ with the associated
distributional CMH system
$(\mathcal{K},\omega^B_{\mathcal {K}},H_{\mathcal{K}}, F^B_{\mathcal{K}}, u^B_{\mathcal{K}})$
and the nonholonomic reduced distributional CMH system
$(\bar{\mathcal{K}},\omega^B_{\bar{\mathcal{K}}},h_{\bar{\mathcal{K}}}, f^B_{\bar{\mathcal{K}}}, u^B_{\bar{\mathcal{K}}})$, assume that
$\gamma: Q \rightarrow T^*Q$ is a one-form on $Q$,
and that $\lambda= \gamma \cdot \pi_{Q}: T^* Q \rightarrow T^* Q, $
and that for any $G$-invariant
symplectic map $\varepsilon: T^* Q \rightarrow T^* Q $
with respect to $\omega^B$, denote that
$\tilde{X}^\gamma = T\pi_{Q}\cdot \tilde{X}\cdot \gamma$ and
$\tilde{X}^\varepsilon = T\pi_{Q}\cdot \tilde{X}\cdot \varepsilon$,
where $\tilde{X}=X^B_{(\mathcal{K},\omega^B_{\mathcal{K}},
H_{\mathcal{K}}, F^B_{\mathcal{K}}, u^B_{\mathcal{K}})}
=X^B_\mathcal {K}+ F^B_{\mathcal{K}}+u^B_{\mathcal{K}}$
is the dynamical vector field of the distributional CMH system
corresponding to the nonholonomic reducible CMH system with symmetry
$(T^*Q,G,\omega^B,\mathcal{D},H,F,u)$. Moreover,
assume that $\textmd{Im}(\gamma)\subset \mathcal{M}, $ and that $\gamma$ is
$G$-invariant, and that $\varepsilon(\mathcal{M})\subset \mathcal{M}$,
and that $ \textmd{Im}(T\gamma)\subset \mathcal{K}. $ Denote that
$\bar{\gamma}=\pi_{/G}(\gamma): Q \rightarrow T^* Q/G ,$
and that
$\bar{\lambda}=\pi_{/G}(\lambda): T^* Q \rightarrow T^* Q/G, $ and that
$\bar{\varepsilon}=\pi_{/G}(\varepsilon): T^* Q \rightarrow T^* Q/G. $
Then the following two assertions hold:\\
\noindent $(\mathbf{i})$
If the one-form $\gamma: Q \rightarrow T^*Q $ satisfies the condition,
$\mathbf{d}\gamma=-B $ on $\mathcal{D}$ with respect to
$T\pi_Q: TT^* Q \rightarrow TQ, $ then $\bar{\gamma}$ is a solution
of the Type I Hamilton-Jacobi equation $T\bar{\gamma}\cdot \tilde{X}^ \gamma =
X^B_{\bar{\mathcal{K}}}\cdot \bar{\gamma} $
for the nonholonomic reduced distributional CMH system
$(\bar{\mathcal{K}},\omega^B_{\bar{\mathcal{K}}},h_{\bar{\mathcal{K}}}, f^B_{\bar{\mathcal{K}}}, u^B_{\bar{\mathcal{K}}})$, where
$X^B_{\bar{\mathcal{K}}}$ is the nonholonomic reduced dynamical vector field. \\
\noindent $(\mathbf{ii})$
The $\varepsilon$ and $\bar{\varepsilon}$ satisfy the Type II Hamilton-Jacobi equation
$T\bar{\gamma}\cdot \tilde{X}^\varepsilon
= X^B_{\bar{\mathcal{K}}}\cdot \bar{\varepsilon}$
if and only if they satisfy the equation
$\tau_{\bar{\mathcal{K}}}\cdot T\bar{\varepsilon}\cdot X^B_{h_{\bar{\mathcal{K}}}\cdot
\bar{\varepsilon}}= T\bar{\lambda} \cdot \tilde{X}\cdot \varepsilon, $
where $ X^B_{h_{\bar{\mathcal{K}}} \cdot\bar{\varepsilon}}$ is the magnetic Hamiltonian
vector field of the function $h_{\bar{\mathcal{K}}}\cdot \bar{\varepsilon}: T^* Q\rightarrow
\mathbb{R}. $
\end{theo}
See the proof and the more details in Wang \cite{wa21b}.\\

For a given nonholonomic reducible CMH system
$(T^*Q,G,\omega^B,\mathcal{D},H,F,u)$ with the associated
the distributional CMH system $(\mathcal{K},\omega^B_{\mathcal{K}},
H_{\mathcal{K}}, F^B_{\mathcal{K}}, u^B_{\mathcal{K}})$
and the nonholonomic reduced distributional CMH system
$(\bar{\mathcal{K}},\omega^B_{\bar{\mathcal{K}}},h_{\bar{\mathcal{K}}},
f^B_{\bar{\mathcal{K}}}, u^B_{\bar{\mathcal{K}}})$,
we know that the nonholonomic dynamical vector field
$X^B_{\mathcal{K}}$ and the nonholonomic reduced dynamical vector field
$X^B_{\bar{\mathcal{K}}}$ are $\pi_{/G}$-related; that is,
$X^B_{\bar{\mathcal{K}}}\cdot \pi_{/G}=T\pi_{/G}\cdot X^B_{\mathcal{K}}.$
Then we can prove the following Theorem 7.15,
which states the relationship between the solutions of the Type II
Hamilton-Jacobi equations and nonholonomic reduction.

\begin{theo}
For a given nonholonomic reducible CMH system
$(T^*Q,G,\omega^B,\mathcal{D},H,F,u)$ with the associated
the distributional CMH system $(\mathcal{K},\omega^B_{\mathcal{K}},
H_{\mathcal{K}}, F^B_{\mathcal{K}}, u^B_{\mathcal{K}})$
and the nonholonomic reduced distributional CMH system
$(\bar{\mathcal{K}},\omega^B_{\bar{\mathcal{K}}},h_{\bar{\mathcal{K}}},
f^B_{\bar{\mathcal{K}}}, u^B_{\bar{\mathcal{K}}})$, assume that
$\gamma: Q \rightarrow T^*Q$ is a one-form on $Q$,
and that $\varepsilon: T^* Q \rightarrow T^* Q $ is a $G$-invariant
symplectic map with respect to $\omega^B$. Denote that
$\bar{\gamma}=\pi_{/G}(\gamma): Q \rightarrow T^* Q/G $, and that
$\bar{\varepsilon}=\pi_{/G}(\varepsilon): T^* Q \rightarrow T^* Q/G. $
Under the hypotheses and notations of Theorem 7.14, then we have that
$\varepsilon$ is a solution of the Type II Hamilton-Jacobi equation $T\gamma\cdot
\tilde{X}^\varepsilon= X^B_{\mathcal{K}}\cdot \varepsilon $ for the distributional
CMH system $(\mathcal{K},\omega^B_{\mathcal{K}},H_{\mathcal{K}},
F^B_{\mathcal{K}}, u^B_{\mathcal{K}})$ if and only if
$\varepsilon$ and $\bar{\varepsilon}$ satisfy the Type II
Hamilton-Jacobi equation $T\bar{\gamma}\cdot \tilde{X}^\varepsilon =
X^B_{\bar{\mathcal{K}}}\cdot \bar{\varepsilon} $ for the nonholonomic reduced
distributional CMH system $ (\bar{\mathcal{K}},
\omega^B_{\bar{\mathcal{K}}}, h_{\bar{\mathcal{K}}},
f^B_{\bar{\mathcal{K}}}, u^B_{\bar{\mathcal{K}}} ). $
\end{theo}
See the proof and the more details in Wang \cite{wa21b}.

\begin{rema}
It is worthy of noting that
the Type I Hamilton-Jacobi equation
$T\bar{\gamma}\cdot \tilde{X}^ \gamma = X^B_{\bar{\mathcal{K}}}\cdot
\bar{\gamma} $ is the equation of the nonholonomic
reduced differential one-form $\bar{\gamma}$, and that
the Type II Hamilton-Jacobi equation $T\bar{\gamma}\cdot \tilde{X}^\varepsilon
= X^B_{\bar{\mathcal{K}}}\cdot \bar{\varepsilon}$ is the equation of the symplectic
diffeomorphism map $\varepsilon$ and the nonholonomic reduced symplectic
diffeomorphism map $\bar{\varepsilon}$.
If the nonholonomic CMH system with symmetry we considered
has not any the external force and the control; that is, $F=0 $ and $W=\emptyset$,
in this case, the nonholonomic CMH system with symmetry
$(T^*Q,G,\omega^B,\mathcal{D},H,F,W)$ is just the nonholonomic magnetic Hamiltonian system
with symmetry $(T^*Q,G,\omega^B,\mathcal{D},H)$,
and with the magnetic symplectic form $\omega^B$ on $T^*Q$.
From the Hamilton-Jacobi theorem; that is,
Theorem 7.14, we can get the Theorems 5.2 and 5.3 given in Wang \cite{wa21c}.
It shows that Theorem 7.14 can be regarded as an extension of the two types of
Hamilton-Jacobi theorem for the nonholonomic magnetic Hamiltonian system with symmetry
to that for the system with the external force and the control.
In particular, in this case, if $B=0$, then
the magnetic symplectic form $\omega^B$
is just the canonical symplectic form $\omega$ on $T^*Q$,
from Theorem 7.14, we can also get the Theorems 4.2
and 4.3 given in Le\'{o}n and Wang \cite{lewa15}.
It shows that the Theorem 7.14 can be regarded as an extension of the two types of
Hamilton-Jacobi theorem for the nonholonomic Hamiltonian system with symmetry
to that for the system with the magnetic,
the external force and the control.
\end{rema}

In addition, we also give some generalizations of the above
results from the viewpoint of change of the geometrical structures.
A natural problem is what and how we could do, if we define a controlled
Hamiltonian system on the cotangent bundle $T^*Q$ by using a Poisson
structure, and if the symplectic reduction procedure given by Marsden \textit{et
al.} \cite{mawazh10} does not work or is not efficient enough. In
Wang and Zhang \cite{wazh12}, we study the optimal reduction theory
of a CH system with Poisson structure and symmetry, by using the
optimal momentum map and the reduced Poisson tensor (resp. the reduced
symplectic form). We prove the optimal point reduction,
the optimal orbit reduction, and the
regular Poisson reduction theorems for the CH system, and explain the
relationships between the OpCH-equivalence, the OoCH-equivalence,
the RPR-CH-equivalence for the optimal reducible CH systems with symmetries
and the CH-equivalence for the associated optimal reduced CH systems;
see Wang and Zhang \cite{wazh12} for the more details.\\

It is worthy of noting that if there is no momentum map of the Lie group action
for our considered system, then the reduction procedures given in
Marsden et al. \cite{mawazh10} and Wang and Zhang \cite{wazh12} can
not work. One must look for a new way. On the other hand, motivated
by the work of Poisson reductions by distribution for Poisson
manifolds, see Marsden and Ratiu \cite{mara86}, we note that the phase space
$T^*Q$ of the CH system is also a Poisson manifold, and its control
subset $W\subset T^*Q$ is a fiber submanifold. If we assume that
$D\subset TT^*Q |_W$ is a controllability distribution of the CH
system, then we can study naturally the Poisson
reduction by controllability distribution for the CH system. For a
symmetric CH system, and its control subset $W\subset T^*Q$ is a
$G$-invariant fiber submanifold, if we assume that $D\subset TT^*Q
|_W$ is a $G$-invariant controllability distribution of the symmetric CH
system, then we can give the Poisson reducible conditions
by controllability distribution for this CH system, and prove the Poisson reducible
property for the CH system and it is kept invariant under the CH-equivalence.
We also study the relationship between Poisson reduction by
$G$-invariant controllability distribution for the
regular (resp. singular) Poisson reducible CH system and
Poisson reduction by the reduced controllability distribution for the associated
reduced CH system. In addition, we can also develop the singular
Poisson reduction and SPR-CH-equivalence for the
CH system with symmetry, and prove the singular Poisson reduction theorem.
See Ratiu and Wang \cite{rawa12} for more details.

\section{Some Applications}

Now, it is a natural problem
if there is a practical RCH system and how to show the effect on the controls
in regular point reduction and Hamilton-Jacobi theory of the system.
In the following we shall give two examples of the application for the RCH system;
that is, the controlled rigid spacecraft -rotor system and
the controlled underwater vehicle-rotor system, and
give the Type I and Type II Hamilton-Jacobi equations
for the two systems.
See Wang \cite{wa20a, wa13e} for more details.

\subsection{The Controlled Rigid Spacecraft -Rotor System }

We first describe a rigid spacecraft carrying an internal "non-mass" rotor,
which is called a carrier body, where "non-mass" means that the mass of a rotor
is very very small comparing with the mass of the rigid spacecraft.
We first assume that the external forces
and the torques acting on the rigid spacecraft-rotor system are due to
buoyancy and gravity, and that the rigid spacecraft-rotor system with non-coincident
centers of buoyancy and gravity. We put a rotor within the rigid spacecraft,
and assume that the rotor
spins under the influence of a control torque $u$ acting on the rotor;
see Marsden \cite{ma92} and Leonard and Marsden \cite{lema97}.
In this case, the configuration space is
$Q=\textmd{SO}(3)\circledS \mathbb{R}^3\times S^1 \cong
\textmd{SE}(3)\times S^1$, with the first factor being the attitude
of the rigid spacecraft and the drift of the rigid spacecraft in the rotational
process and the second factor being the angle of rotor.
For convenience, using the local left trivialization, we denote uniformly that, locally,
$Q= \textmd{SE}(3)\times \mathbb{R}, $ and
$T^* Q= T^*(\textmd{SE}(3)\times \mathbb{R})\cong \textmd{SE}(3)\times \mathfrak{se}^\ast(3)
\times \mathbb{R} \times \mathbb{R}^{*}$, and
the Hamiltonian $H(A,c,\Pi,\Gamma,\alpha,l): T^* Q \cong
\textmd{SE}(3)\times \mathfrak{se}^\ast
(3)\times\mathbb{R}\times\mathbb{R}^* \to \mathbb{R}$ is given by
\begin{align}
 H(A,c,\Pi,\Gamma,\alpha,l) &=\Omega\cdot \Pi+\dot{\alpha}\cdot
l-L(A,c,\Omega,\Gamma,\alpha,\dot{\alpha}) \nonumber \\
&=\frac{1}{2}[\frac{\Pi_1^2}{\bar{I}_1}+\frac{\Pi_2^2}{\bar{I}_2}
+\frac{(\Pi_3-l)^2}{\bar{I}_3}+\frac{l^2}{J_3}]+ gh\Gamma \cdot
\chi \; .
\end{align}
In the following as the application of the theoretical result,
we shall regard the rigid spacecraft-rotor system
with the control torque $u$ acting on the rotor and
with the non-coincident centers of the buoyancy and the gravity
as a regular point reducible RCH
system on the generalization of the Euclidean group
$Q=\textmd{SE}(3)\times \mathbb{R}$,
in this case, we can give the regular point reduction
and the two types of Hamilton-Jacobi equations
for the controlled rigid
spacecraft-rotor system with the non-coincident centers of
the buoyancy and the gravity.\\

Assume that Lie group $G=\textmd{SE}(3)$ acts freely and properly on
$Q=\textmd{SE}(3)\times \mathbb{R}$ by the left translation on the first
factor $\textmd{SE}(3)$ and the trivial action on the second factor
$\mathbb{R}$, and that the action of
$\textmd{SE}(3)$ on phase space $T^\ast Q= T^\ast
\textmd{SE}(3)\times T^\ast \mathbb{R}$ is by the cotangent lift of the left
translation on $\textmd{SE}(3)$ at the identity,
and that the orbit space $(T^* Q)/ \textmd{SE}(3)$
is a smooth manifold and $\pi: T^*Q \rightarrow (T^*Q )/
\textmd{SE}(3) $ is a smooth submersion. Since $\textmd{SE}(3)$ acts
trivially on $\mathfrak{se}^\ast(3)$ and $\mathbb{R}\times \mathbb{R}^*$, it follows
that $(T^\ast Q)/ \textmd{SE}(3)$ is diffeomorphic to
$\mathfrak{se}^\ast(3) \times \mathbb{R}
\times \mathbb{R}^*$. For $(\mu,a) \in \mathfrak{se}^\ast(3)$, the
co-adjoint orbit $\mathcal{O}_{(\mu,a)} \subset
\mathfrak{se}^\ast(3)$ has the induced orbit symplectic form
$\omega^{-}_{\mathcal{O}_{(\mu,a)}}$, which coincides with the
restriction of the heavy top Lie-Poisson bracket on $\mathfrak{se}^\ast(3)$ to
the co-adjoint orbit $\mathcal{O}_{(\mu,a)}$. From Abraham and
Marsden \cite{abma78}, we know that $((T^\ast
\textmd{SE}(3))_{(\mu,a)},\omega_{(\mu,a)})$ is symplectically
diffeomorphic to
$(\mathcal{O}_{(\mu,a)},\omega_{\mathcal{O}_{(\mu,a)}}^{-})$, and
hence we have that the reduced space $((T^\ast Q)_{(\mu,a)},\omega_{(\mu,a)})$ is
symplectically diffeomorphic to $(\mathcal{O}_{(\mu,a)} \times
\mathbb{R}\times \mathbb{R}^*,\tilde{\omega}_{\mathcal{O}_{(\mu,a)}
\times \mathbb{R} \times \mathbb{R}^*}^{-}). $\\

From the expression $(8.1)$ of the Hamiltonian, we know that
$H(A,c,\Pi,\Gamma,\alpha,l)$ is invariant under the cotangent lift of the left
$\textmd{SE}(3)$-action $\Phi^{T*}: \textmd{SE}(3)\times T^\ast Q \to
T^\ast Q$. Assume that the control torque $u: T^\ast
Q \to W $ is invariant under the cotangent lift
$\Phi^{T*}$ of the left $\textmd{SE}(3)$-action, and
the dynamical vector field of the regular point reducible
controlled spacecraft-rotor system $(T^\ast Q,\textmd{SE}(3),\omega_Q,H,u)$
can be expressed by
\begin{align}
\tilde{X}= X_{(T^\ast Q,\textmd{SE}(3),\omega_Q,H,u)}= X_H+ \textnormal{vlift}(u),
\label{8.2} \end{align}
where $\textnormal{vlift}(u)= \textnormal{vlift}(u)\cdot X_H $
is the change of $X_H$ under the action of the control torque $u$.
For the point $(\mu,a)\in
\mathfrak{se}^\ast(3)$ is the regular value of the momentum map $\mathbf{J}_Q$, we
have the $R_p$-reduced Hamiltonian
$h_{(\mu,a)}(\Pi,\Gamma, \alpha,l): \mathcal{O}_{(\mu,a)}
\times\mathbb{R}\times\mathbb{R}^*(\subset \mathfrak{se}^\ast
(3)\times\mathbb{R}\times\mathbb{R}^*)\to \mathbb{R}$ given by
$$h_{(\mu,a)}(\Pi,\Gamma,\alpha,l)\cdot
\pi_{(\mu,a)}= H(A,c,\Pi,\Gamma,\alpha,l)|_{\mathcal{O}_{(\mu,a)}
\times\mathbb{R}\times\mathbb{R}^*},$$
and the  $R_p$-reduced control torque $u_{(\mu,a)}:
\mathcal{O}_{(\mu,a)} \times\mathbb{R}\times\mathbb{R}^{*} \to W_{(\mu,a)}
(\subset \mathcal{O}_{(\mu,a)} \times\mathbb{R}\times\mathbb{R}^{*}) $ is
given by
$$u_{(\mu,a)}(\Pi,\Gamma,\alpha,l)\cdot \pi_{(\mu,a)}=
u(A,c,\Pi,\Gamma, \alpha,l )|_{\mathcal{O}_{(\mu,a)}
\times\mathbb{R}\times\mathbb{R}^{*} }. $$
The $R_p$-reduced controlled spacecraft-rotor
system is the 4-tuple $(\mathcal{O}_{(\mu,a)} \times \mathbb{R}\times
\mathbb{R}^{*},\tilde{\omega}_{\mathcal{O}_{(\mu,a)} \times \mathbb{R}
\times \mathbb{R}^{*}}^{-},h_{(\mu,a)}, \\ u_{(\mu,a)}). $ \\

For convenience, the maps involved in
the following theorem are shown in Diagram-11.
\begin{center}
\hskip 0cm \xymatrix{ \mathbf{J}_Q^{-1}(\mu,a) \ar[r]^{i_{(\mu,a)}} & T^* Q
\ar[d]_{X_{H\cdot \varepsilon}} \ar[dr]^{\tilde{X}^\varepsilon} \ar[r]^{\pi_Q}
& Q \ar[d]^{\tilde{X}^\gamma} \ar[r]^{\gamma}
& T^*Q \ar[d]_{\tilde{X}} \ar[dr]_{X_{h_{(\mu,a)} \cdot\bar{\varepsilon}}} \ar[r]^{\pi_{(\mu,a)}}
& \;\;\; \mathcal{O}_{(\mu,a)}\times \mathbb{R}\times \mathbb{R}^{*} \ar[d]^{X_{h_{(\mu,a)}}} \\
& T(T^*Q)  & TQ \ar[l]^{T\gamma}
& T(T^*Q) \ar[l]^{T\pi_Q} \ar[r]_{T\pi_{(\mu,a)}}
& \;\;\; T(\mathcal{O}_{(\mu,a)}\times \mathbb{R}\times \mathbb{R}^{*})}
\end{center}
$$\mbox{Diagram-11}$$

\begin{theo}
In the case of non-coincident centers of buoyancy and gravity,
if the 5-tuple $(T^\ast Q, \\ \textmd{SE}(3),\omega_Q,H,u), $ where $Q=
\textmd{SE}(3)\times \mathbb{R}, $ is a regular point reducible
rigid spacecraft-rotor system with the control torque $u$ acting on the rotor,
then for a point $(\mu,a) \in \mathfrak{se}^\ast(3)$, the regular
value of the momentum map $\mathbf{J}_Q: \textmd{SE}(3)\times
\mathfrak{se}^\ast(3) \times \mathbb{R}\times \mathbb{R}^{*} \to
\mathfrak{se}^\ast(3)$, the $R_p$-reduced controlled rigid spacecraft-rotor system is the 4-tuple
$(\mathcal{O}_{(\mu,a)} \times \mathbb{R}\times
\mathbb{R}^{*},\tilde{\omega}_{\mathcal{O}_{(\mu,a)} \times \mathbb{R}
\times \mathbb{R}^{*}}^{-},h_{(\mu,a)},u_{(\mu,a)}). $ Assume that
$\gamma: \textmd{SE}(3)\times \mathbb{R}\rightarrow
T^*(\textmd{SE}(3)\times \mathbb{R})$ is a one-form on
$\textmd{SE}(3)\times \mathbb{R}$,
and that $\lambda=\gamma \cdot \pi_{(\textmd{SE}(3)\times \mathbb{R})}:
T^* (\textmd{SE}(3)\times \mathbb{R} )\rightarrow T^* (\textmd{SE}(3)\times \mathbb{R}), $
and that $\varepsilon:
T^* (\textmd{SE}(3) \times \mathbb{R})\rightarrow T^* (\textmd{SE}(3)\times \mathbb{R}) $ is an
$\textmd{SE}(3)_{(\mu,a)}$-invariant symplectic map.
Denote that
$\tilde{X}^\gamma = T\pi_{(\textmd{SE}(3)\times \mathbb{R})}\cdot \tilde{X}\cdot \gamma$, and that
$\tilde{X}^\varepsilon = T\pi_{(\textmd{SE}(3)\times \mathbb{R})}\cdot \tilde{X}\cdot \varepsilon$,
where $\tilde{X}=X_{(T^\ast Q,\textmd{SE}(3),\omega_Q,H,u)}$ is the dynamical vector field of
the controlled rigid spacecraft-rotor system $(T^\ast Q,\textmd{SE}(3),\omega_Q,H,u)$.
Moreover, assume that $\textmd{Im}(\gamma)\subset \mathbf{J}_Q^{-1}(\mu,a), $ and that $\gamma$ is
$\textmd{SE}(3)_{(\mu,a)}$-invariant,
and that $\varepsilon(\mathbf{J}_Q^{-1}(\mu,a))\subset \mathbf{J}_Q^{-1}(\mu,a). $
Denote that $\bar{\gamma}=\pi_{(\mu,a)}(\gamma):
\textmd{SE}(3)\times \mathbb{R} \rightarrow
\mathcal{O}_{(\mu,a)}\times \mathbb{R}\times \mathbb{R}^{*}, $ and that
$\bar{\lambda}=\pi_{(\mu,a)}(\lambda): T^* (\textmd{SE}(3)\times \mathbb{R}) \rightarrow
\mathcal{O}_{(\mu,a)}\times \mathbb{R}\times \mathbb{R}^{*}, $ and that
$\bar{\varepsilon}=\pi_{(\mu,a)}(\varepsilon): \mathbf{J}_Q^{-1}(\mu,a)\rightarrow
\mathcal{O}_{(\mu,a)}\times \mathbb{R}\times \mathbb{R}^{*}. $
Then the following two assertions hold:\\
\noindent $(\mathbf{i})$
If the one-form $\gamma: \textmd{SE}(3)\times \mathbb{R} \rightarrow
T^*(\textmd{SE}(3)\times \mathbb{R}) $ is closed with respect to
$T\pi_{(\textmd{SE}(3)\times \mathbb{R})}: TT^* (\textmd{SE}(3)\times \mathbb{R})
\rightarrow T(\textmd{SE}(3)\times \mathbb{R}), $
then $\bar{\gamma}$ is a solution of the Type I Hamilton-Jacobi equation
$T\bar{\gamma}\cdot \tilde{X}^\gamma= X_{h_{(\mu,a)}}\cdot \bar{\gamma}. $\\
\noindent $(\mathbf{ii})$
The $\varepsilon$ and $\bar{\varepsilon} $ satisfy the Type II Hamilton-Jacobi equation
$T\bar{\gamma}\cdot \tilde{X}^\varepsilon= X_{h_{(\mu,a)}}\cdot \bar{\varepsilon} $
if and only if they satisfy
the equation $T\bar{\varepsilon}\cdot(X_{h_{(\mu,a)} \cdot \bar{\varepsilon}})
= T\bar{\lambda}\cdot \tilde{X}\cdot\varepsilon. $
\end{theo}
See Wang \cite{wa20a} for more details.

\subsection{The Controlled Underwater Vehicle-Rotor System}

We first describe a underwater vehicle carrying two internal "non-mass" rotors;
that is called a carrier body, moving in a given fluid,
where "non-mass" means that the mass of a rotor
is very very small compared to the mass of the underwater vehicle.
We consider that the underwater vehicle-rotor system has a
neutrally buoyant, and that rigid body (often ellipsoidal) submerged in an
infinitely large volume of the incompressible, the inviscid, the irrotational
fluid which is at rest at infinity. The dynamics of the carrier body-fluid
system are described by the Kirchhoff's equations.
We first assume that the external forces
and the torques acting on the underwater vehicle-rotor system are due to
the buoyancy and the gravity, and that the underwater vehicle-rotor system with non-coincident centers
of the buoyancy and the gravity. We put two rotors within the underwater vehicle,
and assume that the rotors spin under the influence of a control torque $u$
acting on the rotors; see Leonard and Marsden \cite{lema97}.
In this case, the configuration space is $Q=W \times V$, where
$W= \textmd{SE}(3)\circledS \mathbb{R}^3= (\textmd{SO}(3)\circledS
\mathbb{R}^3)\circledS \mathbb{R}^3$ is a double semidirect product Lie group
and $V=S^1\times S^1$, with the first factor being the attitude and the
position of the underwater vehicle as well as the drift of the underwater
vehicle-rotor in the rotational and the translational process, and the second
factor being the angles of the rotors. For convenience,
using the local left trivialization, we denote uniformly that, locally,
$Q= \textmd{SE}(3)\circledS \mathbb{R}^3 \times \mathbb{R}^2, $ and
$T^* Q= T^*(\textmd{SE}(3)\circledS \mathbb{R}^3 \times \mathbb{R}^2)
\cong \textmd{SE}(3)\circledS \mathbb{R}^3\times \mathfrak{se}^\ast(3)\circledS \mathbb{R}^{3*}
\times \mathbb{R}^2 \times \mathbb{R}^{2*}$, and
the Hamiltonian $H(A,c,b,\Pi,\Gamma,P,\theta,l):T^* Q\cong
(\textmd{SE}(3)\circledS \mathbb{R}^3)\times (\mathfrak{se}^*(3)\circledS \mathbb{R}^{3*})
\times\mathbb{R}^2\times\mathbb{R}^{2*}\to \mathbb{R}$
is given by
\begin{align} &H(A,c,b,\Pi,\Gamma,P,\theta,l)=\Omega\cdot
\Pi+ v\cdot P+ \dot{\theta}\cdot l-L(A,c,b,\Omega,\Gamma,v,\theta,\dot{\theta}) \nonumber \\
&=\frac{1}{2}[\frac{(\Pi_1-l_1)^2}{\bar{I}_1}+\frac{(\Pi_2-l_2)^2}{\bar{I}_2}
+\frac{\Pi_3^2}{\bar{I}_3}+
\frac{P_1^2}{m_1}+\frac{P_2^2}{m_2}+\frac{P_3^2}{m_3}+\frac{l_1^2}{J_1}+\frac{l_2^2}{J_2}]+
gh\Gamma \cdot \chi \;.
 \label{8.3}
\end{align}
In the following as the application of the theoretical result,
we shall regard the controlled underwater vehicle-rotor system
with the control torque $u$ acting on the rotors and
with the non-coincident centers of the buoyancy and the gravity
as a regular point reducible RCH
system on the generalization of the semidirect product Lie group
$Q=W\times V$, in this case, we can give the regular point reduction
and the two types of Hamilton-Jacobi equations for the controlled underwater
vehicle-rotor system with the non-coincident centers of the buoyancy and the gravity.\\

Assume that semidirect product Lie group $W=\textmd{SE}(3)\circledS
\mathbb{R}^3 $ acts freely and properly on $Q=
(\textmd{SE}(3)\circledS \mathbb{R}^3 )\times \mathbb{R}^2 $ by the left
translation on $\textmd{SE}(3)\circledS \mathbb{R}^3 $, and that the
action of $W$ on the phase space $T^\ast Q$ is by cotangent lift of
left translation on $Q$ at the identity. If $(\Pi,w_1,w_2) \in
\mathfrak{se}^\ast(3)\circledS \mathbb{R}^{3*}$ is a regular value of
the momentum map $\mathbf{J}_Q$, then the $R_p$-reduced space
$$(T^\ast Q)_{(\Pi,w_1,w_2)}=
\mathbf{J}^{-1}_Q(\Pi,w_1,w_2)/(\textmd{SE}(3)\circledS \mathbb{R}^3
)_{(\Pi,w_1,w_2)}$$
is symplectically diffeomorphic to the
orbit space $\mathcal{O}_{(\Pi,w_1,w_2)} \times \mathbb{R}^2 \times
\mathbb{R}^{2*} \subset \mathfrak{se}^\ast(3)\circledS \mathbb{R}^{3*}
\times \mathbb{R}^2 \times \mathbb{R}^{2*} $ with the induced orbit
symplectic form $\tilde{\omega}_{\mathcal{O}_{(\mu,a_1,a_2)} \times \mathbb{R}^2
\times \mathbb{R}^{2*}} ^{-}.$ \\

From the expression $(8.3)$ of the Hamiltonian, we know that
$H(A,c,b,\Pi,\Gamma, P, \theta,l)$ is invariant
under the cotangent lift of the left
$\textmd{SE}(3)\circledS \mathbb{R}^3$-action
$\Phi^{T^*}:\textmd{SE}(3)\circledS \mathbb{R}^3\times T^\ast Q\to
T^\ast Q$. Assume that the control torque
$u: T^\ast Q \to W$ acting on the rotors is invariant under the cotangent lift
of the left $\textmd{SE}(3)\circledS \mathbb{R}^3$-action, and that
the dynamical vector field of the regular point reducible
controlled underwater vehicle-rotor system
$(T^\ast Q,\textmd{SE}(3)\circledS \mathbb{R}^3,\omega_Q,H,u)$
can be expressed by
\begin{align}
\tilde{X}= X_{(T^\ast Q,\textmd{SE}(3)\circledS \mathbb{R}^3,\omega_Q,H,u)}= X_H+ \textnormal{vlift}(u),
\label{5.3}
\end{align}
where $\textnormal{vlift}(u)= \textnormal{vlift}(u)\cdot X_H $
is the change of $X_H$ under the action of the control torque $u$.
For the point $(\mu,a_1,a_2)\in \mathfrak{se}^\ast(3)\circledS \mathbb{R}^{3*}$,
the regular value of the momentum map $\mathbf{J}_Q$, we
have the $R_p$-reduced Hamiltonian
$h_{(\mu,a_1,a_2)}(\Pi,\Gamma,P,\theta,l):\mathcal{O}_{(\mu,a_1,a_2)}\times\mathbb{R}^2
\times \mathbb{R}^{2*} (\subset \mathfrak{se}^\ast (3)\circledS \mathbb{R}^{3*}\times
\mathbb{R}^2\times \mathbb{R}^{2*}) \to \mathbb{R}$ given by
$$h_{(\mu,a_1,a_2)}(\Pi,\Gamma,P,\theta,l)\cdot \pi_{(\mu,a_1,a_2)}
=H(A,c,b,\Pi,\Gamma,P,\theta,l)|_{\mathcal{O}_{(\mu,a_1,a_2)}\times
\mathbb{R}^2\times \mathbb{R}^{2*}}, $$
and that the  $R_p$-reduced control torque $u_{(\mu,a_1,a_2)}:
\mathcal{O}_{(\mu,a_1,a_2)} \times\mathbb{R}^2\times\mathbb{R}^{2*} \to \mathcal{C}_{(\mu,a_1,a_2)}
(\subset \mathcal{O}_{(\mu,a_1,a_2)} \times\mathbb{R}^2\times\mathbb{R}^{2*}) $ is
given by $$u_{(\mu,a_1,a_2)}(\Pi,\Gamma,P,\theta,l)\cdot \pi_{(\mu,a_1,a_2)}=
u(A,c,b,\Pi,\Gamma, P,\theta,l )|_{\mathcal{O}_{(\mu,a_1,a_2)}
\times\mathbb{R}^2\times\mathbb{R}^{2*} }. $$
The $R_p$-reduced controlled underwater vehicle-rotor
system is the 4-tuple $(\mathcal{O}_{(\mu,a_1,a_2)} \times \mathbb{R}^2 \times
\mathbb{R}^{2*}, \\ \tilde{\omega}_{\mathcal{O}_{(\mu,a_1,a_2)} \times \mathbb{R}^2
\times \mathbb{R}^{2*}}^{-},h_{(\mu,a_1,a_2)},u_{(\mu,a_1,a_2)}). $ \\

For convenience, the maps involved in
the following theorem are shown in Diagram-12.
\begin{center}
\hskip 0cm \xymatrix{ \mathbf{J}_Q^{-1}(\mu,a_1,a_2) \ar[r]^{i_{(\mu,a_1,a_2)}} & T^* Q
\ar[d]_{X_{H\cdot \varepsilon}} \ar[dr]^{\tilde{X}^\varepsilon} \ar[r]^{\pi_Q}
& Q \ar[d]^{\tilde{X}^\gamma} \ar[r]^{\gamma}
& T^*Q \ar[d]_{\tilde{X}} \ar[dr]^{X_{h_{(\mu,a_1,a_2)} \cdot\bar{\varepsilon}}} \ar[r]^{\pi_{(\mu,a_1,a_2)}}
& \;\;\;\;\;\; \mathcal{O}_{(\mu,a_1,a_2)}\times \mathbb{R}^2\times \mathbb{R}^{2*} \ar[d]^{X_{h_{(\mu,a_1,a_2)}}} \\
& T(T^*Q)  & TQ \ar[l]^{T\gamma}
& T(T^*Q) \ar[l]^{T\pi_Q} \ar[r]_{T\pi_{(\mu,a_1,a_2)}}
& \;\;\;\;\;\; T(\mathcal{O}_{(\mu,a_1,a_2)}\times \mathbb{R}^2\times \mathbb{R}^{2*})}
\end{center}
$$\mbox{Diagram-12}$$

\begin{theo}
In the case of non-coincident centers of the buoyancy and the gravity,
if the 5-tuple $(T^\ast Q,\textmd{SE}(3)\circledS \mathbb{R}^3,\omega_Q,H,u), $ where $Q=
\textmd{SE}(3)\circledS \mathbb{R}^3\times \mathbb{R}^2, $ is a regular point reducible
underwater vehicle-rotor system with the control torque $u$ acting on the rotors,
then, for a point $(\mu,a_1,a_2)\in \mathfrak{se}^\ast(3)\circledS
\mathbb{R}^{3*}$, the regular value of the momentum map $\mathbf{J}_Q:
T^* Q \cong \textmd{SE}(3)\circledS \mathbb{R}^3 \times
\mathfrak{se}^\ast(3)\circledS \mathbb{R}^{3*} \times
\mathbb{R}^2\times \mathbb{R}^{2*} \to \mathfrak{se}^\ast(3)\circledS
\mathbb{R}^{3*}$, the $R_p$-reduced controlled underwater vehicle-rotor system is the 4-tuple
$(\mathcal{O}_{(\mu,a_1,a_2)} \times \mathbb{R}^2 \times \mathbb{R}^{2*},
\tilde{\omega}^{-}_{\mathcal{O}_{(\mu,a_1,a_2)} \times
\mathbb{R}^2 \times \mathbb{R}^{2*}},h_{(\mu,a_1,a_2)},u_{(\mu,a_1,a_2)})$,
where $ \mathcal{O}_{(\mu,a_1,a_2)}
\subset \mathfrak{se}^\ast(3)\circledS \mathbb{R}^{3*}$ is the co-adjoint orbit
of the semidirect product Lie group $\textmd{SE}(3)\circledS \mathbb{R}^3$.
Assume that
$\gamma: Q \rightarrow T^* Q $ is a one-form on $Q=\textmd{SE}(3)\circledS
\mathbb{R}^3 \times \mathbb{R}^2 $, and that $\lambda=\gamma \cdot \pi_Q:
T^* Q \rightarrow T^* Q, $ and that $\varepsilon: T^* Q \rightarrow T^* Q $ is an
$\textmd{(SE}(3)\circledS \mathbb{R}^3)_{(\mu,a_1,a_2)}$-invariant symplectic map,
where $(\textmd{SE}(3)\circledS \mathbb{R}^3)_{(\mu,a_1,a_2)}$ is
the isotropy subgroup of the co-adjoint $\textmd{SE}(3)\circledS
\mathbb{R}^3$-action at the point $(\mu,a_1,a_2)$.
Denote that
$\tilde{X}^\gamma = T\pi_Q \cdot \tilde{X}\cdot \gamma$, and that
$\tilde{X}^\varepsilon = T\pi_Q \cdot \tilde{X}\cdot \varepsilon$,
where $\tilde{X}=X_{(T^\ast Q,\textmd{SE}(3)\circledS
\mathbb{R}^3,\omega_Q,H,u)}$ is the dynamical vector field of
the controlled underwater vehicle-rotor system $(T^\ast Q,\textmd{SE}(3)\circledS
\mathbb{R}^3,\omega_Q,H,u)$.
Moreover, assume that
$\textmd{Im}(\gamma)\subset
\mathbf{J}_Q^{-1}((\mu,a_1,a_2)), $ and that $\gamma$ is
$(\textmd{SE}(3)\circledS \mathbb{R}^3)_{(\mu,a_1,a_2)}$-invariant,
and that $\varepsilon (\mathbf{J}_Q^{-1}(\mu,a_1,a_2))\subset \mathbf{J}_Q^{-1}(\mu,a_1,a_2). $
 Denote that
$\bar{\gamma}=\pi_{(\mu,a_1,a_2)}(\gamma): Q=\textmd{SE}(3)\circledS
\mathbb{R}^3 \times \mathbb{R}^2 \rightarrow \tilde{\mathcal{O}}, $
and that $\bar{\lambda}=\pi_{(\mu,a_1,a_2)}(\lambda):
T^* Q \rightarrow \tilde{\mathcal{O}} $,
 $\bar{\varepsilon}=\pi_{(\mu,a_1,a_2)}(\varepsilon): \mathbf{J}_Q^{-1}(\mu,a_1,a_2)\rightarrow
\tilde{\mathcal{O}}. $
Then the following two assertions hold:\\
\noindent $(\mathbf{i})$
If the one-form $\gamma: Q \rightarrow T^* Q $ is closed with respect to
$T\pi_Q: TT^* Q \rightarrow TQ, $
then $\bar{\gamma}$ is a solution of the Type I Hamilton-Jacobi equation
$T\bar{\gamma}\cdot \tilde{X}^\gamma= X_{h_{(\mu,a_1,a_2)}}\cdot \bar{\gamma}. $\\
\noindent $(\mathbf{ii})$
The $\varepsilon$ and the $\bar{\varepsilon} $ satisfy the Type II Hamilton-Jacobi equation
$T\bar{\gamma}\cdot \tilde{X}^\varepsilon= X_{h_{(\mu,a_1,a_2)}}\cdot \bar{\varepsilon} $
if and only if they satisfy
the equation $T\bar{\varepsilon}\cdot(X_{h_{(\mu,a_1,a_2)} \cdot \bar{\varepsilon}})
= T\bar{\lambda}\cdot \tilde{X}\cdot\varepsilon. $
\end{theo}
See Wang \cite{wa13e} for more details.\\

It is worthy of noting that
the motions of the controlled rigid spacecraft-rotor system
and the controlled underwater vehicle-rotor system are different, and
the configuration spaces, the Hamiltonian functions, the actions of Lie groups,
the $R_p$-reduced symplectic forms and the $R_p$-reduced systems of
the controlled rigid spacecraft-rotor system
and the controlled underwater vehicle-rotor system,
all of them are also different. However,
the two types of Hamilton-Jacobi equations given by calculation in detail
are same, that is, the internal rules are same by comparing
Theorems 8.1 and 8.2. This is very important!
In particular, the method of calculations for the regular point
reductions and the two types of Hamilton-Jacobi equations for the rigid spacecraft with
internal rotor and the underwater vehicle with internal rotors is very important and efficient,
and it is generalized and used to the more practical controlled Hamiltonian systems.
Finally, we also note that there have been a
lot of beautiful results of reduction theory of Hamiltonian systems
in celestial mechanics, hydrodynamics and plasma physics. So, it is
an important topic to study the application of symmetric reduction and
Hamilton-Jacobi theory for the regular controlled Hamiltonian systems
in celestial mechanics, hydrodynamics and
plasma physics. These are our goals in future research.\\

It is well-known that the symmetric reduction and
Hamilton-Jacobi theory for Hamiltonian systems
are very important research subjects
in mathematics and analytical mechanics.
Following the theoretical development of geometric mechanics, a lot
of important problems involving the subjects are being explored and
studied; see Wang and Zhang \cite{wazh12},
Ratiu and Wang \cite{rawa12}, and Wang \cite{wa18, wa15a, wa21a}.
In addation, it is very important for a rigorous theoretical work to
offer uniformly a composition of the research results
from a global view point.
In particular, it is the key thought of the researches of geometrical mechanics
of the Professor Jerrold E. Marsden to explore and reveal the deeply internal
relationship between the geometrical structure of phase space and the dynamical
vector field of a mechanical system. It is also our goal of pursuing and inheriting.\\

\noindent {\bf Acknowledgments:}
Jerrold E. Marsden (1942-08-17---2010-09-21) is a great mathematician,
mechanician and physician. He is a model of the excellent scientists.
The author would like to thank Professor Jerrold E. Marsden
for his understanding, support and help
in the study of geometric mechanics, and to
dedicate the article to the 80th anniversary of the birth
of Professor Jerrold E. Marsden.\\

\end{document}